\documentclass[10pt,twocolumn,twoside]{IEEEtran}

\usepackage{graphicx}%
\usepackage{multirow}%
\usepackage{amsmath,amssymb,amsfonts}%
\usepackage{amsthm}%
\usepackage{mathrsfs}%
\usepackage{xcolor}%
\usepackage{textcomp}%
\usepackage{manyfoot}%
\usepackage{booktabs}%
\usepackage{algorithm}%
\usepackage{algorithmicx}%
\usepackage{algpseudocode}%

\usepackage{threeparttable}

\usepackage{epstopdf}%

\usepackage{here}
\usepackage{bm}
\usepackage{algorithm}
\usepackage{algpseudocode}
\usepackage{comment}
\usepackage{cases}
\usepackage{mathtools}
\usepackage{enumitem}
\usepackage{mathabx}
\usepackage{hyperref}
\usepackage{hyperref}
\hypersetup{
    hidelinks=true,
}
\usepackage{stmaryrd}
\usepackage{tikz}
\usetikzlibrary{arrows.meta}
\usetikzlibrary{positioning}
\usetikzlibrary{decorations.pathreplacing,calc}
\usepackage{csvsimple}
\usepackage{lscape}

\usepackage{etoolbox}

\usepackage{tabularx}
\usepackage{booktabs}

\newcommand{\argmin}{\mathop{\mathrm{argmin}}}

\newcommand{\inprod}[2]{{\left\langle #1,#2 \right\rangle}}

\newcommand{\TT}{\mathsf{T}}

\newcommand{\prox}[1]{\mathrm{prox}_{#1}}
\newcommand{\moreau}[2]{{}^{#2}#1}

\newcommand{\exR}{\mathbb{R}\cup \{+\infty\}}
\newcommand{\dom}[1]{\mathrm{dom}(#1)}

\newcommand{\Limsup}{\mathop{\mathrm{Lim~sup}}\limits}

\newcommand{\norm}[1]{\left\lVert #1 \right\rVert}
\newcommand{\Id}{\mathrm{Id}}

\newcommand{\Fsubdiff}{\partial_{\rm F}}
\newcommand{\Lsubdiff}{\partial_{\rm L}}
\newcommand{\subdiff}{\partial}
\newcommand{\dist}{\mathop{\mathrm{dist}}}

\makeatletter
\newcommand{\doublewidetilde}[1]{{%
  \mathpalette\double@widetilde{#1}%
}}
\newcommand{\double@widetilde}[2]{%
  \sbox\z@{$\m@th#1\widetilde{#2}$}%
  \ht\z@=.9\ht\z@
  \widetilde{\box\z@}%
}
\makeatother

\theoremstyle{plain}%
\newtheorem{theorem}{Theorem}[section]
\newtheorem{lemma}[theorem]{Lemma}

\newtheorem{proposition}[theorem]{Proposition}

\theoremstyle{definition}
\newtheorem{definition}[theorem]{Definition}
\newtheorem{fact}[theorem]{Fact}
\newtheorem{example}[theorem]{Example}
\newtheorem{assumption}[theorem]{Assumption}
\newtheorem{problem}[theorem]{Problem}
\newtheorem{remark}[theorem]{Remark}

\theoremstyle{remark}

\DeclarePairedDelimiter{\abs}{\lvert}{\rvert}

\mathtoolsset{showonlyrefs=true}

\newcommand{\param}{\bm{\varphi}}
\newcommand{\cost}{F}
\newcommand{\opnorm}[1]{\norm{#1}_{2}}
\newcommand{\stdnorm}[1]{\norm{#1}_{2}}

\hypersetup{
    hidelinks=true,
}
\usepackage{amsmath,amsfonts}
\usepackage{array}
\usepackage[caption=false,font=normalsize,labelfont=sf,textfont=sf]{subfig}
\usepackage{textcomp}
\usepackage{stfloats}
\usepackage{url}
\usepackage{verbatim}
\usepackage{graphicx}
\usepackage{cite}
\begin{document}

\mathtoolsset{showonlyrefs=true}
\title{A Proximal Variable Smoothing for Minimization of Nonlinearly Composite Nonsmooth Function   \\ -- Finite-Max Minimization and MIMO Applications}

\author{Keita~{\sc Kume},~\IEEEmembership{Member,~IEEE,}~
Isao~{\sc Yamada},~\IEEEmembership{Fellow,~IEEE}
\thanks{
  K.~Kume and I.~Yamada are with
  the Dept. of Information and Communications Engineering,
  Institute of Science Tokyo,
  2-12-1-S3-60, O-okayama,
  Meguro-ku, Tokyo 152-8550, Japan
  (e-mail: kume@sp.ict.e.titech.ac.jp).
  This work was partially supported by JSPS Grants-in-Aid (19H04134, 24K23885).
  Compared to a conference paper~\cite{Kume-Yamada25A}, this paper contains new contributions:
  (i) complete proofs of mathematical results in~\cite{Kume-Yamada25A}; 
  (ii) convergence rate analysis of the proposed algorithm;
  (iii) applicability of the proposed algorithm to finite-max nonlinear minimization and MIMO applications.
}
}

\maketitle

\begin{abstract}
  We propose a proximal variable smoothing algorithm for a nonsmooth optimization problem whose cost function is the sum of three functions including a weakly convex composite function.
  The proposed algorithm has a single-loop structure inspired by a proximal gradient-type method.
  More precisely, the proposed algorithm consists of two steps: (i)  a gradient descent of a time-varying smoothed surrogate function
  designed partially with the Moreau envelope of the weakly convex function; (ii) an application of the proximity operator of the remaining function
  not covered by the smoothed surrogate function.
  For the proposed algorithm, we present a subsequential convergence guarantee in terms of a stationary point,  and a convergence rate
  $\mathcal{O}(\epsilon^{-3})$
  for achieving an $\epsilon$-stationary point.
  Numerical experiments demonstrate the effectiveness of the proposed algorithm in two scenarios:
  (i)
  robust target localization and
  (ii)
  multiple-input-multiple-output (MIMO) signal detection.
\end{abstract}

\begin{IEEEkeywords}
  Nonsmooth nonconvex optimization, Moreau envelope, variable smoothing, constrained finite-max nonlinear minimization,
 MIMO detection
\end{IEEEkeywords}

\section{Introduction}

Sparsity-aware or robust signal estimation problems have been tackled as typical inverse problems and have been fundamental tasks in signal processing and machine learning.
Such estimation problems include, e.g.,
regularized least squares estimation~\cite{Tibshirani96,Chen-Gu14,Selesnick17,Abe-Yamagishi-Yamada20,Yata-Yamagishi-Yamada22,Kuroda-Kitahara22},
sparse spectral clustering~\cite{Lu-Yan-Lin16,Wang-Liu-Chen-Ma-Xue-Zhao22,Kume-Yamada24A}, robust phase retrieval~\cite{Duchi-Ruan18,Zhen-Ma-Xue24,Yazawa-Kume-Yamada26}, and robust matrix recovery~\cite{Charisopoulos-Chen-Davis-Diaz-Ding-Drusvyatskiy21,Wang-So-Zoubir23}.
For such estimation problems, the most common approach is to formulate the task as the following nonsmooth (possibly nonconvex) optimization problem, where its minimizer (or, alternatively, its {\em stationary point} in Definition~\ref{definition:optimality}~\ref{enum:definition:optimality:stationary}) is assigned to an estimate of the target signal.
\begin{problem}\label{problem:origin}
Let
$\mathcal{X}$
and
$\mathcal{Z}$
be Euclidean spaces, i.e.,
finite-dimensional real Hilbert spaces.
Then,
\begin{equation}
  \mathrm{find} \ \bm{x}^{\star} \in \argmin_{\bm{x}\in \mathcal{X}} (h+g\circ \mathfrak{S} + \phi)(\bm{x}) (\eqqcolon (\cost+\phi)(\bm{x})),
\end{equation}
where
$\phi$,
$h$,
$\mathfrak{S}$,
$g$,
and
$\cost\coloneqq h+g\circ\mathfrak{S}:\mathcal{X}\to\mathbb{R}$
satisfy the following conditions:
\begin{enumerate}[label=(\roman*)]
  \item
        $\phi:\mathcal{X}\to \exR$
        is a proper lower semicontinuous convex function, i.e.,
        $\phi \in \Gamma_{0}(\mathcal{X})$,
        and~{\em prox-friendly}, i.e., the proximity operator
        $\prox{\gamma \phi}\ (\gamma \in \mathbb{R}_{++})$
        (see~\eqref{eq:prox_definition})
        is available as a computable tool;
  \item \label{enum:problem:origin:h}
        $h:\mathcal{X} \to \mathbb{R}$
        is continuously differentiable, and its gradient
        $\nabla h:\mathcal{X}\to\mathcal{X}$
        is $L_{\nabla h}(>0)$-Lipschitz continuous over
        $\dom{\phi}\coloneqq\{\bm{x} \in \mathcal{X} \mid \phi(\bm{x}) < +\infty\}$,
        i.e.,
        $\stdnorm{\nabla h(\bm{x}_{1}) - \nabla h(\bm{x}_{2})}\leq L_{\nabla h}\stdnorm{\bm{x}_{1} - \bm{x}_{2}}$
        holds for
        $\bm{x}_{1},\bm{x}_{2} \in \dom{\phi}$;
  \item \label{enum:problem:origin:S}
        $\mathfrak{S}: \mathcal{X} \to \mathcal{Z}$
        is a continuously
        differentiable (possibly nonlinear) mapping such that the derivative
        $\mathrm{D}\mathfrak{S}$
        is $L_{\mathrm{D}\mathfrak{S}}(>0)$-Lipschitz continuous over
        $\dom{\phi}$,
        i.e.,
        $\opnorm{\mathrm{D}\mathfrak{S}(\bm{x}_{1}) - \mathrm{D}\mathfrak{S}(\bm{x}_{2})} \leq L_{\mathrm{D}\mathfrak{S}}\stdnorm{\bm{x}_{1}-\bm{x}_{2}}$
        holds for
        $\bm{x}_{1},\bm{x}_{2} \in \dom{\phi}$;
  \item
        $g:\mathcal{Z} \to \mathbb{R}$
        is (a) $L_{g}(>0)$-Lipschitz continuous
        (possibly nonsmooth),
        (b)~$\eta$-weakly convex with
        $\eta > 0$,
        i.e.,
        $g+\frac{\eta}{2}\stdnorm{\cdot}^{2}$
        is convex,
        and (c) prox-friendly;
  \item \label{enum:problem:origin:minimizer}
        $\cost + \phi$
        is bounded below, i.e.,
        $\inf_{\bm{x}\in\mathcal{X}}(\cost+\phi)(\bm{x}) > -\infty$.
\end{enumerate}
\end{problem}

In typical applications of Problem~\ref{problem:origin},
$\phi \in \Gamma_{0}(\mathcal{X})$
is given as an {\em indicator function}
\begin{equation}
  \iota_{C}:\mathcal{X}\to \exR:\bm{x}\mapsto
  \left\{\begin{array}{cl} 0,  & \mathrm{if}\ \bm{x}\in C;    \\
             +\infty, & \mathrm{if}\ \bm{x}\notin C,
  \end{array}\right. \label{eq:indicator}
\end{equation}
associated with a nonempty closed convex set
$C \subset \mathcal{X}$,
which works as a convex constraint set.
The function
$h$
is given as a data-fidelity function, e.g.,
$\mathcal{X}\ni \bm{x}\mapsto\frac{1}{2}\stdnorm{\bm{y} - \bm{A}\bm{x}}^{2}$
for the linear regression model
$\bm{y} = \bm{A}\bm{x}^{\star} + \bm{\epsilon} \in \mathcal{Y}$
with a linear operator
$\bm{A}:\mathcal{X} \to \mathcal{Y}$
and noise
$\bm{\epsilon} \in \mathcal{Y}$
(see, e.g.,~\cite{Tibshirani96,Selesnick17,Abe-Yamagishi-Yamada20,Yata-Yamagishi-Yamada22,Kuroda-Kitahara22}),
where
$\mathcal{Y}$
is a Euclidean space.
The function
$g\circ\mathfrak{S}$
works as, e.g., a certain regularizer.
For example,
$g$
serves as a sparsity-promoting function in sparsity-aware applications~\cite{Chen-Gu14,Selesnick17,Abe-Yamagishi-Yamada20,Yata-Yamagishi-Yamada22,Kuroda-Kitahara22,Lu-Yan-Lin16,Wang-Liu-Chen-Ma-Xue-Zhao22,Kume-Yamada24A} or a robust loss function in robust estimation~\cite{Duchi-Ruan18,Zhen-Ma-Xue24,Yazawa-Kume-Yamada26,Charisopoulos-Chen-Davis-Diaz-Ding-Drusvyatskiy21,Wang-So-Zoubir23,Yang-Chen-Ma-Chen-Gu-So19,Suzuki-Yukawa23}, where such functions
$g$
include, e.g.,
$\ell_{1}$-norm, Minimax Concave Penalty (MCP)~\cite{Zhang10} and Smoothly Clipped Absolute Deviation (SCAD)~\cite{Fan-Li01}
(see also~\cite[Sect. 2.2]{Kume-Yamada24A} for other possible examples).
Moreover,
$\mathfrak{S}$
is designed according to a prior knowledge on the target signal
$\bm{x}^{\star}$
(see~\cite{Abe-Yamagishi-Yamada20,Yata-Yamagishi-Yamada22,Kuroda-Kitahara22} for the linear case of
$\mathfrak{S}$,
and~\cite{Lu-Yan-Lin16,Wang-Liu-Chen-Ma-Xue-Zhao22,Kume-Yamada24A,Duchi-Ruan18,Zhen-Ma-Xue24,Yazawa-Kume-Yamada26,Charisopoulos-Chen-Davis-Diaz-Ding-Drusvyatskiy21,Wang-So-Zoubir23} for the nonlinear case of
$\mathfrak{S}$).

Nevertheless, Problem~\ref{problem:origin} is a challenging nonconvex nonsmooth optimization problem due to the composition of the nonsmooth function
$g$
and
the nonlinear mapping
$\mathfrak{S}$.
As shown in Table~\ref{table:summary},
existing methods related to Problem~\ref{problem:origin} involve a trade-off among algorithmic simplicity, subsequential convergence guarantees, and nonasymptotic convergence rates  regarding
a ($\epsilon$-) stationary point of
$\cost + \phi$
in Problem~\ref{problem:origin}.

For Problem~\ref{problem:origin},
{\em proximal linear methods (PLM)}, e.g.,~\cite{Lewis-Wright16,Duchi-Ruan18,Drusvyatskiy-Paquette19,Zhen-Ma-Xue24}, are known to be the most prominent algorithms where a subproblem involving linearization of
$\mathfrak{S}$
in
$g\circ\mathfrak{S}$
is solved at each iteration.
For PLM, a subsequential convergence guarantee of a stationary point of
$\cost + \phi$
is established in~\cite{Lewis-Wright16},
and a convergence rate
$\mathcal{O}(\epsilon^{-3}\log(\epsilon^{-1}))$
for achieving an $\epsilon$-stationary point of
$\cost + \phi$
is derived in~\cite{Drusvyatskiy-Paquette19},
where 
$\mathcal{O}(\cdot)$
stands for Landau's big-O notation.
Here, a convergence rate is measured in the number of evaluations of {\em basic operations}; see~\cite[Sect. 6]{Drusvyatskiy-Paquette19} and the end of Section~\ref{sec:measure}.
Further theoretical developments on PLM have been reported in~\cite{Duchi-Ruan18,Zhen-Ma-Xue24} for a specific application, called robust phase retrieval.
Multiplier-based methods, exemplified by
{\em augmented Lagrangian methods (ALM)}, e.g.,~\cite{DeMarchi-Jia-Kanzow-Mehlitz23,Bourkhissi-Necoara-Patrinos-TranDinh25},
provide a general-purpose framework for nonlinear constrained optimization to which Problem~\ref{problem:origin} can be reformulated with an auxiliary variable
  $\bm{z}=\mathfrak{S}(\bm{x})$
  (see~\eqref{eq:ALM} and~\cite[Sect. 2]{Bourkhissi-Necoara-Patrinos-TranDinh25}).
For ALM, a subsequential convergence guarantee and an iteration complexity
$\mathcal{O}(\epsilon^{-3})$
for finding a ($\epsilon$-) {\em KKT point} (corresponding to a stationary point) have been established respectively in~\cite{DeMarchi-Jia-Kanzow-Mehlitz23} and~\cite{Bourkhissi-Necoara-Patrinos-TranDinh25}.
We note that such guarantees of ALM typically require some additional assumption on the generated sequence (see, e.g.,~\cite[Prop. 3.2 (i)-(ii)]{DeMarchi-Jia-Kanzow-Mehlitz23},~\cite[Thm. 1]{Bourkhissi-Necoara-Patrinos-TranDinh25}).
PLM and ALM typically require an iterative solver for solving subproblems, exactly or to sufficient accuracy, at every iteration.
However, such an exact solution is not available in general within the finite number of iterations via an iterative solver on which their numerical performance depends.
In contrast, 
{\em subgradient methods (SGM)}, e.g.~\cite{Davis-Drusvyatskiy-MacPhee18,Li-Zhu-So-Man-Cho-Lee19,Zhu-Zhao-Zhang23}, are single-loop and structurally simple, but the currently available convergence rates
$\mathcal{O}(\epsilon^{-4})$
of SGM
do not reach the best-known rate
$\mathcal{O}(\epsilon^{-3})$
(see Table~\ref{table:summary}).

The main contribution of this paper is to establish, for Problem~\ref{problem:origin},
a single-loop algorithm that simultaneously enjoys (i) a subsequential convergence guarantee to a stationary point of
$\cost + \phi$
and (ii) a best-known convergence rate of
$\mathcal{O}(\epsilon^{-3})$
for finding an $\epsilon$-stationary point of
$\cost +\phi$.

More precisely, we extend a {\em proximal variable smoothing} algorithm~\cite{Liu-Xia24}, which has been proposed for Problem~\ref{problem:origin} in a special case where
$\mathfrak{S}$
is a surjective linear operator onto
$\mathcal{Z}$.
The proposed algorithm is designed as a proximal gradient-type method for time-varying functions
$h+\moreau{g}{\mu_{n}}\circ \mathfrak{S} + \phi$
with the {\em Moreau envelope}
$\moreau{g}{\mu}:\mathcal{Z} \to \mathbb{R}$
(see~\eqref{eq:moreau})
of
$g$
with
$\mu \in (0,\eta^{-1})$,
where
$\mu_{n} \in (0,\eta^{-1})$
converges to zero.
We note that the proposed algorithm is single-loop due to the prox-friendliness of
$\phi$
and
$g$
assumed in Problem~\ref{problem:origin}.
We emphasize that a convergence analysis in~\cite{Liu-Xia24} relies on
the surjectivity and linearity of
$\mathfrak{S}$
(see~\eqref{eq:bound_surjective_linear}),
hence we cannot apply the discussion in~\cite{Liu-Xia24} to our nonlinear setting of
$\mathfrak{S}$.
To address this gap, we investigate further properties, beyond~\eqref{eq:bound_surjective_linear}, of a gradient mapping-type stationarity measure employed in~\cite{Liu-Xia24} (see Propositions~\ref{proposition:stationary_relation} and~\ref{proposition:measure}).
By leveraging these properties,
regarding a stationary point,
we establish a subsequential convergence of the proposed algorithm under both diminishing stepsize and stepsize obtained by a standard {\em backtracking algorithm} (see Theorem~\ref{theorem:convergence_extension}).
For a diminishing stepsize, we additionally derive a
convergence rate
$\mathcal{O}(\epsilon^{-3})$
for finding an $\epsilon$-stationary point of Problem~\ref{problem:origin}
(see~\eqref{eq:basic_operation} and Remark~\ref{remark:reasonable_mu}).

As selected applications of Problem~\ref{problem:origin}, we present brief introductions in Section~\ref{sec:application} to a constrained finite-max nonlinear minimization (e.g., {\em maxmin dispersion problem}~\cite{Haines-Loeppky-Tseng-Wang13,Pan-Shen-Xu21,Lopes-Peres-Bilches25} and {\em robust target localization}~\cite{Domingos-Xavier24,Domingos-Xavier25}), and a {\em multiple-input-multiple-output (MIMO) signal detection}~\cite{Yang-Hanzo15,Chen19,Hayakawa-Hayashi20,Shoji-Yata-Kume-Yamada25}.
In particular, for MIMO signal detection, we present a new optimization model capturing a certain periodic property of the target signal by the nonlinear sine function.
Numerical experiments in Section~\ref{sec:numerical} demonstrate the effectiveness of the proposed algorithm in these two application scenarios.

\begin{table}
  \caption{Summary of existing algorithms for Problem~\ref{problem:origin}}\label{table:summary}
  \centering
  \begin{threeparttable}
    {
      \begin{tabular}[c]{c@{\hspace{1pt}}c@{}|c@{\hspace{3pt}}c@{\hspace{3pt}}c|@{}c@{\hspace{1em}}c}
        \toprule
        \multirow{2.0}{*}{Method}                                  &\multirow{2.5}{*}{Ref.}
                                           &
        \multicolumn{3}{c|}{Assumptions\tnote{a}}    &
        \multicolumn{2}{c}{Convergence result\tnote{b}}
        \\
        \cmidrule{3-5}
        \cmidrule{6-7}
        \multirow{1.2}{*}{(Structure)}       &
                  &
        \multicolumn{1}{c}{$g$}            &
        \multicolumn{1}{c}{$\mathfrak{S}$} &
        \multicolumn{1}{c|}{$\phi$}        &
        \multicolumn{1}{c}{Subseq. conv.}   &
        \multicolumn{1}{c}{\hspace{0em}Rate}
        \\
        \midrule
        \multirow{1}{*}{PLM}

                                                & \cite{Drusvyatskiy-Paquette19}     & cvx                             & ${C^{1}}$                  & {$\Gamma_{0}(\mathcal{X})$}                 & -                       & \hspace{0em}$\widetilde{\mathcal{O}}(\epsilon^{-3})$ \\
        \multirow{1}{*}{(Double-loop)}
                                                & \cite{Lewis-Wright16}              & {w-c}                  & ${C^{1}}$                  & {$\Gamma_{0}(\mathcal{X})$}                 &  \checkmark                  & \hspace{0em}-                                               \\
                                                \midrule
                                                \multirow{1}{*}{ALM\tnote{c}}
                                                & \cite{DeMarchi-Jia-Kanzow-Mehlitz23}         & lsc & ${C^{1}}$                  & { lsc}         &  \checkmark                  & \hspace{0em}-                                               \\

                                                \multirow{1}{*}{(Double-loop)}
                                                & \cite{Bourkhissi-Necoara-Patrinos-TranDinh25}            & lsc                             & ${C^{1}}$                  & { lsc}                       & -                       & \hspace{0em}$\mathcal{O}(\epsilon^{-3})$\tnote{d}                    \\
        \midrule
        \multirow{1}{*}{SGM\tnote{e}}
                                           & \multirow{2}{*}{\cite{Davis-Drusvyatskiy-MacPhee18,Li-Zhu-So-Man-Cho-Lee19,Zhu-Zhao-Zhang23}}   & \multirow{2}{*}{cvx}                  & \multirow{2}{*}{${C^{1}}$}                         & \multirow{2}{*}{$\Gamma_{0}(\mathcal{X})$}        & \multirow{2}{*}{-}                       & \multirow{2}{*}{\hspace{0em}$\mathcal{O}(\epsilon^{-4})$}                   \\
        \multirow{1}{*}{(Single-loop)} & & & & & & \\
        \midrule
        \multirow{5}{*}{Variable}
                                           & \cite{Bohm-Wright21}               & {w-c}                  & lin                           & 0                 & -                       & \hspace{0em}$\mathcal{O}(\epsilon^{-3})$                    \\
                                           \multirow{5}{*}{smoothing}& \cite{Kume-Yamada24A}             & {{w-c}} & {${C^{1}}$} & 0 & { \checkmark} & \hspace{0em}{$\mathcal{O}(\epsilon^{-3})$}   \\
                                           \multirow{5}{*}{(Single-loop)}& \cite{Yazawa-Kume-Yamada26}             & {{DC}} & {${C^{1}}$} & 0 & { \checkmark} & \hspace{0em} -   \\
                                           & \cite{Lopes-Peres-Bilches25}               & {w-c} & lin                           & $\iota_{C}$\tnote{f}        & - &\hspace{0em} $\mathcal{O}(\epsilon^{-3})$                    \\
                                           & \cite{Liu-Xia24}           & {w-c}                  & sur, lin                           & {$\Gamma_{0}(\mathcal{X})$}        & -                       & \hspace{0em}$\mathcal{O}(\epsilon^{-3})$                    \\
                                           & Alg.~\ref{alg:proposed}             & {{w-c}} & {${C^{1}}$} & $\Gamma_{0}(\mathcal{X})$ & { \checkmark} & \hspace{0em}{$\mathcal{O}(\epsilon^{-3})$\tnote{g}}   \\
        \bottomrule
      \end{tabular}
      \begin{tablenotes}[flushleft]
        \item[a] ``cvx'',  ``w-c'' and `DC'' denote convex, weakly convex, and difference-of-convex function respectively.
    ``$ C^{1}$'', ``lin'' and ``sur'' denote continuously differentiable mapping, linear operator and surjective respectively.
    ``lsc'' denotes proper lower semicontinuous functions.
    ``$\Gamma_{0}(\mathcal{X})$'' denotes  proper lower semicontinuous convex function.
      \item[b]
        ``Subseq. conv.'' stands for whether a subsequential convergence guarantee of a stationary point is shown or not.
        ``Rate'' denotes the number of evaluations of {\em basic operations}~\cite{Drusvyatskiy-Paquette19}  for achieving an $\epsilon$-stationary point, where
        $\widetilde{\mathcal{O}}(\cdot)$
        hides logarithmic term in Landau's big-O $\mathcal{O}(\cdot)$.
  \item[c]
    Convergence guarantee of ALM~\cite{DeMarchi-Jia-Kanzow-Mehlitz23,Bourkhissi-Necoara-Patrinos-TranDinh25} requires some additional assumption on the generated sequence.
  \item[d] The rate in~\cite{Bourkhissi-Necoara-Patrinos-TranDinh25} means the number of outer iterations for an $\epsilon$-stationary point.
  \item[e]
    Subgradient methods~\cite{Davis-Drusvyatskiy-MacPhee18,Li-Zhu-So-Man-Cho-Lee19,Zhu-Zhao-Zhang23} are applicable to Problem~\ref{problem:origin} if
    $\cost$
    is weakly convex, which is guaranteed, e.g., under the convexity of
    $g$;
    see Remark~\ref{remark:composite_non_weak_convexity}.
  \item [f]
          $C \subset \mathcal{X}$ is a nonempty closed subspace of
          $\mathcal{X}$.
          See~\eqref{eq:indicator} for
        $\iota_{C}$.

  \item [g]
        This convergence rate is established for the proposed algorithm with
        diminishing stepsizes in Example~\ref{example:stepsize:diminishing}.
      \end{tablenotes}
    }
  \end{threeparttable}
\end{table}

\textbf{Related works on Moreau-envelope-based smoothing}:
The Moreau envelope of a nonsmooth (weakly) convex function has been used in variable smoothing-type algorithms~\cite{Bot-Hendrich15,Bohm-Wright21,Kume-Yamada24,Liu-Xia24,Kume-Yamada24A,Lopes-Peres-Bilches25,Yazawa-Kume-Yamada26}
for designing a smoothed surrogate function of the cost function.
Table~\ref{table:summary} summarizes the applicability of each variable smoothing-type algorithm.
Compared with variable smoothing-type algorithms,
the proposed algorithm admits
a smooth nonlinear mapping
$\mathfrak{S}$
in
$g\circ\mathfrak{S}$
together with a general prox-friendly convex function
$\phi\in \Gamma_{0}(\mathcal{X})$,
while retaining the best-known
$\mathcal{O}(\epsilon^{-3})$
convergence rate.
We note that the authors recently proposed a variable smoothing-type algorithm~\cite{Yazawa-Kume-Yamada26} where
$\phi\equiv 0$
but
$g$
admits a wider class, namely difference-of-convex functions, than weakly convex functions.
We also note that the minimization of
$h+g\circ\mathfrak{S}$
in Problem~\ref{problem:origin} over a nonempty closed nonconvex set
$C\subset \mathcal{X}$
has been tackled by a variable smoothing with a {\em parametrization strategy}~\cite{Kume-Yamada24,Kume-Yamada24A} in a special case where
$C=\{\param(\bm{y}) \in \mathcal{X}\mid \bm{y} \in \mathcal{Y}\}$
is given with a continuously differentiable mapping
$\param:\mathcal{Y}\to \mathcal{X}$
defined over the Euclidean space
$\mathcal{Y}$
(see, e.g.~\cite{Kume-Yamada21,Kume-Yamada22,Levin-Kileel-Boumal24}, for parametrization strategy).

The idea of Moreau envelope based-smoothing is not limited to variable smoothing-type algorithms.
In~\cite{Drusvyatskiy-Paquette19}, a smoothing variant of proximal linear methods has also been proposed, where the outer convex function
$g$
is replaced by its convex Moreau envelope
$\moreau{g}{\mu}$
with a carefully fixed
$\mu > 0$.
In~\cite{Drusvyatskiy-Paquette19}, although a convergence rate
$\mathcal{O}(\epsilon^{-3}\log(\epsilon^{-1}))$
is derived for a smoothing proximal linear method,
the analysis therein does not naturally lead to a subsequential convergence guarantee to a stationary point of the original cost function because a fixed smoothing parameter
$\mu$
depends on tolerance
$\epsilon>0$.
In contrast, thanks to a variable setting of
$\mu$,
the proposed algorithm achieves a subsequential convergence guarantee together with a slightly better convergence rate
$\mathcal{O}(\epsilon^{-3})$
than~\cite{Drusvyatskiy-Paquette19}.

  {\bf Notation}:
$\mathbb{N}$,
$\mathbb{R}$,
$\mathbb{R}_{+}$
and
$\mathbb{R}_{++}$,
denote respectively the sets of all positive integers, all real numbers, all nonnegative real numbers, and all positive real numbers.
$\stdnorm{\cdot}$
and
$\inprod{\cdot}{\cdot}$
are respectively the Euclidean norm and the standard inner product for vectors, and
$\Id$
stands for the identity operator.
For a linear operator
$A:\mathcal{X}\to \mathcal{Y}$,
its adjoint
$A^{*}:\mathcal{Y}\to \mathcal{X}$
is defined as the linear operator satisfying
$\inprod{\bm{x}}{A^{*}\bm{y}} = \inprod{A\bm{x}}{\bm{y}}\ (\forall \bm{x}\in\mathcal{X},\forall \bm{y}\in\mathcal{Y})$,
and $\opnorm{A}\coloneqq \sup_{\stdnorm{\bm{x}}\leq 1, \bm{x}\in \mathcal{X}} \stdnorm{A\bm{x}}$
denotes the operator norm of
$A$,
i.e., the largest singular value of
$A$
in this case.
For a differentiable mapping
$\mathcal{F}:\mathcal{X}\to \mathcal{Z}$,
the linear operator
$\mathrm{D}\mathcal{F}(\widebar{\bm{x}}):\mathcal{X} \to \mathcal{Z}:\bm{v}\mapsto \lim_{\mathbb{R}\setminus\{0\} \ni t\to 0} t^{-1}(\mathcal{F}(\widebar{\bm{x}}+t\bm{v}) -\mathcal{F}(\widebar{\bm{x}}))$
denotes the G\^{a}teaux derivative of
$\mathcal{F}$
at
$\widebar{\bm{x}}\in \mathcal{X}$.
For a differentiable function
$J:\mathcal{X} \to \mathbb{R}$,
$\nabla J(\widebar{\bm{x}}) \in \mathcal{X}$
denotes the gradient of
$J$
at
$\widebar{\bm{x}} \in \mathcal{X}$
if
$\mathrm{D}J(\widebar{\bm{x}})[\bm{v}] = \inprod{\nabla J(\widebar{\bm{x}})}{\bm{v}}\ (\bm{v} \in \mathcal{X})$.

\section{Preliminary on Nonsmooth Analysis}
We review nonsmooth analysis in~\cite{Rockafellar-Wets98}
(see a recent comprehensive survey~\cite{Li-So-Ma20} in Signal Processing Magazine).

\subsection{Subdifferentials and Stationary Point} \label{sec:preliminary:subdifferential}
A function
$J:\mathcal{X} \to \exR$
is called (a) proper if
$\dom{J}\coloneqq \{\bm{x} \in \mathcal{X} \mid J(\bm{x}) < +\infty\} \neq \emptyset$;
(b)~lower semicontinuous if
$\{(\bm{x},a) \in \mathcal{X}\times \mathbb{R} \mid J(\bm{x})\leq a\} \subset \mathcal{X}\times \mathbb{R}$
is closed in
$\mathcal{X} \times \mathbb{R}$;
(c) convex if
$J(t\bm{x}_{1}+(1-t)\bm{x}_{2}) \leq t J(\bm{x}_{1}) + (1-t)J(\bm{x}_{2})\ (\forall \bm{x}_{1},\bm{x}_{2} \in \mathcal{X}, \forall t \in (0,1))$.
Moreover,
$\Gamma_{0}(\mathcal{X})$
denotes the set of all proper lower semicontinuous convex functions defined over
$\mathcal{X}$.
For a convex function
$J\in \Gamma_{0}(\mathcal{X})$,
the convex subdifferential of
$J$
is defined by
\begin{equation}
  \subdiff J(\widebar{\bm{x}})\coloneqq \{\bm{v} \in \mathcal{X} \mid J(\bm{x}) \geq J(\widebar{\bm{x}}) + \inprod{\bm{v}}{\bm{x}-\widebar{\bm{x}}}\ (\forall \bm{x} \in\mathcal{X})\}.
\end{equation}
The convex subdifferential is extended to the nonconvex case.
\begin{definition}[Subdifferentials of nonconvex function~{\cite[Def. 8.3]{Rockafellar-Wets98}}]\label{definition:subdifferential}
  For
  a proper lower semicontinuous function
  $J:\mathcal{X} \to \exR$,
  the {\em limiting (or general) subdifferential}
  $\Lsubdiff J:\mathcal{X} \rightrightarrows \mathcal{X}$
  of
  $J$
  at
  $\widebar{\bm{x}} \in \dom{J}$
  is defined by
    {
      \thickmuskip=0.1\thickmuskip
      \medmuskip=0.1\medmuskip
      \thinmuskip=0.1\thinmuskip
      \arraycolsep=0.1\arraycolsep
      \begin{equation}
        \Lsubdiff J(\widebar{\bm{x}})
        \coloneqq
        \left\{\bm{v}\in\mathcal{X} \;\middle|\;
        \begin{array}{l}
          \exists (\bm{x}_{n})_{n=1}^{\infty}(\subset \mathcal{X})\to \widebar{\bm{x}},\
          \exists \bm{v}_{n} \in \Fsubdiff J(\bm{x}_{n}) \\
          {\rm s.t.}\
          \lim_{n\to\infty}J(\bm{x}_{n})=J(\widebar{\bm{x}}),\
          \lim_{n\to\infty}\bm{v}_{n}=\bm{v}
        \end{array}
        \right\}.
        \hspace{-1em}
        \label{eq:general_subdifferential}
      \end{equation}}%
  Here, the {\em Fr\'echet (or regular) subdifferential}
  $\Fsubdiff J: \mathcal{X} \rightrightarrows \mathcal{X}$
  at
  $\widehat{\bm{x}} \in \dom{J}$
  is defined
  as the set of all
  $\bm{w} \in \mathcal{X}$
  satisfying
  \begin{equation}
    \sup\limits_{\epsilon>0} \left(\inf\limits_{0<\stdnorm{\bm{x}-\widehat{\bm{x}}}< \epsilon}\frac{J(\bm{x})-J(\widehat{\bm{x}}) - \inprod{\bm{w}}{\bm{x}-\widehat{\bm{x}}}}{\stdnorm{\bm{x}-\widehat{\bm{x}}}}\right) \geq 0.
  \end{equation}
  For
  $\widebar{\bm{x}} \notin \dom{J}$,
  $\Fsubdiff J(\widebar{\bm{x}})$.
  and
  $\Lsubdiff  J(\widebar{\bm{x}})$
  are defined as
  $\emptyset$.
\end{definition}
For special functions
$J$,
we get simple expressions of
$\Fsubdiff J(\widebar{\bm{x}})$
and
$\Lsubdiff J(\widebar{\bm{x}})\ (\widebar{\bm{x}} \in \mathcal{X})$.
Under the continuous differentiability of
$J$,
$\Fsubdiff J(\widebar{\bm{x}}) = \Lsubdiff J(\widebar{\bm{x}})= \{\nabla J(\widebar{\bm{x}})\}$
holds~\cite[Exe. 8.8 (b)]{Rockafellar-Wets98}.
For an
$\widetilde{\eta}(\geq 0)$-weakly convex function
$J$,
$\widetilde{J}\coloneqq J+\frac{\widetilde{\eta}}{2}\stdnorm{\cdot}^{2}$
is convex, and we get
$\Fsubdiff J(\widebar{\bm{x}}) = \Lsubdiff J(\widebar{\bm{x}}) = \subdiff \widetilde{J} (\widebar{\bm{x}}) - \widetilde{\eta} \widebar{\bm{x}} \coloneqq \{\bm{v}-\widetilde{\eta}\widebar{\bm{x}}\mid \bm{v}\in\subdiff\widetilde{J}(\widebar{\bm{x}})\}$~\cite[Prop. 6.3]{Bauschke-Moursi-Wang21}.

For
$\cost$
in Problem~\ref{problem:origin},
$\Fsubdiff \cost = \Lsubdiff \cost$
holds by~\cite[Lemma 2.4]{Kume-Yamada24A} with~\cite[Cor. 8.11]{Rockafellar-Wets98}.
Since useful asymptotic properties have been investigated for the limiting subdifferential (see~\cite[Ch. 8]{Rockafellar-Wets98}), we use
$\Lsubdiff \cost$
in the subsequent discussion.
We also note that
$\Fsubdiff \cost(\widebar{\bm{x}}) = \Lsubdiff \cost(\widebar{\bm{x}})$
is not identical to
$\subdiff (\cost+\frac{\widetilde{\eta}}{2}\stdnorm{\cdot}^{2})(\widebar{\bm{x}})-\widetilde{\eta} \widebar{\bm{x}}$
with any
$\widetilde{\eta} \geq 0$
in general
because
$\cost=h+g\circ\mathfrak{S}$
is not necessarily weakly convex
(see Remark~\ref{remark:composite_non_weak_convexity}).
Nevertheless, by using the subdifferential
$\Lsubdiff g(\bm{z}) = \subdiff (g+\frac{\eta}{2}\stdnorm{\cdot}^{2})(\bm{z}) - \eta\bm{z}\ (\bm{z} \in \mathcal{Z})$
of the $\eta$-weakly convex function
$g$,
we have
\begin{equation}
(\widebar{\bm{x}} \in \mathcal{X}) 
   \quad \Lsubdiff \cost(\widebar{\bm{x}}) = \nabla h(\widebar{\bm{x}}) + (\mathrm{D}\mathfrak{S}(\widebar{\bm{x}}))^{*}[\Lsubdiff g(\mathfrak{S}(\widebar{\bm{x}}))] 
   \label{eq:chain_rule}
\end{equation}
from the sum rule~\cite[Exe. 8.8]{Rockafellar-Wets98}, and the chain rule~\cite[Thm. 10.6]{Rockafellar-Wets98} with~\cite[Thm. 9.13]{Rockafellar-Wets98} and the Lipschitz continuity of
$g$,
where
$(\mathrm{D}\mathfrak{S}(\widebar{\bm{x}}))^{*}[\Lsubdiff g(\mathfrak{S}(\widebar{\bm{x}}))] \coloneqq \{(\mathrm{D}\mathfrak{S}(\widebar{\bm{x}}))^{*}[\bm{d}]\in \mathcal{X} \mid \bm{d} \in \Lsubdiff g(\mathfrak{S}(\widebar{\bm{x}}))\}$.
Moreover, for Problem~\ref{problem:origin},
we have
\begin{equation}
   (\widebar{\bm{x}} \in \mathcal{X})
    \quad
   \Lsubdiff (\cost+\phi)(\widebar{\bm{x}}) = \Lsubdiff \cost(\widebar{\bm{x}}) + \subdiff \phi(\widebar{\bm{x}}) \label{eq:sum_rule}
\end{equation}
from the sum rule~\cite[Cor. 10.9 and Exe. 10.10]{Rockafellar-Wets98} (see a similar discussion in, e.g.,~\cite[Lemma 2.4]{Kume-Yamada24A}),
where we used Minkowski sum
$\Lsubdiff \cost(\widebar{\bm{x}}) + \subdiff \phi(\widebar{\bm{x}})  \coloneqq \{\bm{v}_{1}+\bm{v}_{2} \mid \bm{v}_{1}\in \Lsubdiff \cost(\widebar{\bm{x}}), \bm{v}_{2} \in \subdiff\phi(\widebar{\bm{x}})\} \subset \mathcal{X}$
(see, e.g.,~\cite[p.25]{Rockafellar-Wets98}).

\begin{definition}[Stationarity] \label{definition:optimality}
  For
  $\cost + \phi$
  in Problem~\ref{problem:origin},
  \begin{enumerate}[label=(\alph*)]
    \item  \label{enum:definition:optimality:stationary}
          $\widebar{\bm{x}}\in \dom{\phi}$
          is said to be
          a {\em stationary point}
          of
          $\cost + \phi$
          if
          $\Lsubdiff (\cost+\phi)(\widebar{\bm{x}}) \ni \bm{0}$.
    \item  \label{enum:definition:optimality:epsilon}
          $\widebar{\bm{x}}\in \dom{\phi}$
          is said to be
          an {\em $\epsilon$-stationary point} of
          $\cost + \phi$
          with
          $\epsilon \geq 0$
          if there exists
          $\widebar{\bm{z}} \in \mathcal{Z}$
          satisfying
          \begin{equation}
            \hspace{-2em}
            \begin{cases}
              \dist\left(\bm{0}, \nabla h(\widebar{\bm{x}}) + (\mathrm{D}\mathfrak{S}(\widebar{\bm{x}}))^{*}[\Lsubdiff g(\widebar{\bm{z}})]+\subdiff \phi(\widebar{\bm{x}})\right) \leq \epsilon; \\
              \stdnorm{\mathfrak{S}(\widebar{\bm{x}}) - \widebar{\bm{z}}} \leq \epsilon,
            \end{cases}\hspace{-2em} \label{eq:KKT}
          \end{equation}
          where
          $\dist:\mathcal{X}\times 2^{\mathcal{X}} \to [0,+\infty]$
          stands for the distance
            $\dist(\widebar{\bm{x}},E)\coloneqq
            \left\{
            \begin{array}{cl}
              \inf\{\stdnorm{\bm{v}-\widebar{\bm{x}}}\mid \bm{v} \in E\}, & \mathrm{if}\  E\neq \emptyset; \\
              +\infty,                                                 & \mathrm{if}\  E=\emptyset.
            \end{array}
            \right.$
  \end{enumerate}
\end{definition}

By Fermat's rule~\cite[Thm. 10.1]{Rockafellar-Wets98},
every local minimizer of
$\cost + \phi$
is a stationary point of
$\cost+\phi$.
Hence, the stationarity has been widely used as a tractable necessary condition for local optimality in the literature, e.g.,~\cite{Lewis-Wright16,Duchi-Ruan18,Drusvyatskiy-Paquette19,Zhen-Ma-Xue24,DeMarchi-Jia-Kanzow-Mehlitz23,Bourkhissi-Necoara-Patrinos-TranDinh25,Davis-Drusvyatskiy-MacPhee18,Li-Zhu-So-Man-Cho-Lee19,Zhu-Zhao-Zhang23,Liu-Xia24,Lopes-Peres-Bilches25}.

\begin{remark}[On $\epsilon$-stationary point]
  If
  $\widebar{\bm{x}}\in \dom{\phi}$
  is an $\epsilon$-stationary point of
  $\cost + \phi$
  with
  $\epsilon = 0$,
  then
  $\widebar{\bm{x}}$
  is a stationary point of
  $\cost + \phi$
  by~\eqref{eq:chain_rule} and~\eqref{eq:sum_rule} with
  $\widebar{\bm{z}} = \mathfrak{S}(\widebar{\bm{x}})$
  in~\eqref{eq:KKT}.
  Moreover, the condition~\eqref{eq:KKT} with
  $\epsilon = 0$
  coincides with a KKT-type condition, e.g.,~\cite[Def. 2.3]{DeMarchi-Jia-Kanzow-Mehlitz23}, for a given
  $(\widebar{\bm{x}},\widebar{\bm{z}})\in \mathcal{X}\times \mathcal{Z}$:
  \begin{equation}
    (\exists \bm{d} \in \mathcal{Z}) \quad
    \begin{cases}
      \bm{0} \in \nabla h(\widebar{\bm{x}}) + (\mathrm{D}\mathfrak{S}(\widebar{\bm{x}}))^{*}[\bm{d}] + \subdiff \phi(\widebar{\bm{x}}), \\
      \bm{d} \in \Lsubdiff g(\widebar{\bm{z}}),\ \mathrm{and}\ 
      \mathfrak{S}(\widebar{\bm{x}}) - \widebar{\bm{z}} = \bm{0}
    \end{cases}
  \end{equation}
  of
  the optimization problem (see also~\cite[Sect. 2]{Bourkhissi-Necoara-Patrinos-TranDinh25})
  \begin{equation}
    \mathop{\mathrm{minimize}}\limits_{\bm{x}\in\mathcal{X},\bm{z}\in\mathcal{Z}}\
    h(\bm{x}) + g(\bm{z}) + \phi(\bm{x})\ \mathrm{s.t.}\ \mathfrak{S}(\bm{x}) - \bm{z} = \bm{0}. \label{eq:ALM}
  \end{equation}
  Hence, our $\epsilon$-stationarity in Def.~\ref{definition:optimality}~\ref{enum:definition:optimality:epsilon} can be viewed as an approximate KKT condition~\cite[Def. 4]{Bourkhissi-Necoara-Patrinos-TranDinh25}.
  The condition~\eqref{eq:KKT}
  also appears in~\cite{Bohm-Wright21} under the linearity of
  $\mathfrak{S}$
  and
  $\phi\equiv 0$.
\end{remark}

\begin{remark}[$\cost$ is not necessarily weakly convex] \label{remark:composite_non_weak_convexity}
  In a special case where
  $g$
  is convex,
  $\cost=h+g\circ\mathfrak{S}$
  in Problem~\ref{problem:origin} is weakly convex~\cite[Lemma 4.2]{Drusvyatskiy-Paquette19}.
  On the other hand, under the weak convexity of
  $g$,
  $\cost$
  is not weakly convex in general.
  As a counterexample, consider
  $h\equiv 0$,
  $\mathfrak{S}:\mathbb{R}^{2}\to \mathbb{R}:(x,y)^{\TT}\mapsto xy$,
  and
  MCP ($1$-weakly convex function)~\cite{Zhang10}
  $g(t)\coloneqq \abs{t}-t^{2}/2$
  if
  $\abs{t} \leq 1$;
  $g(t)\coloneqq 2^{-1}$
  otherwise.
  Clearly,
  $\cost=h+g\circ\mathfrak{S}$
  satisfies the conditions in Problem~\ref{problem:origin}.
  Assume contrarily that
  $\cost$
  is $\widetilde{\eta}$-weakly convex with some
  $\widetilde{\eta} > 0$,
  i.e.,
  $\widetilde{\cost}\coloneqq \cost+\frac{\widetilde{\eta}}{2}\stdnorm{\cdot}^{2}$
  is convex.
  Consider a region
  $\mathcal{D}\coloneqq \{(x,y)^{\TT}\in \mathbb{R}^{2} \mid 0 < xy <1\}$
  on which
  $\widetilde{\cost}(x,y) = xy - \frac{(xy)^{2}}{2}+\frac{\widetilde{\eta}(x^{2}+y^{2})}{2}$
  is twice continuously differentiable.
  By the convexity of
  $\widetilde{\cost}$,
  its Hessian matrix
  $\nabla^{2}\widetilde{\cost}(x,y) =\begin{bmatrix} -y^{2}+\widetilde{\eta} & 1-xy \\ 1-xy & -x^{2} + \widetilde{\eta} \end{bmatrix}$
  must be positive semidefinite over
  $\mathcal{D}$~\cite[Thm. 2.14]{Rockafellar-Wets98}.
  However,
  for
  $(\widebar{x},\widebar{y})^{\TT}\coloneqq ((4\sqrt{\widetilde{\eta}})^{-1}, 2\sqrt{\widetilde{\eta}})^{\TT} \in \mathcal{D}$,
  and
  $\bm{v}\coloneqq (1,0)^{\TT} \in \mathbb{R}^{2}$,
  we have
  $\bm{v}^{\TT}\nabla^{2} \widetilde{\cost}(\widebar{x},\widebar{y}) \bm{v} = -\widebar{y}^{2} + \widetilde{\eta} = -3\widetilde{\eta} < 0$,
  which is absurd.
\end{remark}

\vspace{-0.3em}
\subsection{Proximity Operator and Moreau Envelope}\label{sec:preliminary:Moreau}
{\em The proximity operator} and {\em the Moreau envelope} have been used as computational tools for nonsmooth optimization, e.g.,~\cite{Rockafellar-Wets98,Yamada-Yukawa-Yamagishi11,Bauschke-Combettes17,Bohm-Wright21,Kume-Yamada24A,Liu-Xia24}.
For a proper lower semicontinuous
$\eta(>0)$-weakly
convex function
$J:\mathcal{X}\to\exR$,
its proximity operator and Moreau envelope of index
$\mu \in (0,\eta^{-1})$
are respectively defined by
\begin{align}
  \hspace{-1em}
  \prox{\mu J}:\mathcal{X}\to\mathcal{X}:   & \widebar{\bm{x}} \mapsto    \argmin_{\bm{x} \in \mathcal{X}} \left(J(\bm{x}) + \frac{1}{2\mu}\|\bm{x}-\widebar{\bm{x}}\|^{2}\right); \label{eq:prox_definition} \\
  \moreau{J}{\mu}:\mathcal{X}\to\mathbb{R}: & \widebar{\bm{x}} \mapsto \min_{\bm{x} \in \mathcal{X}} \left(J(\bm{x}) + \frac{1}{2\mu}\|\bm{x}-\widebar{\bm{x}}\|^{2}\right),
  \label{eq:moreau}
\end{align}
where
$\prox{\mu J}(\widebar{\bm{x}}) \in \dom{J} \subset \mathcal{X}$
for
$\widebar{\bm{x}} \in \mathcal{X}$
is single-valued due to the strong convexity of
$J + (2\mu)^{-1}\|\cdot-\widebar{\bm{x}}\|^{2}$.
Commonly-used sparsity-promoting functions or robust loss functions
$J$
have the closed-form expressions of
$\prox{\mu J}$
and
$\moreau{J}{\mu}$.
Such functions
$J$
include, e.g.,
$\ell_{1}$-norm, MCP, and SCAD (see, e.g.,~\cite{Bohm-Wright21} and~\cite[Sect. 2.2]{Kume-Yamada24A}).
Moreover,
the metric projection
$P_{C}:\mathcal{X}\to\mathcal{X}:\widebar{\bm{x}}\mapsto \argmin_{\bm{x} \in C}\stdnorm{\bm{x}-\widebar{\bm{x}}}$
onto a nonempty closed convex set
$C\subset \mathcal{X}$
is given by
$P_{C} = \prox{\gamma \iota_{C}}\ (\gamma > 0)$.

The Moreau envelope has been used for designing a smoothed surrogate function that approximates smoothly a nonsmooth function, e.g.,~\cite{Bot-Hendrich15,Burke-Hoheisel17,Bohm-Wright21,Liu-Xia24,Kume-Yamada24A,Lopes-Peres-Bilches25}.
Along this direction, we define a smoothed surrogate function of
$\cost=h+g\circ\mathfrak{S}$
in Problem~\ref{problem:origin}
as
\begin{equation}
  (\mu \in (0,\eta^{-1}), \widebar{\bm{x}} \in \mathcal{X}) \quad
  \cost^{\langle \mu \rangle}(\widebar{\bm{x}}) \coloneqq (h+\moreau{g}{\mu}\circ\mathfrak{S})(\widebar{\bm{x}}) \label{eq:surrogate_function}
\end{equation}
with the gradient (by using the chain rule)
\begin{equation}
  \nabla \cost^{\langle \mu \rangle}(\widebar{\bm{x}}) =
  \nabla h(\widebar{\bm{x}}) + (\mathrm{D}\mathfrak{S}(\widebar{\bm{x}}))^{*}[\nabla \moreau{g}{\mu}(\mathfrak{S}(\widebar{\bm{x}}))]. \label{eq:gradient_expression}
\end{equation}
For each
$\mu \in (0,\eta^{-1})$,
the gradient of
$\moreau{g}{\mu}$
\begin{equation}
  (\widebar{\bm{z}} \in \mathcal{Z}) \quad \nabla \moreau{g}{\mu}(\widebar{\bm{z}}) = \mu^{-1}(\widebar{\bm{z}}-\prox{\mu g}(\widebar{\bm{z}})), \label{eq:Moreau_grad}
\end{equation}
is
$\max\{\mu^{-1},\frac{\eta}{1-\eta\mu^{-1}}\}$-Lipschitz continuous~\cite[Cor.~3.4 b)]{Hoheisel-Laborde-Oberman20}.
Fact~\ref{lemma:basic} below collects key properties of
$\cost^{\langle \mu \rangle}$ for our analysis.
\begin{fact}[Basic properties of $\cost^{\langle \mu \rangle}$~{\cite[Lemma 4.1]{Kume-Yamada24A}}]
  \label{lemma:basic}
  Consider
  $\cost$
  in Problem~\ref{problem:origin},
  and
  $\cost^{\langle \mu \rangle}$
  in~\eqref{eq:gradient_expression}.
  Then, we have
  \begin{equation}
    \thickmuskip=0.3\thickmuskip
    \medmuskip=0.3\medmuskip
    \thinmuskip=0.3\thinmuskip
    (\mu\in(0,\eta^{-1}),\ \widebar{\bm{x}}\in\mathcal{X})\
    \cost^{\langle \mu \rangle}(\widebar{\bm{x}}) \leq \cost(\widebar{\bm{x}}) \leq \cost^{\langle \mu \rangle} (\widebar{\bm{x}}) +  L_{g}^{2}\mu, \label{eq:smoothed_uniformly}
  \end{equation}
  \begin{align}
    & (\mu_{1},\mu_{2}\in (0,\eta^{-1})\ \mathrm{s.t.}\ \mu_{2}\leq \mu_{1},\ \widebar{\bm{x}}\in\mathcal{X})\\
    & \quad \cost^{\langle \mu_{1} \rangle}(\widebar{\bm{x}})
    \leq \cost^{\langle \mu_{2} \rangle}(\widebar{\bm{x}})
    \leq \cost^{\langle \mu_{1} \rangle}(\widebar{\bm{x}}) + L_{g}^{2}(\mu_{1}-\mu_{2}). \label{eq:smoothed_uniformly_different_parameter}
  \end{align}
\end{fact}

Although all derivatives of
$h$,
$\mathfrak{S}$,
and
$\moreau{g}{\mu}$
are Lipschitz continuous (see also Problem~\ref{problem:origin}~\ref{enum:problem:origin:h}-\ref{enum:problem:origin:S}, and~\cite[Cor.~3.4~b)]{Hoheisel-Laborde-Oberman20}),
the Lipschitz continuity of the gradient
$\nabla \cost^{\langle \mu \rangle}$
with fixed
$\mu \in (0,\eta^{-1})$
is not guaranteed automatically unless
$\mathfrak{S}$
is a linear operator.
The Lipschitz continuity of the gradient plays a key role even in smooth optimization, e.g.,~\cite{Beck17}.
In this paper,
we explicitly assume the Lipschitz continuity of
$\nabla \cost^{\langle \mu \rangle}$.

\begin{assumption}[Lipschitz continuity of $\nabla \cost^{\langle \mu \rangle}$] \label{assumption:Lipschitz_smoothed}
  For each
  $\mu \in (0,(2\eta)^{-1}]$,
  $\nabla \cost^{\langle \mu \rangle}$
  in~\eqref{eq:gradient_expression} is
  $L_{\nabla \cost^{\langle \mu \rangle}}$-Lipschitz continuous over
  $\dom{\phi}$,
  where
  $L_{\nabla \cost^{\langle \mu \rangle}} \coloneqq \varpi_{1}+\varpi_{2}\mu^{-1} > 0$
  is given with constants
  $\varpi_{1} \geq 0$
  and
  $\varpi_{2} > 0$
  independent of
  $\mu$.
\end{assumption}

\begin{example}[On Assumption~\ref{assumption:Lipschitz_smoothed}] \label{example:assumption:Lipschitz_smoothed}
  Assumption~\ref{assumption:Lipschitz_smoothed}
  is satisfied if
  $\mathfrak{S}$
  is
  $\kappa_{\mathfrak{S}}(>0)$-Lipschitz continuous over
  $\dom{\phi}$
  (or equivalently,
  $\sup_{\bm{x} \in \dom{\phi}}\opnorm{\mathrm{D}\mathfrak{S}(\bm{x})}\leq \kappa_{\mathfrak{S}}$)
  by~\cite[Prop.~4.5~(a)]{Kume-Yamada24A},
  where constants are given by
  $\varpi_{1} = L_{\nabla h} + L_{g}L_{\mathrm{D}\mathfrak{S}}$
  and
  $\varpi_{2} = \kappa_{\mathfrak{S}}^{2}$.
  Hence, Assumption~\ref{assumption:Lipschitz_smoothed} is satisfied if
  (i)
  $\mathfrak{S}$
  is linear;
  or
  (ii)
  $\mathfrak{S}$
  is twice continuously differentiable and
  $\dom{\phi}$
  is compact.
  We note sufficient conditions for Assumption~\ref{assumption:Lipschitz_smoothed} are not limited to the Lipschitz continuity of
  $\mathfrak{S}$
  over
  $\dom{\phi}$ (see~\cite[Lemma C.2 and Prop. III.3 (a)]{Yazawa-Kume-Yamada26}).
\end{example}

\section{Proximal Variable Smoothing for Problem~\ref{problem:origin}} \label{sec:proposed}
\begin{algorithm}[t]
  \caption{Proximal variable smoothing for Problem~\ref{problem:origin}}
  \label{alg:proposed}
  {
    \begin{algorithmic}[0]
      \Require
      $\bm{x}_{1}\in \dom{\phi},c\in(0,1)$
      \For{$n=1,2,\ldots$}
      \State
      Set
      {
 \thickmuskip=0\thickmuskip
\medmuskip=0\medmuskip
\thinmuskip=0\thinmuskip
\arraycolsep=0\arraycolsep
$\cost_{n}\coloneqq\cost^{\langle \mu_{n} \rangle}= h+\moreau{g}{\mu_{n}}\circ \mathfrak{S}$}
      with $\mu_{n} \in (0,\frac{1}{2\eta}]$
      in~\eqref{eq:nonsummable}
      \State
      Find
      $\gamma_{n} > 0$
      satisfying~\eqref{eq:Armijo} with $\gamma\coloneqq \gamma_{n}$ and $c$
      \State
      $\bm{x}_{n+1} \leftarrow \prox{\gamma_{n}\phi}(\bm{x}_{n}-\gamma_{n}\nabla \cost_{n}(\bm{x}_{n})) \in \dom{\phi}$
      \EndFor
      \Ensure
      $\bm{x}_{n} \in \dom{\phi}$
    \end{algorithmic}
  }
\end{algorithm}
\begin{algorithm}[t]
  \caption{Backtracking algorithm to find $\gamma_{n}$ satisfying~\eqref{eq:Armijo} with $\gamma\coloneqq \gamma_{n}$ at $n$th iteration in Algorithm~\ref{alg:proposed} (see Exm.~\ref{example:stepsize:backtracking})}
  \label{alg:backtracking}
  \begin{algorithmic}
    \Require
    $c \in (0,1),\rho \in (0, 1),\ \gamma_{\rm init} > 0$
    \State
    $\gamma \leftarrow \gamma_{\rm init}$
    \While{Condition~\eqref{eq:Armijo} with $\gamma$ and $c$ does not hold}
    \State $\gamma \leftarrow \rho \gamma$
    \EndWhile
    \State
    $\gamma_{n}\leftarrow \gamma$
    \Ensure \thickmuskip=0\thickmuskip
    \medmuskip=0\medmuskip
    \thinmuskip=0\thinmuskip
    \arraycolsep=0\arraycolsep
    \small
    $\gamma_{n}= \max\{\gamma_{\rm init}\rho^{l} \mid l\in \mathbb{N}\cup \{0\},\ \mathrm{\eqref{eq:Armijo}\ holds\ with\ }\gamma \coloneqq \gamma_{\rm init}\rho^{l}\}$
  \end{algorithmic}
\end{algorithm}

\subsection{Algorithmic Design} \label{sec:algorithmic_design}
To find a stationary point of Problem~\ref{problem:origin},
we present a proximal variable smoothing algorithm in Algorithm~\ref{alg:proposed} as an extension of~\cite[Alg. 1]{Liu-Xia24}.
Although~\cite[Alg. 1]{Liu-Xia24} admits only a special case of Problem~\ref{problem:origin} where
$\mathfrak{S}$
is a surjective linear operator, the proposed algorithm can employ a nonlinear differentiable mapping
$\mathfrak{S}$
(see the paragraph containing~\eqref{eq:bound_surjective_linear} in Section~\ref{sec:measure} for gap between~\cite{Liu-Xia24} and ours).
This subsection focuses on an algorithmic design of the proposed algorithm, and theoretical analyses are deferred to Section~\ref{sec:measure}.

The basic update rule of Algorithm~\ref{alg:proposed} at $n$th iteration is
\begin{equation}
  \bm{x}_{n+1} \coloneqq \prox{\gamma_{n}\phi}(\bm{x}_{n}-\gamma_{n}\nabla \cost_{n}(\bm{x}_{n})) \in \dom{\phi}
\end{equation}
with the current estimate
$\bm{x}_{n} \in \dom{\phi}$
of a stationary point and
a smoothed surrogate function
$\cost_{n}\coloneqq \cost^{\langle\mu_{n}\rangle}=h+\moreau{g}{\mu_{n}}\circ\mathfrak{S}$
of
$\cost\coloneqq h+g\circ\mathfrak{S}$
(see~\eqref{eq:surrogate_function}),
where
$\mu_{n} \in (0,\eta^{-1})$
and
$\gamma_{n} \in \mathbb{R}_{++}$.
Since
$g$
and
$\phi$
in Problem~\ref{problem:origin} are prox-friendly, the proposed Algorithm~\ref{alg:proposed} can be implemented as a single-loop algorithm.

Here, we employ smoothing parameters
$(\mu_{n})_{n=1}^{\infty} \subset (0,(2\eta)^{-1}]$
satisfying the following conditions
\begin{equation}
  \left\{
  \begin{array}{l}
    {\rm(i)}\ \lim_{n\to\infty} \mu_{n} = 0, \quad
    {\rm(ii)}\ \sum_{n=1}^{\infty} \mu_{n} = +\infty, \\
    {\rm(iii)}\ (\forall n\in \mathbb{N})\ \mu_{n+1}\leq \mu_{n},
  \end{array}\right.
  \label{eq:nonsummable}
\end{equation}
which is introduced in~\cite[Eq. (26)]{Kume-Yamada24A} for establishing a subsequential convergence guarantee of a variable smoothing-type algorithm~\cite[Alg. 1]{Kume-Yamada24A} developed for Problem~\ref{problem:origin} with
$\phi\equiv 0$.
For example,
$(\mu_{n})_{n=1}^{\infty}\coloneqq ((2\eta)^{-1}n^{-1/\alpha})$
satisfies~\eqref{eq:nonsummable} for each
$\alpha \geq 1$,
and
$\alpha=3$
is a reasonable choice for achieving a better convergence rate (see Remark~\ref{remark:reasonable_mu}).

For every
$n\in\mathbb{N}$,
we use a stepsize
$\gamma_{n} \in \mathbb{R}_{++}$
satisfying the following {\em Armijo condition} with
$\gamma\coloneqq \gamma_{n}$:
{ %
\thickmuskip=0.3\thickmuskip
\medmuskip=0.3\medmuskip
\thinmuskip=0.3\thinmuskip
\begin{equation}
  \begin{array}{l}
    \arraycolsep=0.3\arraycolsep
    (\cost_{n}+\phi)\left(\prox{\gamma\phi}(\bm{x}_{n} - \gamma \nabla \cost_{n}(\bm{x}_{n}))\right) \\
    \leq (\cost_{n}+\phi)(\bm{x}_{n}) - c\gamma\stdnorm{\frac{\bm{x}_{n}-\prox{\gamma \phi}(\bm{x}_{n}-\gamma \nabla\cost_{n}(\bm{x}_{n}))}{\gamma}}^{2},
  \end{array}\hspace{-1em}
  \label{eq:Armijo}
\end{equation}}%
where a constant
$c \in (0,1)$
is typically chosen as a small value (e.g.,
$10^{-4}$).
A stepsize
$\gamma_{n}$
satisfying~\eqref{eq:Armijo} with
$\gamma \coloneqq \gamma_{n}$
is reasonable because the sufficient decrease
$(\cost_{n}+\phi)(\bm{p}) < (\cost_{n}+\phi)(\bm{x}_{n})$
is implied if
$\bm{x}_{n}\neq  \prox{\gamma_{n} \phi}(\bm{x}_{n}-\gamma_{n} \nabla\cost_{n}(\bm{x}_{n}))$.

  Under Assumption~\ref{assumption:Lipschitz_smoothed}, we have mainly two choices of
  $\gamma_{n}$
  satisfying the Armijo condition~\eqref{eq:Armijo} with
  $\gamma \coloneqq \gamma_{n}$ (see Examples~\ref{example:stepsize:diminishing} and~\ref{example:stepsize:backtracking} below).
  These choices are motivated by the fact that~\eqref{eq:Armijo} is satisfied for every
  $\gamma \in (0,2(1-c)L_{\nabla \cost_{n}}^{-1}]$
  by~\cite[Lemma 10.4]{Beck17} together with the convexity of
  $\phi$
  and
  the $L_{\nabla \cost_{n}}$-Lipschitz continuity of
  $\nabla \cost_{n}$
  over
  $\dom{\phi}$.

  \begin{example}[Diminishing stepsizes]\label{example:stepsize:diminishing}
          $\gamma_{n}\coloneqq 2(1-c)L_{\nabla \cost_{n}}^{-1}$
          is a naive choice of stepsize
          satisfying the Armijo condition~\eqref{eq:Armijo} with
          $\gamma \coloneqq \gamma_{n}$
          in a case where
          $L_{\nabla \cost_{n}}$
          in Assumption~\ref{assumption:Lipschitz_smoothed}
          is known.
          In this case,
          $(\gamma_{n})_{n=1}^{\infty}$
          is diminishing stepsizes because
          $L_{\nabla \cost_{n}}=\varpi_{1}+\varpi_{2}\mu_{n}^{-1}$
          is monotonically increasing due to the monotonically decreasing of
          $\mu_{n}$
          by the condition (iii) in~\eqref{eq:nonsummable}.
  \end{example}

  \begin{example}[Backtracking-based stepsizes]\label{example:stepsize:backtracking}
          The so-called {\em backtracking algorithm} (see~\cite[Sect. 10.3.3]{Beck17})
          can be used for finding
          $\gamma_{n}$
          satisfying the Armijo condition~\eqref{eq:Armijo} with
          $\gamma\coloneqq \gamma_{n}$
          even if
          $L_{\nabla \cost_{n}}$
          in Assumption~\ref{assumption:Lipschitz_smoothed}
          is unknown.
          The backtracking algorithm illustrated in Algorithm~\ref{alg:backtracking} (with
          $c \in (0,1)$,
          $\rho \in (0,1)$,
          and
          $\gamma_{\rm init} \in \mathbb{R}_{++}$)
          repeats to increase
          $l \in \mathbb{N}\cup \{0\}$
          from
          $0$
          until
          $\gamma \coloneqq \gamma_{\rm init}\rho^{l}$
          satisfies the Armijo condition~\eqref{eq:Armijo}.
          Since the Armijo condition~\eqref{eq:Armijo}
          holds for every
          $\gamma \in (0,2(1-c)L_{\nabla \cost_{n}}^{-1}]$~\cite[Lemma 10.4]{Beck17},
          Algorithm~\ref{alg:backtracking} terminates within
          $\max\{0,\lceil \log_{\rho}(2(1-c)L_{\nabla \cost_{n}}^{-1}\gamma_{\rm init}^{-1})\rceil\}$
          iterations at most with the ceiling function
          $\lceil \cdot \rceil$,
          and thus we have
          \begin{equation}
            \gamma_{\rm init} \geq \gamma_{n} \geq \min\left\{\gamma_{\rm init}, 2\rho(1-c)L_{\nabla \cost_{n}}^{-1}\right\}.
            \label{eq:stepsize_lower_bound}
          \end{equation}
  \end{example}

  \begin{remark} \label{remark:stepsize}
    With both stepsizes in Examples~\ref{example:stepsize:diminishing} and~\ref{example:stepsize:backtracking}, a subsequential convergence guarantee of Algorithm~\ref{alg:proposed} is established in Theorem~\ref{theorem:convergence_extension}~\ref{enum:theorem:convergence:subsequential} in Section~\ref{sec:measure}.
    In addition, by employing stepsizes in Example~\ref{example:stepsize:diminishing}, we derive a convergence rate
    $\mathcal{O}(\epsilon^{-3})$
    of an $\epsilon$-stationary point of
    $\cost + \phi$
    (see Theorem~\ref{theorem:convergence_extension}~\ref{enum:theorem:convergence:rate} and Remark~\ref{remark:reasonable_mu}).
    Nevertheless, compared with diminishing stepsizes in Example~\ref{example:stepsize:diminishing},
    backtracking-based stepsizes in Example~\ref{example:stepsize:backtracking} tend to yield faster practical convergence of Alg.~\ref{alg:proposed} (see also Fig.~\ref{fig:maxmin-four-sizes} for numerical comparisons) because Alg.~\ref{alg:backtracking} exploits local information around the current estimate
          $\bm{x}_{n}$
          to get larger stepsizes than Example~\ref{example:stepsize:diminishing}.
  \end{remark}

\subsection{Convergence Analysis via Gradient Mapping} \label{sec:measure}
In this section, we present a convergence analysis of Algorithm~\ref{alg:proposed} in terms of a subsequential convergence guarantee and a convergence rate to a
($\epsilon$-)
stationary point of Problem~\ref{problem:origin}.
For our convergence analysis, we introduce a useful quantity with
$\cost$
and
$\phi$
in Problem~\ref{problem:origin},
and with
$\gamma > 0$:
\begin{equation}
  \mathcal{M}_{\gamma}^{\cost,\phi}
  : \mathcal{X} \to \mathbb{R}
  :\widebar{\bm{x}}
  \mapsto
  \dist\left(\bm{0},\frac{\widebar{\bm{x}} - \prox{\gamma\phi}\left(\widebar{\bm{x}}-\gamma \Lsubdiff \cost\left(\widebar{\bm{x}}\right)\right)}{\gamma}\right),
  \label{eq:measure}
\end{equation}
and its smoothed variant replaced
$\cost$
by
$\cost^{\langle \mu \rangle}$
in~\eqref{eq:surrogate_function}:
\begin{equation}
  \thickmuskip=0.3\thickmuskip
  \medmuskip=0.3\medmuskip
  \thinmuskip=0.3\thinmuskip
  (\bm{x}\in\mathcal{X}) \ \mathcal{M}_{\gamma}^{\cost^{\langle \mu \rangle},\phi}(\bm{x}) = \stdnorm{\frac{\bm{x}-\prox{\gamma \phi}(\bm{x}-\gamma\nabla \cost^{\langle \mu \rangle})(\bm{x})}{\gamma}}, \hspace{-1em}\label{eq:measure_smooth}
\end{equation}
where
$\prox{\gamma\phi}(\widebar{\bm{x}}-\gamma \Lsubdiff \cost(\widebar{\bm{x}}))\coloneqq \{\prox{\gamma\phi}(\widebar{\bm{x}}-\gamma \bm{v}) \in \mathcal{X} \mid \bm{v} \in \Lsubdiff \cost(\widebar{\bm{x}})\} \subset \mathcal{X}$
(see Remark~\ref{remark:stationarity} for a reason why we consider these quantities).
By the update of Algorithm~\ref{alg:proposed}, we have
\begin{equation}
  (n\in\mathbb{N})\quad
  \mathcal{M}_{\gamma_{n}}^{\cost^{\langle \mu_{n}\rangle},\phi}(\bm{x}_{n}) = \stdnorm{\gamma_{n}^{-1}(\bm{x}_{n}-\bm{x}_{n+1})}. \label{eq:residual}
\end{equation}

In a special case where
$\cost$
is continuously differentiable,
the quantity
$\mathcal{M}_{\gamma}^{\cost,\phi}(\widebar{\bm{x}})$
is given as the norm
$\stdnorm{\gamma^{-1}(\widebar{\bm{x}} - \prox{\gamma\phi}(\widebar{\bm{x}} - \gamma \nabla \cost(\widebar{\bm{x}})))}$
of the {\em gradient mapping}~\cite[Sect. 10.3.2]{Beck17},
or {\em forward-backward residual}~\cite{Themelis-Stella-Patrinos18},
due to
$\Lsubdiff \cost = \nabla \cost$~\cite[Exe. 8.8 (b)]{Rockafellar-Wets98}.
Moreover,
$\mathcal{M}_{\gamma}^{\cost,\phi}$
coincides with~\cite[$\dist(0, \mathcal{G}_{\lambda,\mathcal{A}}(\cdot))$ in (11) with
$\lambda\coloneqq\gamma$
and
$\mathcal{A}\coloneqq \mathfrak{S}$]{Liu-Xia24} if
the component
$\mathfrak{S}$
of
$\cost = h+g\circ\mathfrak{S}$
is linear and surjective.

The following lemma states that a stationary point
$\widebar{\bm{x}} \in \dom{\phi}$
of
$\cost + \phi$
can be characterized by
$\mathcal{M}_{\gamma}^{\cost,\phi}(\widebar{\bm{x}}) = 0$
for any
$\gamma > 0$.
Hence, the gradient mapping-type quantity
$\mathcal{M}_{\gamma}^{\cost,\phi}$
can be used as a stationarity measure of Problem~\ref{problem:origin}.

\begin{lemma} \label{lemma:optimality}
  Consider
  $\cost= h+g\circ\mathfrak{S}$
  and
  $\phi \in \Gamma_{0}(\mathcal{X})$
  in Problem~\ref{problem:origin}.
  Then, for
  $\widebar{\bm{x}} \in\mathcal{X}$
  and
  $\gamma \in \mathbb{R}_{++}$,
  we have
  \begin{enumerate}[label=(\alph*)]
    \item \label{enum:lemma:optimality:existence}
          $\mathcal{M}_{\gamma}^{\cost,\phi}(\widebar{\bm{x}}) = \stdnorm{\frac{\widebar{\bm{x}}-\prox{\gamma\phi}(\widebar{\bm{x}}-\gamma\widebar{\bm{v}})}{\gamma}}$
          for some
          $\widebar{\bm{v}} \in \Lsubdiff \cost(\widebar{\bm{x}})$;
    \item \label{enum:lemma:optimality:stationarity}
          $\mathcal{M}_{\gamma}^{\cost,\phi}(\widebar{\bm{x}})=0 \Leftrightarrow$
          $\widebar{\bm{x}}$
          is a stationary point of
          $\cost + \phi$.
  \end{enumerate}
\end{lemma}
\begin{proof}
  See~\ref{appendix:lemma:optimality}.
\end{proof}

In a special case where
$\mathfrak{S}$
is a surjective linear operator
(hence
$\mathfrak{S}\circ\mathfrak{S}^{*}$
is invertible),
the following upper bound 
for
$\bm{x} \in \mathcal{X}$
and
$\mu \in (0,\eta^{-1})$
is derived in~\cite[Lemma 7]{Liu-Xia24}:
\begin{equation}
  \begin{array}{l}
  \mathcal{M}_{\widebar{\gamma}}^{\cost,\phi}(\widebar{\bm{x}})
  \leq
  \mathcal{M}_{\gamma}^{\cost^{\langle \mu \rangle},\phi}(\bm{x}) + \frac{\mu L_{g}}{\sigma_{\min}(\mathfrak{S})}\left(\frac{2}{\widebar{\gamma}} + L_{\nabla h}\right)~\mathrm{with}  \\
  \widebar{\bm{x}} \coloneqq \bm{x}-\mathfrak{S}^{*}\circ(\mathfrak{S}\circ\mathfrak{S}^{*})^{-1}(\mathfrak{S}(\bm{x}) - \prox{\mu g}(\mathfrak{S}(\bm{x}))) \in \mathcal{X}
  \end{array}  \label{eq:bound_surjective_linear}
\end{equation}
for
$\gamma,\widebar{\gamma} \in \mathbb{R}_{++}$
with
$\gamma \leq \widebar{\gamma}$
and the minimum positive singular value
$\sigma_{\min}(\mathfrak{S})$
of
$\mathfrak{S}$.
By exploiting~\eqref{eq:bound_surjective_linear} (together with Lemma~\ref{lemma:optimality}~\ref{enum:lemma:optimality:stationarity})
under the surjectivity and linearity of
$\mathfrak{S}$,
\cite{Liu-Xia24} establishes a convergence analysis of a proximal variable smoothing algorithm, in Alg.~\ref{alg:proposed},
that produces
$(\bm{x}_{n})_{n=1}^{\infty}$
vanishing the RHS of~\eqref{eq:bound_surjective_linear} along
$\bm{x}\coloneqq \bm{x}_{n}$
and
$\mu \coloneqq \mu_{n}$.
In contrast, since
$\mathfrak{S}$
in Problem~\ref{problem:origin} is neither surjective nor linear (hence
$\widebar{\bm{x}}$
in~\eqref{eq:bound_surjective_linear} is not well-defined),
the bound~\eqref{eq:bound_surjective_linear} cannot apply to our convergence analysis of Alg.~\ref{alg:proposed}.
Instead, we derive alternative links between a stationarity of
$\cost + \phi$
and
$\mathcal{M}_{\gamma}^{\cost^{\langle \mu \rangle},\phi}$.

\begin{proposition}
  \label{proposition:stationary_relation}
  Consider
  $\cost= h+g\circ\mathfrak{S}$
  and
  $\phi \in \Gamma_{0}(\mathcal{X})$
  in Problem~\ref{problem:origin} under Assumption~\ref{assumption:Lipschitz_smoothed}.
  Let
  $\widebar{\bm{x}} \in \mathcal{X}$,
  $\mu \in (0,\eta^{-1})$,
  and
  $L_{\nabla \cost^{\langle \mu \rangle}} > 0$
  be a Lipschitz constant of
  $\nabla \cost^{\langle \mu \rangle}$.
  Then, for
  $\gamma \in \mathbb{R}_{++}$
  and
  $\widebar{\bm{p}} = \prox{\gamma \phi}(\widebar{\bm{x}} - \gamma \nabla \cost^{\langle \mu \rangle}(\widebar{\bm{x}}))$,
  we have
    {
      \thickmuskip=0.3\thickmuskip
      \medmuskip=0.3\medmuskip
      \thinmuskip=0.3\thinmuskip
      \arraycolsep=0.3\arraycolsep
      \begin{align}
         & \dist(\bm{0}, \nabla h(\widebar{\bm{p}}) + (\mathrm{D}\mathfrak{S}(\widebar{\bm{p}}))^{*}[\Lsubdiff g(\prox{\mu g}(\mathfrak{S}(\widebar{\bm{p}})))] + \subdiff \phi(\widebar{\bm{p}}))         \\
         & \hspace{10em} \leq (1+\gamma L_{\nabla \cost^{\langle \mu \rangle}}) \mathcal{M}_{\gamma}^{\cost^{\langle \mu \rangle}, \phi}(\widebar{\bm{x}}), \label{eq:approximate_stationarity_inequality} \\
         & \stdnorm{\prox{\mu g}(\mathfrak{S}(\widebar{\bm{p}})) - \mathfrak{S}(\widebar{\bm{p}})} \leq \mu L_{g}. \label{eq:bound_prox}
      \end{align}}%
\end{proposition}
\begin{proof}
  See~\ref{appendix:stationary_relation}.
\end{proof}

\begin{proposition} \label{proposition:measure}
  Consider
  $\cost= h+g\circ\mathfrak{S}$
  and
  $\phi \in \Gamma_{0}(\mathcal{X})$
  in Problem~\ref{problem:origin}.
  Let
  $\widebar{\bm{x}} \in \mathcal{X}$
  and
  $(\bm{x}_{n})_{n=1}^{\infty} \subset\mathcal{X}$
  satisfy
  $\lim_{n\to\infty}\bm{x}_{n}=\widebar{\bm{x}}$,
  and
  $(\mu_{n})_{n=1}^{\infty} \subset (0,\eta^{-1})\searrow 0$.
  Then,
  \begin{equation}
    (\gamma \in \mathbb{R}_{++}) \quad \liminf_{n\to\infty} \mathcal{M}_{\gamma}^{\cost_{n}, \phi}(\bm{x}_{n}) \geq  \mathcal{M}_{\gamma}^{\cost,\phi}(\widebar{\bm{x}}) \label{eq:measure_smoothing_origin}
  \end{equation}
  holds with
  $\cost_{n} \coloneqq \cost^{\langle \mu \rangle} = h+\moreau{g}{\mu_{n}}\circ\mathfrak{S}$,
  where
  $\liminf_{n\to\infty} a_{n} = \sup_{n\in\mathbb{N}}\inf_{k\geq n} a_{k}$
  is the limit infimum of a real sequence
  $(a_{n})_{n=1}^{\infty} \subset \mathbb{R}$.
  Moreover,
  by combining~\eqref{eq:measure_smoothing_origin} with Lemma~\ref{lemma:optimality}~\ref{enum:lemma:optimality:stationarity},
  $\widebar{\bm{x}}$
  is a stationary point of
  $\cost+\phi$
  if
  $\liminf\limits_{n\to\infty} \mathcal{M}_{\gamma}^{\cost_{n},\phi}(\bm{x}_{n}) = 0$
  holds with some
  $\gamma \in \mathbb{R}_{++}$.
\end{proposition}
\begin{proof}
  See~\ref{appendix:lemma:set_convergence}.
\end{proof}

\begin{remark}[Comparison to another stationarity measure] \label{remark:stationarity}
  The norm 
  $\stdnorm{\nabla \moreau{(\cost+\phi)}{\mu}(\cdot)}$
  of gradient of the Moreau envelope of
  $\cost + \phi$
  is used as a stationarity measure~\cite{Davis-Drusvyatskiy-MacPhee18},
  while it requires the evaluation of
  $\prox{\mu}(\cost + \phi)$
  (see~\eqref{eq:Moreau_grad}), which is unavailable in our setting.
  In contrast,
  $\mathcal{M}_{\gamma}^{\cost^{\langle \mu \rangle},\phi}$
  in~\eqref{eq:measure_smooth} is computable, and naturally appears 
  in Alg.~\ref{alg:proposed} (see~\eqref{eq:residual}).
  Prop.~\ref{proposition:stationary_relation} shows that
  $\mathcal{M}_{\gamma}^{\cost^{\langle \mu \rangle},\phi}$
  controls an approximate stationarity (see Def.~\ref{definition:optimality}~\ref{enum:definition:optimality:epsilon}), and Prop.~\ref{proposition:measure} links an asymptotic vanishing of
  $\mathcal{M}_{\gamma}^{\cost^{\langle \mu_{n} \rangle},\phi}(\bm{x}_{n})$
  to the stationarity of
  $\cost + \phi$.
\end{remark}

By leveraging Propositions~\ref{proposition:stationary_relation} and~\ref{proposition:measure}, we present a convergence analysis of Algorithm~\ref{alg:proposed}.

\begin{theorem}[Convergence analysis of Algorithm~\ref{alg:proposed}]\label{theorem:convergence_extension}
  Consider
  $\cost+\phi$
  in Problem~\ref{problem:origin} under Assumption~\ref{assumption:Lipschitz_smoothed}.
  Choose arbitrarily
  $\bm{x}_{1} \in \dom{\phi}$,
  $c \in (0,1)$,
  $\rho \in (0,1)$
  and
  $\gamma_{\rm init}> 0$.
  Let
  $(\bm{x}_{n})_{n=1}^{\infty} \subset \mathcal{X}$
  be generated by Algorithm~\ref{alg:proposed} with
  $(\gamma_{n})_{n=1}^{\infty}$
  satisfying one of the following
  (see Examples~\ref{example:stepsize:diminishing} and~\ref{example:stepsize:backtracking}):
  \begin{enumerate}[label=Case(\Roman*),leftmargin=*,align=left]
    \item \label{enum:case:diminishing}
          $\gamma_{n} \coloneqq 2(1-c)L_{\nabla \cost_{n}}^{-1}\ (n\in\mathbb{N})$
          with the Lipschitz constant
          $L_{\nabla \cost_{n}}$
          of
          $\nabla \cost_{n}$
          with
            $\cost_{n}\coloneqq \cost^{\langle \mu_{n} \rangle}$
            in~\eqref{eq:surrogate_function};
    \item \label{enum:case:backtracking}
          $\gamma_{n}$
          obtained by Alg.~\ref{alg:backtracking} with
          $c$,
          $\rho$,
          and
          $\gamma_{\rm init}$.
  \end{enumerate}
  Then, the following hold with
  $\widebar{\gamma}\coloneqq \sup_{n\in\mathbb{N}}\gamma_{n} \in \mathbb{R}_{++}$:
  \begin{enumerate}[label=(\alph*)]
    \item \label{enum:theorem:convergence:upperbound}
          For any pair
          $(\underline{k}, \overline{k}) \in \mathbb{N}^{2}$
          satisfying
          $\underline{k} \leq \overline{k}$,
          we have
          \begin{equation}
            \hspace{-3em}
            \min_{\underline{k}\leq n \leq \overline{k}}
            \mathcal{M}_{\widebar{\gamma}}^{\cost_{n},\phi}(\bm{x}_{n})
            \leq
            \min_{\underline{k}\leq n \leq \overline{k}}
            \mathcal{M}_{\gamma_{n}}^{\cost_{n},\phi}(\bm{x}_{n})
            \leq
            \sqrt{\frac{\chi}{\sum_{n=\underline{k}}^{\overline{k}}\mu_{n}}}
            \label{eq:convergence_rate}
          \end{equation}
          with\,{\thickmuskip=0.1\thickmuskip
              \medmuskip=0.1\medmuskip
              \thinmuskip=0.1\thinmuskip
              \arraycolsep=0.1\arraycolsep
            $\chi\coloneqq \frac{((\cost_{1}+\phi)(\bm{x}_{1}) - \inf_{\bm{x}\in\mathcal{X}} (\cost_{1}+\phi)(\bm{x}) + \mu_{1}L_{g}^{2})(\varpi_{1} + \eta\varpi_{2})}{c\eta\min\{\gamma_{\rm init}L_{\nabla \cost_{1}}, 2\rho(1-c)\}}>0$}.
    \item
          ${\displaystyle\liminf_{n\to\infty} \mathcal{M}_{\widebar{\gamma}}^{\cost_{n},\phi}(\bm{x}_{n}) = 0}$
          holds.
        \item \label{enum:theorem:convergence:subsequential}
          Let
          $m:\mathbb{N}\to\mathbb{N}$
          be a monotonically increasing function such that
          $m(1) \coloneqq 1$
          and
          $m(l) \coloneqq
            \min \{n \in \mathbb{N} \mid n > m(l-1), \mathcal{M}_{\widebar{\gamma}}^{\cost_{n},\phi}(\bm{x}_{n}) \leq 2^{-l}\}\ (l \geq 2)$.
          Then, every cluster point
          $\widebar{\bm{x}} \in \mathcal{X}$
          of
          $(\bm{x}_{m(l)})_{l=1}^{\infty}$
          is a stationary point of
          $\cost+\phi$.
    \item  \label{enum:theorem:convergence:rate}
          Let
          $\epsilon \in (0,1)$.
          If we employ
          $\mu_{n} \coloneqq \tau n^{-1/\alpha}$
          with
          $\tau \in (0,(2\eta)^{-1}]$
          and
          $\alpha > 1$,
          and~\ref{enum:case:diminishing}
          $\gamma_{n} \coloneqq 2(1-c)L_{\nabla \cost_{n}}^{-1}$,
          then Algorithm~\ref{alg:proposed}
          finds an
          $\epsilon$-stationary point of
          $\cost + \phi$
          within at most
            {
              \thickmuskip=0.1\thickmuskip
              \medmuskip=0.1\medmuskip
              \thinmuskip=0.1\thinmuskip
              \arraycolsep=0.1\arraycolsep
            $\left\lceil\left(\Delta + (\tau L_{g})^{\alpha-1} + 1\right)^{\frac{\alpha}{\alpha-1}} \epsilon^{\min\left\{\frac{-2\alpha}{\alpha-1},-\alpha\right\}}\right\rceil$}
          iterations with
          $\Delta\coloneqq \chi\tau^{-1}(1-\alpha^{-1})(3-2c)^{2}$,
          where
          $\chi \in \mathbb{R}_{++}$
          is given in Theorem~\ref{theorem:convergence_extension}~\ref{enum:theorem:convergence:upperbound}.
  \end{enumerate}
\end{theorem}
\begin{proof}
  See~\ref{appendix:convergence}.
\end{proof}

By Thm.~\ref{theorem:convergence_extension}~\ref{enum:theorem:convergence:rate}, we derive   an upper bound on the number of evaluations of {\em basic operations (ops.)}
for achieving an $\epsilon$-stationary point of
$\cost + \phi$.
By following~\cite[Sect. 6]{Drusvyatskiy-Paquette19}, 
    each of the following is counted as one basic operation:
    one evaluation of
    $h(\bm{x})$,
    $\nabla h(\bm{x})$,
    $\mathfrak{S}(\bm{x})$,
    $g(\bm{z})$,
    $\phi(\bm{x})$,
    matrix-vector multiplication
    $(\mathrm{D}\mathfrak{S}(\bm{x}))^{*}[\bm{w}]$,
    $\prox{\gamma \phi}(\bm{x})$
    or
    $\prox{\mu g}(\bm{z})$
    $(\bm{x}\in \mathcal{X}, \bm{z},\bm{w}\in \mathcal{Z})$.
At
$n$th iteration of Alg.~\ref{alg:proposed}
in the setting of Thm.~\ref{theorem:convergence_extension}~\ref{enum:theorem:convergence:rate},
we only need to evaluate
$\prox{\gamma_{n}\phi}(\bm{x}_{n}-\gamma_{n}\nabla \cost_{n}(\bm{x}_{n}))$.
By~\eqref{eq:gradient_expression} and~\eqref{eq:Moreau_grad},
the evaluation of
$\nabla \cost_{n}(\bm{x}_{n})$
requires $4$ basic ops.:
$\nabla h(\bm{x}_{n})$,
$\bm{z}_{n}\coloneqq\mathfrak{S}(\bm{x}_{n})$,
$\bm{w}_{n}\coloneqq\prox{\mu_{n}g}(\bm{z}_{n})$,
and
$(\mathrm{D}\mathfrak{S}(\bm{x}_{n}))^{*}[\mu_{n}^{-1}(\bm{z}_{n}-\bm{w}_{n})]$.
Together with one evaluation of
$\prox{\gamma_{n}\phi}$,
that of
$\prox{\gamma_{n}\phi}(\bm{x}_{n}-\gamma_{n}\nabla \cost_{n}(\bm{x}_{n}))$
requires $5$ basic ops.
Hence, by Thm.~\ref{theorem:convergence_extension}~\ref{enum:theorem:convergence:rate},
with diminishing stepsizes in Example~\ref{example:stepsize:diminishing} and
$\mu_{n} \coloneqq \tau n^{-1/\alpha}\ (\tau \in (0,(2\eta)^{-1}], \alpha > 1)$,
Alg.~\ref{alg:proposed} finds an
$\epsilon$-stationary point of Problem~\ref{problem:origin}
within
\begin{equation}
  \hspace{-2em}
  5 \times \mathcal{O}(\epsilon^{\min\left\{\frac{-2\alpha}{\alpha-1},-\alpha\right\}})=\mathcal{O}(\epsilon^{\min\left\{\frac{-2\alpha}{\alpha-1},-\alpha\right\}})~\mathrm{basic~ops}.\hspace{-2em} \label{eq:basic_operation}
\end{equation}

\begin{remark}[Convergence rate of Algorithm~\ref{alg:proposed}] \label{remark:reasonable_mu}
  In the setting of Theorem~\ref{theorem:convergence_extension}~\ref{enum:theorem:convergence:rate}, by letting
  $\alpha = 3$
  for
  $\mu_{n}=\tau n^{-1/\alpha}$,
  Algorithm~\ref{alg:proposed}
  finds an $\epsilon$-stationary point of
  $\cost+\phi$
  within
  $\mathcal{O}(\epsilon^{-3})$
  basic operations by~\eqref{eq:basic_operation}.
  Hence, Algorithm~\ref{alg:proposed} possesses the rate
  $\mathcal{O}(\epsilon^{-3})$,
  matching the best-known convergence rate of existing algorithms for Problem~\ref{problem:origin}, e.g.,~\cite{Bourkhissi-Necoara-Patrinos-TranDinh25,Drusvyatskiy-Paquette19} (see Table~\ref{table:summary}), where any improvement of the known convergence rate is not our aim.
  The significance of this result is that Algorithm~\ref{alg:proposed} retains this best-known rate while keeping a single-loop structure and, at the same time, providing a subsequential convergence guarantee to a stationary point of Problem~\ref{problem:origin} (see Theorem~\ref{theorem:convergence_extension}~\ref{enum:theorem:convergence:subsequential}).
  In contrast, a convergence rate analysis of Algorithm~\ref{alg:proposed} with
  $\gamma_{n}$
  obtained by Algorithm~\ref{alg:backtracking} has not been established in this paper; hence is left for future work,  whereas its subsequential convergence result is established in Theorem~\ref{theorem:convergence_extension}~\ref{enum:theorem:convergence:subsequential}.
  We note that the choice
  $\mu_{n}=\tau n^{-1/3}$
  is also commonly used in the literature of variable smoothing~\cite{Bohm-Wright21,Liu-Xia24,Kume-Yamada24A,Lopes-Peres-Bilches25}.
\end{remark}

\section{Selected Applications of Problem~\ref{problem:origin} and Algorithm~\ref{alg:proposed}} \label{sec:application}
\subsection{Finite-Max Nonlinear Minimization} \label{sec:application:finite_max}
In this subsection, we consider the following problem
\begin{align}
  \mathrm{find}\ \bm{x}^{\star} & \in \argmin_{\bm{x} \in C}\max_{j=1,2,\ldots,m} \mathfrak{S}_{j}(\bm{x}), \label{eq:max_problem}
\end{align}
where
$C\subset \mathcal{X}$
is a compact convex set such that
$P_{C}$
is available as a computable tool
and
each
$\mathfrak{S}_{j}:\mathcal{X}\to\mathbb{R}$
is twice continuously differentiable.
Nonlinear mappings
$\mathfrak{S}_{j}$
are designed according to specific applications, e.g.,
(i)
$\mathfrak{S}_{j}(\bm{x}) \coloneqq - w_{j}\stdnorm{\bm{x}-\bm{u}_{j}}^{2}$
with a weight
$w_{j} > 0$
and an anchor point
$\bm{u}_{j} \in \mathcal{X}$
for maxmin dispersion problem~\cite{Haines-Loeppky-Tseng-Wang13,Pan-Shen-Xu21,Lopes-Peres-Bilches25};
and
(ii)
$\mathfrak{S}_{j}(\bm{x}) \coloneqq (y_{j}^{2}-\stdnorm{\bm{x}-\bm{u}_{j}}^{2})^{2}$
with a measurement
$y_{j} \in \mathbb{R}$
and an anchor point
$\bm{u}_{j} \in \mathcal{X}$
for robust target localization~\cite{Domingos-Xavier24,Domingos-Xavier25}
(strictly speaking,
$\abs{y_{j}-\stdnorm{\bm{x}-\bm{u}_{j}}}$
has been used therein).

The problem~\eqref{eq:max_problem} can be reformulated into Problem~\ref{problem:origin} with
and
$\mathcal{Z}\coloneqq \mathbb{R}^{m}$
and with
  {
    \thickmuskip=0.3\thickmuskip
    \medmuskip=0.3\medmuskip
    \thinmuskip=0.3\thinmuskip
    \arraycolsep=0.3\arraycolsep
    \begin{align}
       & h\equiv 0, \quad
       g:\mathcal{Z} \to \mathbb{R}: \bm{z} \mapsto \max\{[\bm{z}]_{1},[\bm{z}]_{2},\ldots,[\bm{z}]_{m}\}, \label{eq:max_problem_reformulated_g}                                                                                                                                                                                                       \\
        & \mathfrak{S}:\mathcal{X} \to \mathcal{Z}: \widebar{\bm{x}} \mapsto [
        \mathfrak{S}_{1}(\widebar{\bm{x}})\ \mathfrak{S}_{2}(\widebar{\bm{x}})\ \cdots \ \mathfrak{S}_{m}(\widebar{\bm{x}})]^{\TT}, \ \phi\coloneqq \iota_{C}, \label{eq:max_problem_reformulated}
    \end{align}}%
where
$[\bm{z}]_{l}$
denotes the $l$th entry of
$\bm{z}$,
and
$g$
and
$\phi$
are prox-friendly (see~\cite[Exm. 24.25]{Bauschke-Combettes17} for $g$).
We note that
$\mathrm{D}\mathfrak{S}$
is
$L_{\mathrm{D}\mathfrak{S}}$-Lipschitz continuous over
$\dom{\phi} = \dom{\iota_{C}} = C$
and
Assumption~\ref{assumption:Lipschitz_smoothed} is satisfied for this reformulation
due to the twice continuous differentiability of
$\mathfrak{S}$
over the compact convex set
$C$
(see Example~\ref{example:assumption:Lipschitz_smoothed}),
where a Lipschitz constant
$L_{\mathrm{D}\mathfrak{S}} > 0$
is given by an upper bound of the operator norm of the second derivative
of
$\mathfrak{S}$
over
$C$~\cite[p.288]{Pugh15}.
Since
$g$
is $1$-Lipschitz continuous, we have
$\varpi_{1} = L_{\mathrm{D}\mathfrak{S}}$
and
$\varpi_{2} = \kappa_{\mathfrak{S}}^{2}$
in Assumption~\ref{assumption:Lipschitz_smoothed} with
$\kappa_{\mathfrak{S}} = \max_{\bm{x}\in C} \opnorm{\mathrm{D}\mathfrak{S}(\bm{x})}$.
In particular, 
we get
$L_{\mathrm{D}\mathfrak{S}} = (\sum_{j=1}^{m}L_{\nabla\mathfrak{S}_{j}}^{2})^{1/2}$
and
$\kappa_{\mathfrak{S}} = (\sum_{j=1}^{m} \kappa_{\mathfrak{S}_{j}}^{2})^{1/2}$
with a Lipschitz constant
$L_{\nabla \mathfrak{S}_{j}}$
of
$\nabla \mathfrak{S}_{j}$
over
$\dom{\phi}$
and
$\kappa_{\mathfrak{S}_{j}} \coloneqq \max_{\bm{x}\in C} \stdnorm{\nabla \mathfrak{S}_{j}(\bm{x})}$.
Lemma~\ref{lemma:application_constant} below presents explicitly these constants in the scenario (ii) for which we will test numerical performance in Section~\ref{sec:numerical:dispersion}.
Hence, a stationary point of the problem~\eqref{eq:max_problem} can be approximated iteratively by Alg.~\ref{alg:proposed} (see Thm.~\ref{theorem:convergence_extension}).

\begin{lemma} \label{lemma:application_constant}
  Consider the problem~\eqref{eq:max_problem} with
  $\mathfrak{S}_{j}(\bm{x}) \coloneqq (y_{j}^{2}-\stdnorm{\bm{x}-\bm{u}_{j}}^{2})^{2}$
  ($j=1,2,\ldots,m$,
  $y_{j} \in \mathbb{R}$,
  $\bm{u}_{j} \in \mathcal{X}$)
  and with
  $C \coloneqq B(\bm{0}, r)\coloneqq \{\bm{x}\in\mathcal{X}\mid \stdnorm{\bm{x}}\leq r\}$
  for some
  $r > 0$.
  For each
  $j \in \{1,2,\ldots,m\}$,
  let
  $q_{j,{\rm low}}\coloneqq \max\{0, \stdnorm{\bm{u}_{j}} - r\}$
  and
  $q_{j,{\rm up}} \coloneqq  \stdnorm{\bm{u}_{j}} + r$.
  Then, 
  we get a Lipschitz constant
  $L_{\nabla\mathfrak{S}_{j}}$
  of
  $\nabla\mathfrak{S}_{j}$
  over
  $C$,
  and an upper bound of
  $\kappa_{\mathfrak{S}_{j}} \coloneqq \max_{\bm{x}\in C} \stdnorm{\nabla \mathfrak{S}_{j}(\bm{x})}$
  as:
  \begin{equation}
  L_{\nabla\mathfrak{S}_{j}}  \coloneqq 4\max_{q=q_{j,{\rm low}},q_{j,{\rm up}}}\Upsilon_{j}^{(1)}(q), \quad
  \kappa_{\mathfrak{S}_{j}}
\leq
4\max_{q \in Q_{j}} \Upsilon_{j}^{(2)}(q),
  \end{equation}
  with
  $\Upsilon_{j}^{(1)}(q)\coloneqq \max\{\abs{3q^{2} - y_{j}^{2}}, \abs{q^{2}-y_{j}^{2}}\},
  \ \Upsilon_{j}^{(2)}(q)\coloneqq q\abs{y_{j}^{2} - q^{2}}$,
  and
  $Q_{j} \coloneqq \{\abs{y_{j}}/\sqrt{3},\ q_{j,{\rm up}} \}$
  if
  $\abs{y_{j}}/\sqrt{3} \in [q_{j,{\rm low}}, q_{j,{\rm up}}]$;
  $Q_{j} \coloneqq \{q_{j,{\rm low}},\ q_{j,{\rm up}}\}$
  otherwise.
  Hence,
  Assumption~\ref{assumption:Lipschitz_smoothed} is satisfied by Example~\ref{example:assumption:Lipschitz_smoothed} (see also the paragraph just before Lemma~\ref{lemma:application_constant}).
  Moreover, the function
  $\cost(\bm{x}) = \max\limits_{j=1,2,\ldots,m}\mathfrak{S}_{j}(\bm{x})$
  is
  $\widetilde{\eta}$-weakly convex with
  $\widetilde{\eta}\coloneqq 4\max\limits_{j=1,2,\ldots,m} y_{j}^{2}$.
\end{lemma}
\begin{proof}
  See Section S.1 in the supplementary material.
\end{proof}

\subsection{MU-MIMO Signal Detection} \label{sec:application:mimo}
A massive multiuser multiple-input-multiple-output (MU-MIMO) signal detection problem with
$M(\in \mathbb{N})$-ary
phase-shift keying (PSK)
(see, e.g.,~\cite{Yang-Hanzo15,Chen19,Hayakawa-Hayashi20,Shoji-Yata-Kume-Yamada25}) is formulated as
\begin{equation}
  \mathrm{find}\ \mathsf{s}^{\star} \in \mathsf{D} \subset \mathbb{C}^{U}
  \ \mathrm{from}\
  \mathsf{y} = \mathsf{H}\mathsf{s}^{\star} + \mathsf{e} \in \mathbb{C}^{B}, \label{eq:MIMO_model}
\end{equation}
where
$\mathbb{C}$
denotes the set of all complex numbers
($i$ stands for the imaginary unit, i.e., $i^{2}=-1$ and
$\mathfrak{R}$
and
$\mathfrak{I}$
stand respectively for the real and imaginary parts).
Here,
$\mathsf{s}^{\star} \in \mathsf{D}\coloneqq \{ \exp(i2\pi m/M) \in \mathbb{C} \mid m = 0, 1,\ldots,M-1\}^{U} \subset \mathbb{C}^{U}$
is the transmitted signal with
a discrete set
$\mathsf{D}$
called
a {\em constellation set},
$\mathsf{y} \in \mathbb{C}^{B}$
is the received signal with a known measurement matrix
$\mathsf{H} \in \mathbb{C}^{B\times U}$
and noise
$\mathsf{e} \in \mathbb{C}^{B}$,
and
$U,B\in \mathbb{N}$
are respectively the numbers of transmitting and receiving antennas.
For the problem~\eqref{eq:MIMO_problem},
we review existing approaches, and then propose a new model based on Problem~\ref{problem:origin}.

\subsubsection{Existing Approaches to MIMO Signal Detection} \label{sec:existing_mimo}
The problem~\eqref{eq:MIMO_model} has been tackled via real-valued optimization problems~\cite{Yang-Hanzo15,Chen19,Hayakawa-Hayashi20,Shoji-Yata-Kume-Yamada25} as
\begin{equation}
  \mathop{\mathrm{minimize}}\limits_{\bm{s} \in C \subset \mathbb{R}^{2U}}\ \frac{1}{2}\stdnorm{\bm{y}-\bm{H}\bm{s}}^{2} + \lambda \psi(\bm{s})
  \label{eq:MIMO_problem}
\end{equation}
with the real-valued expressions of
$\mathsf{y}$
and
$\mathsf{H}$:
\begin{equation}
  \bm{y}\coloneqq \widehat{\textnormal{\ensuremath{\mathsf{y}}}}\coloneqq
  \begin{bmatrix} \mathfrak{R}(\mathsf{y})\\ \mathfrak{I}(\mathsf{y}) \end{bmatrix} \in \mathbb{R}^{2B},\;
  \bm{H} \coloneqq \begin{bmatrix} \mathfrak{R}(\mathsf{H}) & - \mathfrak{I}(\mathsf{H}) \\ \mathfrak{I}(\mathsf{H}) & \mathfrak{R}(\mathsf{H}) \end{bmatrix} \in \mathbb{R}^{2B\times 2U},
  \label{eq:real_valued}
\end{equation}
where a constraint set
$(\emptyset \neq )C \subset \mathbb{R}^{2U}$
is designed  as a superset of the real-valued expression
$\mathfrak{D}\coloneqq \{\widehat{\mathsf{s}} \in \mathbb{R}^{2U} \mid \mathsf{s} \in \mathsf{D}\} \subset \mathbb{R}^{2U}$
of
$\mathsf{D}$,
$\lambda \in \mathbb{R}_{++}$
is a weight
and
$\psi:\mathbb{R}^{2U}\to \mathbb{R}$
is a regularizer.
We briefly review below three existing models based on the model~\eqref{eq:MIMO_problem} for the problem~\eqref{eq:MIMO_model}.

A classical {\em linear minimum mean-square-error estimation (LMMSE)}, e.g.,~\cite{Yang-Hanzo15}, can be seen as an estimation via the model
$\mathop{\mathrm{minimize}}\limits_{\bm{s} \in \mathbb{R}^{2U}}\ \frac{1}{2}\stdnorm{\bm{y}-\bm{H}\bm{s}}^{2} + \frac{\sigma^{2}}{2}\stdnorm{\bm{s}}^{2}, \label{eq:LMMSE}$
where
$\sigma^{2}$
stands for
the variance of noise
$\mathsf{e}$
in~\eqref{eq:MIMO_model}.
Its unique minimizer is explicitly given by
$(\bm{H}^{\TT}\bm{H}+\sigma^{2}\bm{I})^{-1}\bm{H}^{\TT}\bm{y}$.

{\em Modulus-constrained least squares model (Modulus)}~\cite{Chen19} is given with a constraint 
$\mathbb{T}\coloneqq \{z \in \mathbb{C} \mid \abs{z} = 1\} \subset \mathbb{C}$
by
\begin{equation}
  \mathop{\mathrm{minimize}}\limits_{\bm{s} \in \{\widehat{\mathsf{z}} \in \mathbb{R}^{2U} \mid \mathsf{z} \in \mathbb{T}^{U}\}}\ \frac{1}{2}\stdnorm{\bm{y}-\bm{H}\bm{s}}^{2}, \label{eq:modulus}
\end{equation}
where
$\mathbb{T}$
is
a superset of the constellation set
$\mathsf{D}$.
We note that
\cite{Chen19} addressed directly the complex version of the model~\eqref{eq:modulus}.

{\em Sum-of-absolute-values model (SOAV)}~\cite{Hayakawa-Hayashi20} is given by
\begin{equation}
  \mathop{\mathrm{minimize}}\limits_{\bm{s} \in \mathrm{Conv}(\mathfrak{D})}\ \frac{1}{2}\stdnorm{\bm{y}-\bm{H}\bm{s}}^{2} + \lambda\psi_{\rm SOAV}(\bm{s}), \label{eq:SOAV}
\end{equation}
where
$\mathrm{Conv}(\cdot)$
denotes the convex hull of a given set.
The so-called SOAV function
$\psi_{\rm SOAV}\coloneqq\frac{1}{M}\sum_{m=1}^{M} \norm{\cdot- \widehat{\mathsf{s}}_{m}}_{1}$~\cite{Nagahara15}
is designed to penalize a given
$\bm{s} \in \mathbb{R}^{2U}$
that deviates from discrete-valued points in
$\mathfrak{D}$,
where
$\mathsf{s}_{m} \coloneqq \exp(i2\pi m/M)\bm{1} \in \mathsf{D}$
is given with the all one vectors
$\bm{1}$,
and
$\widehat{\mathsf{s}}_{m} \in \mathfrak{D}\ (m=0,1,\ldots,M-1)$
(see~\eqref{eq:real_valued} for $\widehat{(\cdot)}$).
However, the penalization by
$\psi_{\rm SOAV}$
is not contrastive enough to distinguish between points in
$\mathfrak{D}$
and the other points because any point in
$\mathfrak{D}$
is never an isolated (local) minimizer of
$\psi_{\rm SOAV}$
due to its convexity.
Hence, in order to distinguish them, some nonconvex regularizer
$\psi$
is required.
This requirement motivates us to propose a new nonconvex regularizer for~\eqref{eq:MIMO_model} 
as follwos.

\subsubsection{Proposed Approach to MIMO Signal Detection}
For the problem~\eqref{eq:MIMO_model},
we propose the following model with a more contrastive regularizer
$\psi_{\lambda_{r},\lambda_{\theta}}$
in~\eqref{eq:polar_regularizer} below than
$\psi_{\rm SOAV}$:
\begin{equation}
  \mathop{\mathrm{minimize}}\limits_{(\bm{r},\bm{\theta}) \in [\underline{r},1]^{U}\times\mathbb{R}^{U}}
  \frac{1}{2}\stdnorm{\bm{y}-\bm{H}\param(\bm{r},\bm{\theta})}^{2}
  + \psi_{\lambda_{r},\lambda_{\theta}}(\bm{r},\bm{\theta}),
  \label{eq:MIMO_proposed}
\end{equation}
where we used a
polar coordinate-type expression of
$\bm{s} \in \mathbb{R}^{2U}$:
\begin{equation} \label{eq:polar}
  \param:\mathbb{R}^{U}\times \mathbb{R}^{U} \to \mathbb{R}^{2U}:(\bm{r},\bm{\theta}) \mapsto \bm{s}\coloneqq
  \begin{bmatrix}
    \bm{r}\odot \mathrm{\bf cos}(\bm{\theta}) \\ \bm{r} \odot \mathrm{\bf sin}(\bm{\theta})
  \end{bmatrix},
\end{equation}
i.e.,
$\bm{r}$
and
$\bm{\theta}$
denote the amplitude and phase of
$\bm{s}$
respectively,
\begin{equation}
  \thickmuskip=0.0\thickmuskip
  \medmuskip=0.0\medmuskip
  \thinmuskip=0.0\thinmuskip
  \arraycolsep=0.5\arraycolsep
  \mathrm{\bf sin}:\mathbb{R}^{U}\to \mathbb{R}^{U}:\bm{\theta}\mapsto \begin{bmatrix} \sin([\bm{\theta}]_{1}) & \sin([\bm{\theta}]_{2}) & \cdots & \sin([\bm{\theta}]_{U}) \end{bmatrix}^{\TT} \label{eq:sin}
\end{equation}
(so as $\mathrm{\bf cos}:\mathbb{R}^{U}\to\mathbb{R}^{U}$),
and
$\underline{r} \in (0,1]$
is a lower bound
for each
$[\bm{r}]_{u}$
inspired by~\cite{Ramírez-Santamaría06}.
Here, $\odot$
is the entry-wise (Hadamard) product.
The proposed regularizer
$\psi_{\lambda_{r},\lambda_{\theta}}$
is
\begin{align}
  \psi_{\lambda_{r},\lambda_{\theta}}
   & :[\underline{r},1]^{U}\times \mathbb{R}^{U}\to \mathbb{R}: \\
   & :(\bm{r},\bm{\theta}) \mapsto
  \lambda_{r}\sum_{u=1}^{U} [\bm{r}]_{u}^{-1}
  + \lambda_{\theta}\norm{\mathrm{\bf sin}\left(\frac{M\bm{\theta}}{2}\right)}_{1}
  \label{eq:polar_regularizer}
\end{align}
with weights
$\lambda_{r},\lambda_{\theta} \in \mathbb{R}_{+}$.
For
$\lambda_{r},\lambda_{\theta} \in \mathbb{R}_{++}$,
we have
\begin{align}
  \hspace{-1em}
  \param(\bm{r}^{\star},\bm{\theta}^{\star}) \in \mathfrak{D} \Leftrightarrow (\bm{r}^{\star},\bm{\theta}^{\star}) \in &
  \argmin_{(\bm{r},\bm{\theta}) \in [\underline{r},1]^{U}\times \mathbb{R}^{U}} \psi_{\lambda_{r},\lambda_{\theta}}(\bm{r},\bm{\theta}) \label{eq:desired}
\end{align}
by
$\argmin_{(\bm{r},\bm{\theta}) \in [\underline{r},1]^{U}\times \mathbb{R}^{U}} \psi_{\lambda_{r},\lambda_{\theta}}(\bm{r},\bm{\theta}) =\{(\bm{1}, 2\pi \bm{m}/M) \mid \bm{m} \in \mathbb{Z}^{U}\}$,
where
$\mathbb{Z}$
denotes the set of all integers.
Hence, compared to the SOAV function
$\psi_{\rm SOAV}$~\cite{Nagahara15},
$\psi_{\lambda_{r},\lambda_{\theta}}$
is a desirable regularizer
to penalize a point
$\bm{s} = \param(\bm{r},\bm{\theta})$
that deviates from points in the real-valued expression
$\mathfrak{D} = \{\widehat{\mathsf{s}} \in \mathbb{R}^{2U}\mid \mathsf{s} \in \mathsf{D}\} \subset \mathbb{R}^{2U}$
of the constellation set
$\mathsf{D}$.
Such a desirable property in~\eqref{eq:desired} holds even if
$\norm{\cdot}_{1}$
in
$\psi_{\lambda_{r},\lambda_{\theta}}$
is replaced by, e.g.,
MCP~\cite{Zhang10}, and SCAD~\cite{Fan-Li01}.
To examine a potential of the basic idea of the proposed regularizer, we consider the simplest case
$\norm{\cdot}_{1}$
in this paper.

The proposed model~\eqref{eq:MIMO_proposed} can be reformulated into Problem~\ref{problem:origin}  with
$\mathcal{X}\coloneqq \mathbb{R}^{U}\times \mathbb{R}^{U}$,
$\mathcal{Z} \coloneqq \mathbb{R}^{U}$,
and
\vspace{-0.5em}
{
  \thickmuskip=0.1\thickmuskip
  \medmuskip=0.1\medmuskip
  \thinmuskip=0.1\thinmuskip
  \arraycolsep=0.1\arraycolsep
  \begin{align}
     & \hspace{-1em} h:\mathcal{X}\to\mathbb{R}: (\bm{r},\bm{\theta}) \mapsto \frac{1}{2}\stdnorm{\bm{y}-\bm{H}\param(\bm{r},\bm{\theta})}^{2} + \lambda_{r}\sum_{u=1}^{U} \mathfrak{d}([\bm{r}]_{u}), \label{eq:mimo_h} \\
     & \hspace{-1em} g:\mathcal{Z}\to\mathbb{R}:\bm{z} \mapsto \lambda_{\theta}\norm{\bm{z}}_{1}, \
    \mathfrak{S}:\mathcal{X}\to\mathcal{Z}:(\bm{r},\bm{\theta}) \mapsto \mathrm{\bf sin}\left(\frac{M\bm{\theta}}{2}\right), \label{eq:mimo_S}                                                                        \\
     & \hspace{-1em} \phi:\mathcal{X}\to\exR:(\bm{r},\bm{\theta}) \mapsto \iota_{C}(\bm{r},\bm{\theta})\ 
      \mathrm{with}\ 
C\coloneqq [\underline{r},1]^{U} \times \mathbb{R}^{U},
      \label{eq:mimo_phi}
  \end{align}}%
where
$\mathfrak{d}:\mathbb{R}\to\mathbb{R}$
is a differentiable function\footnote{
The function
$t\mapsto t^{-1}$
in the first term of~\eqref{eq:polar_regularizer} is replaced by a function
$\mathfrak{d}$
so that
(I)
$h$
in~\eqref{eq:mimo_h} satisfies the condition~\ref{enum:problem:origin:h} in Problem~\ref{problem:origin};
(II)~$\frac{1}{2}\stdnorm{\bm{y}-\bm{H}\param(\cdot)}^{2} + \psi_{\lambda_{r},\lambda_{\theta}}(\cdot)$
in~\eqref{eq:MIMO_proposed}
and
$h+g\circ\mathfrak{S}$
in~\eqref{eq:mimo_h} and~\eqref{eq:mimo_S}
take the same value on the constraint
$C$.
For example, we can employ
$\mathfrak{d}(t)\coloneqq t^{-1}$
for
$t \in [\underline{r},+\infty)$,
and
$\mathfrak{d}(t)\coloneqq -\underline{r}^{-1}t + 1 + \underline{r}^{-1}$
otherwise.
} such that
$\mathfrak{d}(t) = \frac{1}{t}$
for
$t \in [\underline{r},1]$
and
$\nabla \mathfrak{d}$
is Lipschitz continuous over
$\mathbb{R}$,
and
$g$
and
$\phi$
are prox-friendly (see, e.g.,~\cite[Exm. 24.22]{Bauschke-Combettes17} and~\cite[Lem. 6.26]{Beck17} respectively).
Moreover,
Assumption~\ref{assumption:Lipschitz_smoothed} is satisfied for this reformulation by Lemma~\ref{lemma:application_mimo_constant}.
Hence, a stationary point of the problem~\eqref{eq:MIMO_proposed} can be approximated iteratively to any precision by Algorithm~\ref{alg:proposed} (see Theorem~\ref{theorem:convergence_extension}).
\begin{lemma} \label{lemma:application_mimo_constant}
  The minimization of
  $\cost+\phi=h + g\circ\mathfrak{S} + \phi$
  with~\eqref{eq:mimo_h}-\eqref{eq:mimo_phi}
  is an instance of Problem~\ref{problem:origin}, where
  $\lambda_{\theta},\lambda_{r} \in \mathbb{R}_{+}$
  and
  $\underline{r} \in (0,1]$.
  Moreover, we have
  $\sup_{(\bm{r},\bm{\theta}) \in C}\opnorm{\mathrm{D}\mathfrak{S}(\bm{r},\bm{\theta})} < +\infty$,
  implying
  Assumption~\ref{assumption:Lipschitz_smoothed} is satisfied by Example~\ref{example:assumption:Lipschitz_smoothed}.
\end{lemma}
\begin{proof}
  See Section S.2 in the supplementary material.
\end{proof}

\section{Numerical experiments}\label{sec:numerical}
We demonstrate the effectiveness of Alg.~\ref{alg:proposed} by numerical experiments in two scenarios introduced in Section~\ref{sec:application} performed on MacBook Pro (14-inch, M3, Nov. 2023, 16 GB).
Unless stated otherwise, all experiments were carried by MATLAB.

{\bf Setting of Alg.~\ref{alg:proposed} and~\ref{alg:backtracking}}:
For Algorithm~\ref{alg:proposed}, we employed
$(\mu_{n})_{n=1}^{\infty}\coloneqq (\tau n^{-1/3})_{n=1}^{\infty}$
with
$\tau \in (0,(2\eta)^{-1}]$
(see Remark~\ref{remark:reasonable_mu})
and
$c = 2^{-13}$.
For
$(\gamma_{n})_{n=1}^{\infty}$,
we employed diminishing stepsizes:
$\gamma_{n}\coloneqq 2(1-c)L_{\nabla \cost_{n}}^{-1}=2(1-c)(\varpi_{1} + \tau^{-1}\varpi_{2}n^{1/3})^{-1}$
with
$L_{\nabla \cost_{n}} = \varpi_{1} + \varpi_{2}\mu_{n}^{-1}$
(see Example~\ref{example:stepsize:diminishing}),
or
backtracking-based stepsizes obtained by Algorithm~\ref{alg:backtracking} with
$c$,
$\rho \coloneqq 2^{-1}$
and
$\gamma_{\rm init} \coloneqq 1$
(see Example~\ref{example:stepsize:backtracking}).

\subsection{Application to robust target localization} \label{sec:numerical:dispersion}
We tested the convergence performance of Alg.~\ref{alg:proposed} in a scenario of
robust target localization~\cite{Domingos-Xavier24,Domingos-Xavier25} in~\eqref{eq:max_problem}, where
$C \coloneqq B(\bm{0}, 1)= \{\bm{x}\in\mathcal{X}\mid \stdnorm{\bm{x}}\leq 1\}$
and
$\mathfrak{S}_{j}(\bm{x}) \coloneqq (y_{j}^{2}-\stdnorm{\bm{x}-\bm{u}_{j}}^{2})^{2}\ (j=1,2,\ldots,m)$
was given with
$y_{j} \in \mathbb{R}$
and
$\bm{u}_{j} \in \mathcal{X}$
(see Section~\ref{sec:application:finite_max}).
By letting
$\mathcal{X} \coloneqq \mathbb{R}^{d}$,
for each
$j=1,2,\ldots,m$,
we generated
$\bm{u}_{j} \in \mathbb{R}^{d}$
uniformly from
$[-1,1]^{d}$
and
$y_{j} \coloneqq \stdnorm{\bm{x}^{\star}-\bm{u}_{j}} \in \mathbb{R}$,
where
$\bm{x}^{\star} \in C$
was generated by
$P_{C}(\bm{x}^{\diamond}) \in C$
with
$\bm{x}^{\diamond} \in \mathbb{R}^{d}$
sampled from the uniform distribution
$[-1,1]^{d}$.
With this construction, we have
$\mathfrak{S}_{j}(\bm{x}^{\star}) = 0\ (j=1,2,\ldots,m)$,
from which
$\bm{x}^{\star}$
is a minimizer of the target problem,
and
$\min_{\bm{x}\in \mathcal{X}}(\cost+\phi)(\bm{x}) = 0$.
Hence, the quantity
$(\cost+\phi)(\bm{x}_{n})$
can be used for evaluating the quality of estimated solution
$\bm{x}_{n} \in \mathcal{X}$
by algorithms.

To examine behaviors of Alg.~\ref{alg:proposed}, we employed different settings:
(i)
$\tau \in \{1, 10^{5}\}$
in
$(\mu_{n})_{n=1}^{\infty}= (\tau n^{-1/3})_{n=1}^{\infty}$;
(ii)
diminishing stepsizes
$\gamma_{n}\coloneqq 2(1-c)L_{\nabla \cost_{n}}^{-1}$,
or
stepsizes obtained by the backtracking algorithm in Alg.~\ref{alg:backtracking}.
For comparisons, we employed three algorithms:
(I) smoothing proximal linear method (PLM)~\cite[Alg. 6]{Drusvyatskiy-Paquette19} with a Moreau envelope-based smoothing of
$g$
in the reformulation~\eqref{eq:max_problem_reformulated_g}
(see~\cite[Sect. 6]{Drusvyatskiy-Paquette19}),
(II) augmented Lagrangian method (ALM)~\cite[Alg. 1]{DeMarchi-Jia-Kanzow-Mehlitz23},
and
(III) subgradient method (SGM)~\cite[Eq. (2.4)]{Li-Zhu-So-Man-Cho-Lee19}.
We implemented PLM by MATLAB, where a fast gradient method~\cite[Alg. 5]{Drusvyatskiy-Paquette19} was employed for subproblems in PLM, and the inner loop at
$n$th outer iteration
was terminated when the number of inner iterations exceeded 
$\min\{T_{n}, 10^{4}\}$,
where
$T_{n}\in \mathbb{N}$
is given in~\cite[Ep. (6.12)]{Drusvyatskiy-Paquette19}.
For PLM,
we employed a tolerance
$\epsilon \coloneqq 10^{-2}$
and a smoothing parameter
$\nu\coloneqq \epsilon^{2}/(2L_{g}^{3}L_{\mathrm{D}\mathfrak{S}}) =\epsilon^{2}/(2L_{\mathrm{D}\mathfrak{S}})$~\cite[p.534]{Drusvyatskiy-Paquette19},
where
$L_{\mathrm{D}\mathfrak{S}} = (\sum_{j=1}^{m}L_{\nabla \mathfrak{S}_{j}}^{2})^{1/2}$
is given with
$L_{\nabla \mathfrak{S}_{j}}$
in Lemma~\ref{lemma:application_constant}
(Note: smaller values of
$\epsilon$
led to numerical instability in our implementation).
For ALM, we used the author's open code\footnote{
  Code for ALM~\cite{DeMarchi-Jia-Kanzow-Mehlitz23} is available in \url{https://github.com/aldma/Bazinga.jl}.
} implemented by Julia, where a proximal gradient-type algorithm~\cite[Alg. 2]{DeMarchi-Jia-Kanzow-Mehlitz23} was employed for subproblems in ALM with
default setting for parameters (see~\cite[Sect. 4.1]{DeMarchi-Jia-Kanzow-Mehlitz23}).
We also implemented SGM~\cite[Thm. 1 (b)]{Li-Zhu-So-Man-Cho-Lee19} by MATLAB:
$\bm{x}_{n+1} \in P_{C}(\bm{x}_{n}-\gamma_{n}\bm{v}_{n})$
with
$\gamma_{n} \coloneqq \frac{1}{\widetilde{\eta}\sqrt{n}}$
and
$\bm{v}_{n} \in \Lsubdiff \cost(\bm{x}_{n})$,
where
$\widetilde{\eta}$
is given in Lemma~\ref{lemma:application_constant}.
For each algorithm, we set
$\bm{x}_{1} \coloneqq \bm{0} \in C$.
Each algorithm was terminated when
$\bm{x}_{n}$
achieved
$(\cost+\phi)(\bm{x}_{n})< 10^{-10}$,
or running CPU time exceeded
$100$ seconds.

\begin{figure*}[t]
\centering

\begin{minipage}[t]{0.235\textwidth}
\centering
        \includegraphics[width=\linewidth]{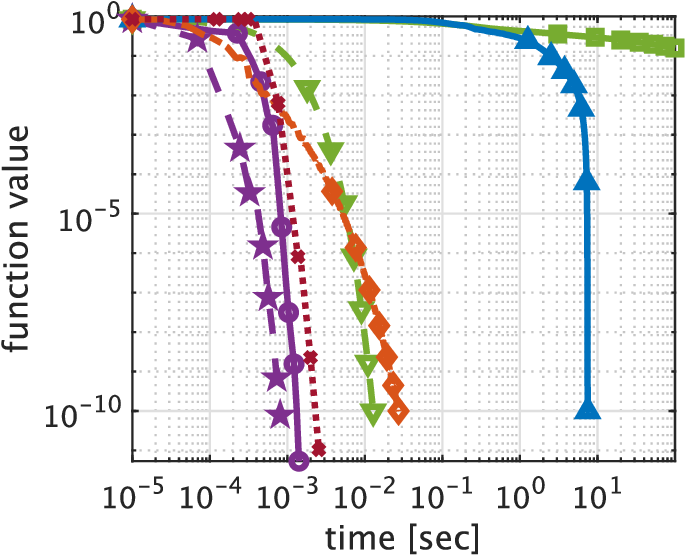}\\[-1mm]
{\footnotesize (a) \((d,m)=(100,10)\), time}
\end{minipage}
\hfill
\begin{minipage}[t]{0.235\textwidth}
\centering
        \includegraphics[width=\linewidth]{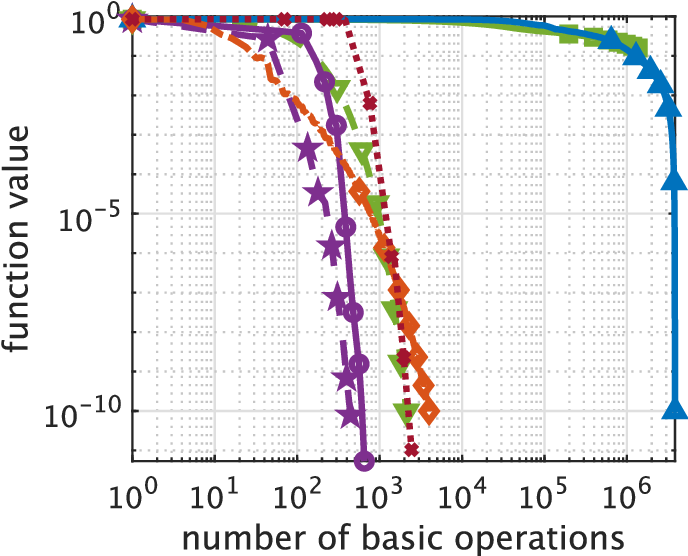}\\[-1mm]
{\footnotesize (b) \((d,m)=(100,10)\), ops.}
\end{minipage}
\hfill
\begin{minipage}[t]{0.235\textwidth}
\centering
        \includegraphics[width=\linewidth]{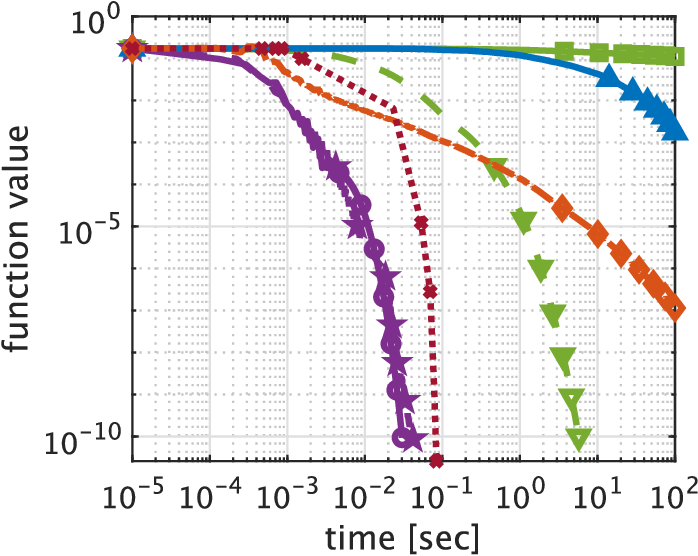}\\[-1mm]
{\footnotesize (c) \((d,m)=(100,50)\), time}
\end{minipage}
\hfill
\begin{minipage}[t]{0.235\textwidth}
\centering
        \includegraphics[width=\linewidth]{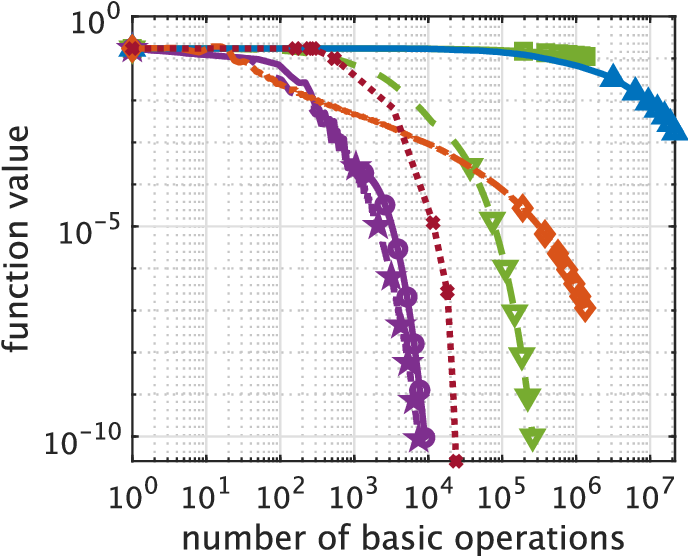}\\[-1mm]
{\footnotesize (d) \((d,m)=(100,50)\), ops.}
\end{minipage}

\vspace{1mm}

\begin{minipage}[t]{0.235\textwidth}
\centering
        \includegraphics[width=\linewidth]{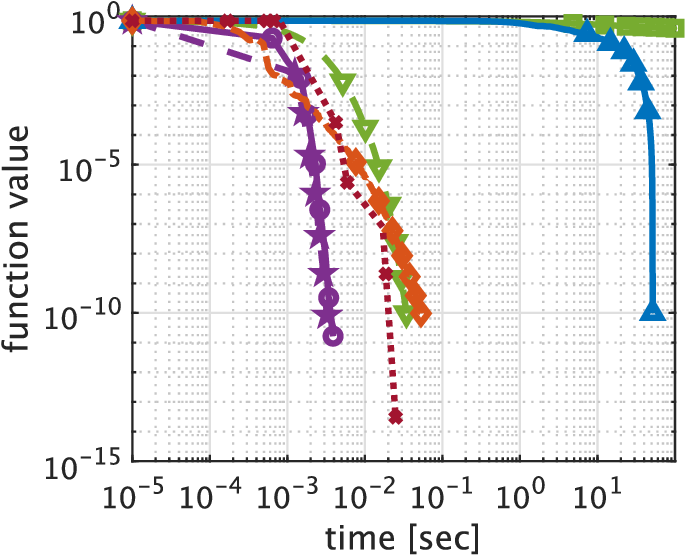}\\[-1mm]
{\footnotesize (e) \((d,m)=(1000,10)\), time}
\end{minipage}
\hfill
\begin{minipage}[t]{0.235\textwidth}
\centering
        \includegraphics[width=\linewidth]{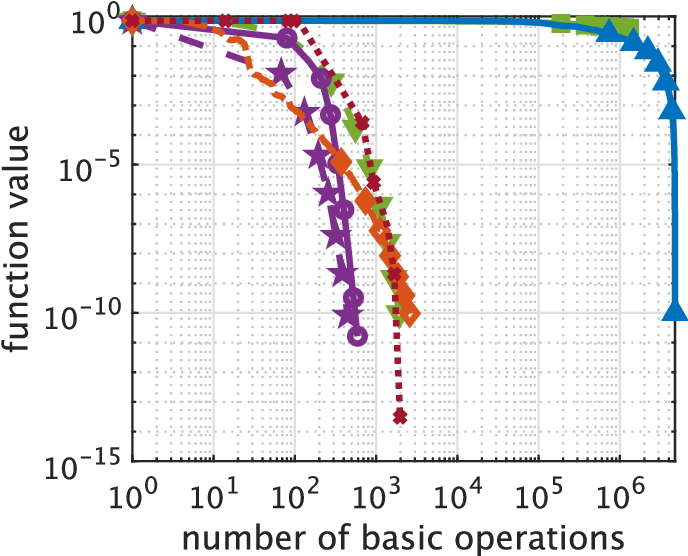}\\[-1mm]
{\footnotesize (f) \((d,m)=(1000,10)\), ops.}
\end{minipage}
\hfill
\begin{minipage}[t]{0.235\textwidth}
\centering
        \includegraphics[width=\linewidth]{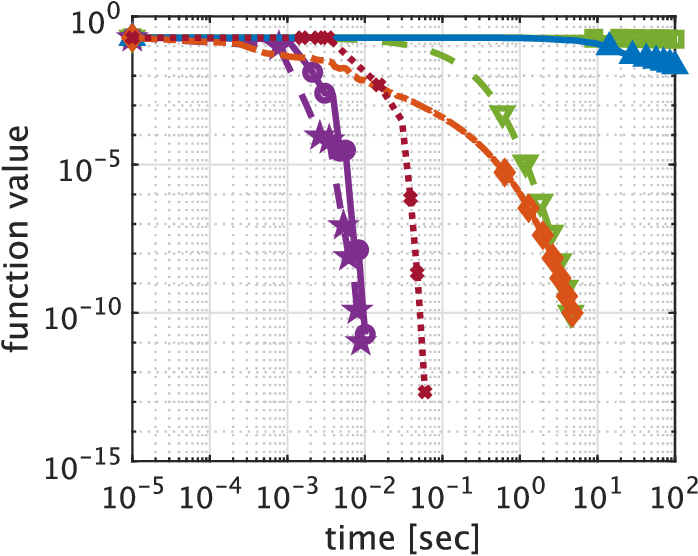}\\[-1mm]
{\footnotesize (g) \((d,m)=(1000,50)\), time}
\end{minipage}
\hfill
\begin{minipage}[t]{0.235\textwidth}
\centering
        \includegraphics[width=\linewidth]{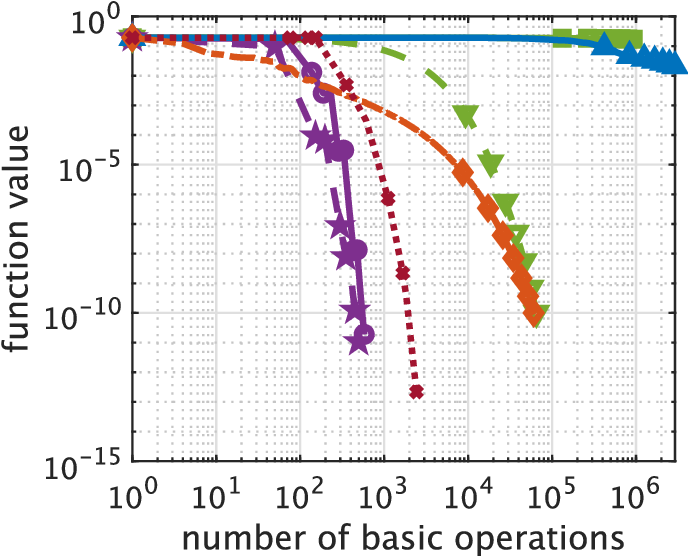}\\[-1mm]
{\footnotesize (h) \((d,m)=(1000,50)\), ops.}
\end{minipage}

\vspace{1mm}

\includegraphics[width=0.82\textwidth]{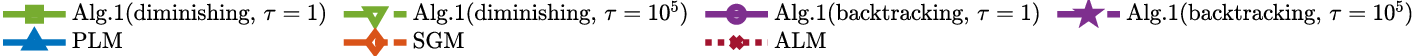}
\vspace{-1em}
\caption{
Comparison of optimization algorithms in a scenario of robust target localization for four problem sizes.
For each problem size, the left figure shows the value of the cost function versus CPU time,
and the right one shows the value of the cost function versus the cumulative number of
evaluations of basic operations (ops.).
For visualization of the results of PLM, we recorded the value of the cost function at every $100$ inner iteration.
}
\label{fig:maxmin-four-sizes}
\vspace{-1.5em}
\end{figure*}

Fig.~\ref{fig:maxmin-four-sizes} shows the values of the cost function against CPU time (sec) and the cumulative number of evaluations of basic operations (ops.; see the end of Section~\ref{sec:measure}) for one problem instance at each of the four problem sizes, where one evaluation of a subgradient of
$g$
was also counted as $1$ basic operation in SGM.
The same problem instance was used for all algorithms at each problem size.
In problem sizes,
$d$
is the dimension of the target
$\bm{x}^{\star}$,
and
$m$
is the number of component functions in the maximum in~\eqref{eq:max_problem}.
To compare single-loop algorithms (Alg.~\ref{alg:proposed} and SGM) with double-loop algorithms (PLM and ALM), the cumulative number of evaluations of basic operations is reported rather than the number of iterations.
Table~\ref{table:deteailed_average} summarizes the averaged value of the cost function, CPU time and the number of evaluations of basic operations, over
$100$
problem instances,
to reach stopping criterion.

\begin{table}[t] 
\centering
\caption{
  Averaged final value of cost function (cost), running CPU time (time) [sec], and number of basic operations (ops.).
}
\label{table:deteailed_average}
\label{tab:average_performance}
\setlength{\tabcolsep}{4pt}
\renewcommand{\arraystretch}{1.08}
\begin{tabular}{llccc}
\toprule
$(d,m)$ & Algorithm
& cost
& time
& ops. \\
\midrule

\multirow{7}{*}{$(100,10)$}
& Alg. 1-bac. $(\tau=1)$          &\underline{4.25E-11}&{\bf 1.16E-03}&\underline{6.28E+02}\\
& Alg. 1-bac. $(\tau=10^5)$       &4.48E-11&\underline{1.95E-03}& {\bf 5.53E+02}\\
& Alg. 1-dim. $(\tau=1)$     &8.55E-02&1.00E+02&1.35E+06\\
& Alg. 1-dim. $(\tau=10^5)$  &9.74E-11&1.54E-02&2.27E+03\\
& PLM                         &9.87E-11&4.65E+00&2.38E+06\\
& SGM                         &9.89E-11&3.22E-02&4.49E+03\\
& ALM                         &{ \bf 3.12E-11}&2.87E-03&2.77E+03\\
\midrule
\multirow{7}{*}{$(100,50)$}
& Alg. 1-bac. $(\tau=1)$          & \underline{8.59E-11}& {\bf 1.63E-02}& {\bf 4.56E+03}\\
& Alg. 1-bac. $(\tau=10^5)$       & 8.62E-11& \underline{2.09E-02}& \underline{5.34E+03}\\
& Alg. 1-dim. $(\tau=1)$     & 9.79E-02&1.00E+02&1.35E+06\\
& Alg. 1-dim. $(\tau=10^5)$  & 1.00E-10&3.57E+00&1.81E+05\\
& PLM                         & 1.87E-03&1.00E+02&2.11E+07\\
& SGM                         & 3.02E-08&9.37E+01&1.23E+06\\
& ALM                         & {\bf 1.38E-11}&6.93E-02&1.98E+04\\
\midrule
\multirow{7}{*}{$(1000,10)$}
& Alg. 1-bac. $(\tau=1)$          &\underline{2.96E-11}&\underline{3.21E-03}&\underline{5.84E+02}\\
& Alg. 1-bac. $(\tau=10^5)$       &3.18E-11& {\bf 2.67E-03}&{\bf 4.37E+02}\\
& Alg. 1-dim. $(\tau=1)$     &3.01E-01&1.00E+02&1.30E+06\\
& Alg. 1-dim. $(\tau=10^5)$  &9.67E-11&3.42E-02&1.73E+03\\
& PLM                         &9.74E-11&4.48E+01&4.30E+06\\
& SGM                         &9.88E-11&5.43E-02&2.23E+03\\
& ALM                         & {\bf 1.72E-11}&2.02E-02&1.44E+03\\
\midrule
\multirow{7}{*}{$(1000,50)$}
& Alg. 1-bac. $(\tau=1)$          &3.60E-11& \underline{1.13E-02}&\underline{6.41E+02}\\
& Alg. 1-bac. $(\tau=10^5)$       &\underline{3.35E-11}&{\bf 9.11E-03}&{\bf 5.01E+02}\\
& Alg. 1-dim. $(\tau=1)$     &1.57E-01&1.00E+02&9.17E+05\\
& Alg. 1-dim. $(\tau=10^5)$  &9.99E-11&4.74E+00&6.97E+04\\
& PLM                         &3.92E-02&1.00E+02&2.93E+06\\
& SGM                         &9.99E-11&5.55E+00&6.83E+04\\
& ALM                         & {\bf 1.15E-11}&5.53E-02&2.05E+03\\
\bottomrule
\end{tabular}
\end{table}

From Fig.~\ref{fig:maxmin-four-sizes} and Table~\ref{table:deteailed_average}, Alg.~\ref{alg:proposed} with backtracking-based stepsizes achieves the stopping criterion
$(\cost+\phi)(\bm{x}_{n}) < 10^{-10}$
within the shortest CPU time and the smallest cumulative number of evaluations of basic operations.
In contrast, Alg.~\ref{alg:proposed} with diminishing stepsizes is more sensitive to the choice
of
$\tau$.
Indeed, a convergence speed of Alg.~\ref{alg:proposed} with diminishing stepsizes with the choice
$\tau=1$
is significantly slower than that with the choice
$\tau=10^{5}$,
and
Alg.~\ref{alg:proposed} with diminishing stepsizes with
$\tau=1$
fails to reach the stopping criterion
$(\cost+\phi)(\bm{x}_{n}) < 10^{-10}$.
A possible reason for this sensitivity is that a small
$\tau$
leads to small stepsizes
$\gamma_{n}=2(1-c)(\varpi_{1} + \tau^{-1}\varpi_{2}n^{1/3})^{-1}$,
which causes slow convergence in practice.
This sensitivity is mitigated by backtracking-based stepsizes because Alg.~\ref{alg:backtracking} finds a stepsize according to the local geometry around the current estimate (see also Remark~\ref{remark:stepsize}).

PLM is substantially slower than Alg.~\ref{alg:proposed} in our experiments, and does not achieve the stopping criterion
$(\cost+\phi)(\bm{x}_{n}) < 10^{-10}$
within
$100$ seconds for the cases with the larger value
$m=50$.
A possible reason is that each outer iteration of PLM requires solving a smoothed proximal linear
subproblem by a fast gradient method~\cite[Alg. 5]{Drusvyatskiy-Paquette19} whose update uses
the gradient of the smoothed and linearized function
$\bm{x}\mapsto \moreau{g}{\nu}(\mathfrak{S}(\bm{x}_{n}) + \mathrm{D}\mathfrak{S}(\bm{x}_{n})[\bm{x}-\bm{x}_{n}])$,
which is
$\opnorm{\mathrm{D}\mathfrak{S}(\bm{x}_{n})}^{2}\nu^{-1}$-Lipschitz continuous,
at the $n$th outer iteration.
A smoothing parameter
$\nu\coloneqq \epsilon^{2}/(2L_{g}^{3}L_{\mathrm{D}\mathfrak{S}}) =\epsilon^{2}/(2L_{\mathrm{D}\mathfrak{S}})$~\cite[p. 534]{Drusvyatskiy-Paquette19}
is chosen to achieve an $\epsilon$-stationarity;
hence the Lipschitz constant of the gradient can be large for a small target tolerance
$\epsilon$
($\epsilon = 10^{-2}$
was used in our experiments).
Since the reciprocal of the Lipschitz constant is used as a stepsize in~\cite[Alg. 5]{Drusvyatskiy-Paquette19}, many inner iterations are required for PLM due to small stepsizes (see also~\cite[p. 535]{Drusvyatskiy-Paquette19}).
Although PLM is a theoretically natural baseline (see Table~\ref{table:summary}), its inner-loop cost can be significant in our numerical experiments.

ALM reaches the stopping criterion
$(\cost+\phi)(\bm{x}_{n})<10^{-10}$
in all tested instances, but it tends to require more evaluations of basic operations than Alg.~\ref{alg:proposed} with backtracking-based stepsizes.
This is because ALM solves subproblems inexactly by an inner iterative solver.
From a viewpoint of CPU time, in a small problem size case
$(d,m)=(100,10)$, 
ALM has comparable performance to Alg.~\ref{alg:proposed} with backtracking-based stepsizes although ALM evaluates approximately $4$ times as many basic operations as Alg.~\ref{alg:proposed} with backtracking-based stepsizes, as shown in Table~\ref{table:deteailed_average}.
This comparable CPU-time performance may be due to the use of different
programming languages:
ALM was implemented by Julia in~\cite{DeMarchi-Jia-Kanzow-Mehlitz23}, whereas the other methods were implemented by MATLAB.
We note that ALM achieves the smallest averaged final value of the cost function from Table~\ref{tab:average_performance} because an inner iterative solver for subproblems in ALM did not terminate after
$(\cost+\phi)(\bm{x}_{n,k})<10^{-10}$
was achieved with an estimate
$\bm{x}_{n,k}$
in the inner loop.

SGM decreases the objective value steadily, but its progress is relatively
slow in particular for the case 
$(d,m)=(100,50)$.
SGM typically requires more CPU time and more basic operations than Alg.~\ref{alg:proposed} with backtracking-based stepsizes.
In contrast,
except for the case
$(d,m)=(100,50)$,
convergence speeds of SGM are comparable to those of Alg.~\ref{alg:proposed} with diminishing stepsizes and
$\tau=10^{5}$,
where similar diminishing stepsizes were employed in SGM.

Overall, these results demonstrate that Alg.~\ref{alg:proposed} with the
backtracking-based stepsizes is effective compared with  PLM, ALM and SGM in this experiment.

\subsection{Application to MU-MIMO Signal Detection} \label{sec:numerical:MIMO}
\captionsetup{aboveskip=0em}
\begin{figure}[t]
\centering

\begin{minipage}[t]{0.235\textwidth}
\centering
    \includegraphics[clip, width=\columnwidth]{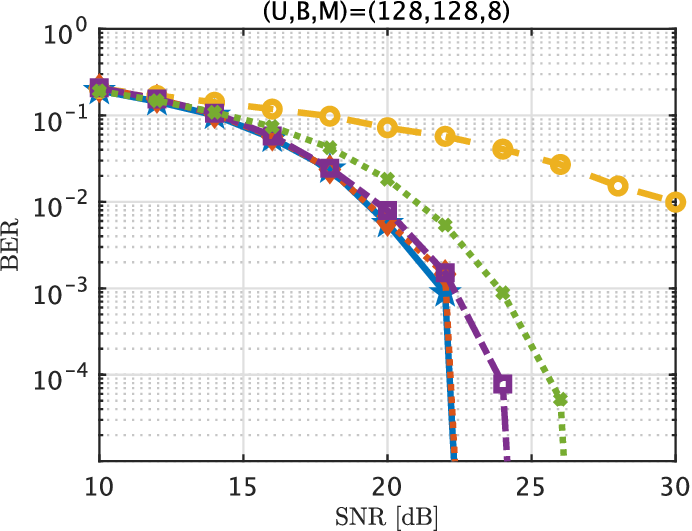}
    {\footnotesize (a) Case $B=U$}
\end{minipage}
\hfill
\begin{minipage}[t]{0.235\textwidth}
\centering
    \includegraphics[clip,width=\columnwidth]{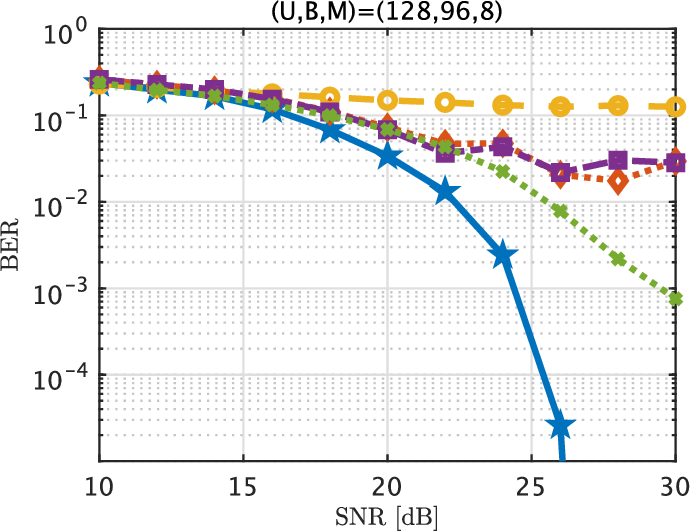}
    {\footnotesize (b) Case $B=\frac{3}{4}U$}
\end{minipage}

\vspace{1mm}

\includegraphics[width=0.4\textwidth]{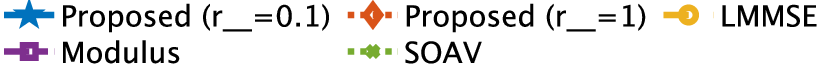}

\caption{
  Comparison (BER vs SNR) of models for MIMO detection~\eqref{eq:MIMO_model}.
}
\label{fig:BER}
\end{figure}

For the MU-MIMO signal detection problem~\eqref{eq:MIMO_model}, we evaluated the estimation accuracy of the proposed model~\eqref{eq:MIMO_proposed} and Algorithm~\ref{alg:proposed}.
In this experiment, according to the setting in~\cite{Chen19},
we randomly chose
(i)
$\mathsf{s}^{\star} \in \mathsf{D}$;
(ii)
$\mathsf{H} \coloneqq \sqrt{\mathsf{R}}\mathsf{G}$
with
$\mathsf{G} \in \mathbb{C}^{B\times U}$
whose entries were sampled from the complex Gaussian distribution
$\mathbb{C}\mathcal{N}(0,1/B)$,
and
$\mathsf{R} \in \mathbb{R}^{B\times B}$
given by
$[\mathsf{R}]_{j,k} = 0.5^{\abs{j-k}}$;
and
(iii)
$\mathsf{e} \in \mathbb{C}^{B}$
whose entries were sampled from
$\mathbb{C}\mathcal{N}(0,\sigma^{2})$,
where
$\sigma^{2}$
was chosen so that
$10\log_{10}\frac{1}{\sigma^{2}}\,\mathrm{(dB)}$
achieved a given signal-to-noise ratio (SNR).

We compared the estimation accuracy of
(i) the proposed model~\eqref{eq:MIMO_proposed} with
$\underline{r} \in \{0.1, 1\}$
solved by Algorithm~\ref{alg:proposed} (with 
$\tau =2^{-1}$
in
$(\mu_{n})_{n=1}^{\infty}\coloneqq (\tau n^{-1/3})_{n=1}^{\infty}$,
and
backtracking-based stepsizes)
with the three models (see Section~\ref{sec:existing_mimo}): (ii) LMMSE~\cite{Yang-Hanzo15}, (iii)
Modulus in~\eqref{eq:modulus} solved by a projected gradient method in~\cite{Chen19},
and (iv) SOAV~\cite{Hayakawa-Hayashi20} in~\eqref{eq:SOAV} solved by a primal-dual splitting algorithm in~\cite{Condat13},
where
$\lambda_{r}, \lambda_{\theta}$
in~\eqref{eq:MIMO_proposed} and~\eqref{eq:polar_regularizer},
and
$\lambda$
in~\eqref{eq:SOAV} were chosen from
$\{10^{k}\mid k=-6,-5,\ldots,0\}$.
Each algorithm except LMMSE was terminated when
$\stdnorm{\bm{x}_{n+1}-\bm{x}_{n}} < 10^{-5}$
was achieved with its generated sequence
$\bm{x}_{n}$,
or
running CPU time exceeded
$5$ seconds.
We evaluated their performance by the averaged bit error rate (BER) of the final estimate over
$100$
trials, where BER was computed by a MATLAB code with
'pskdemod' and 'biterr'.
The lower averaged BER indicates better performance.

Fig.~\ref{fig:BER} demonstrates the averaged BER of each algorithm versus SNR under respectively two settings
$B=U$
and
$B=\frac{3}{4}U$,
where
r\_\_
in Fig.~\ref{fig:BER}
stands for
$\underline{r}$
in the proposed model~\eqref{eq:MIMO_proposed},
and
the averaged BER
$0$
was replaced by the machine epsilon for visualization.
We note that the estimation task~\eqref{eq:MIMO_model} becomes challenging as the ratio
$B/U$
becomes small.
From Fig.~\ref{fig:BER}, we observe that (i) every model except LMMSE achieves better performance than LMMSE;
(ii) the proposed model~\eqref{eq:MIMO_proposed} with
$\underline{r} = 0.1$
outperforms the others for all SNRs and ratios
$B/U$;
(iii) for the
$B=\frac{3}{4}U$
case, the estimation performance of the proposed model with
$\underline{r} = 1$
is deteriorated compared with the proposed model with
$\underline{r} = 0.1$,
although
the proposed models with
$\underline{r} = 1$
and
$\underline{r} = 0.1$
have comparable performance for the
$B=U$
case.
Moreover, from Fig.~\ref{fig:BER} (b) in the challenging setting with
$B=\frac{3}{4}U$,
we see that the proposed model~\eqref{eq:MIMO_proposed} with
$\underline{r} = 0.1$
overwhelms the others.
As we expected, these results show a great potential of the proposed model~\eqref{eq:MIMO_proposed} for MU-MIMO signal detection.

\section{Conclusion}
We proposed a single-loop proximal variable smoothing algorithm for nonsmooth optimization with the cost function
$h+g\circ\mathfrak{S} + \phi$
in Problem~\ref{problem:origin}.
The proposed algorithm consists of (I) a gradient descent step with a smoothed surrogate function
$h+\moreau{g}{\mu}\circ\mathfrak{S}$
with the Moreau envelope
$\moreau{g}{\mu}$
of
$g$,
where
$\mu$
diminishes to zero;
(II) a proximal step with the proximity operator of
$\phi$.
For the proposed algorithm,
we presented a subsequential convergence guarantee and a convergence rate in terms of a stationary point in Theorem~\ref{theorem:convergence_extension} by exploiting properties of a gradient mapping-type stationarity measure (see Propositions~\ref{proposition:stationary_relation} and~\ref{proposition:measure}).
Numerical experiments demonstrate that (i)
for robust target localization,
the proposed algorithm achieves a faster numerical convergence speed compared to existing algorithm; and (ii) for MIMO signal detection, the proposed model achieves a better estimation accuracy than existing models.

\bibliographystyle{IEEEtran}
\bibliography{main}%

\begin{thebibliography}{10}
\providecommand{\url}[1]{#1}
\csname url@samestyle\endcsname
\providecommand{\newblock}{\relax}
\providecommand{\bibinfo}[2]{#2}
\providecommand{\BIBentrySTDinterwordspacing}{\spaceskip=0pt\relax}
\providecommand{\BIBentryALTinterwordstretchfactor}{4}
\providecommand{\BIBentryALTinterwordspacing}{\spaceskip=\fontdimen2\font plus
\BIBentryALTinterwordstretchfactor\fontdimen3\font minus \fontdimen4\font\relax}
\providecommand{\BIBforeignlanguage}[2]{{%
\expandafter\ifx\csname l@#1\endcsname\relax
\typeout{** WARNING: IEEEtran.bst: No hyphenation pattern has been}%
\typeout{** loaded for the language `#1'. Using the pattern for}%
\typeout{** the default language instead.}%
\else
\language=\csname l@#1\endcsname
\fi
#2}}
\providecommand{\BIBdecl}{\relax}
\BIBdecl

\bibitem{Kume-Yamada25A}
K.~Kume and I.~Yamada, ``A proximal variable smoothing for nonsmooth minimization involving weakly convex composite with {MIMO} application,'' in \emph{IEEE ICASSP}, 2025, pp. 1--5.

\bibitem{Tibshirani96}
R.~Tibshirani, ``Regression shrinkage and selection via the lasso,'' \emph{J. R. Stat. Soc. Ser. B Methodol.}, vol.~58, no.~1, pp. 267--288, 1996.

\bibitem{Chen-Gu14}
L.~Chen and Y.~Gu, ``The convergence guarantees of a non-convex approach for sparse recovery,'' \emph{IEEE Trans. Signal Process.}, vol.~62, no.~15, pp. 3754--3767, 2014.

\bibitem{Selesnick17}
I.~Selesnick, ``Sparse regularization via convex analysis,'' \emph{IEEE Trans. Signal Process.}, vol.~65, no.~17, pp. 4481--4494, 2017.

\bibitem{Abe-Yamagishi-Yamada20}
J.~Abe, M.~Yamagishi, and I.~Yamada, ``Linearly involved generalized {M}oreau enhanced models and their proximal splitting algorithm under overall convexity condition,'' \emph{Inverse Problems}, vol.~36, no.~3, 2020.

\bibitem{Yata-Yamagishi-Yamada22}
W.~Yata, M.~Yamagishi, and I.~Yamada, ``{A constrained LiGME model and its proximal splitting algorithm under overall convexity condition},'' \emph{J. Appl. Numer. Optim.}, vol.~4, no.~2, pp. 245--271, 2022.

\bibitem{Kuroda-Kitahara22}
H.~Kuroda and D.~Kitahara, ``Block-sparse recovery with optimal block partition,'' \emph{IEEE Trans. Signal Process.}, vol.~70, pp. 1506--1520, 2022.

\bibitem{Lu-Yan-Lin16}
C.~Lu, S.~Yan, and Z.~Lin, ``Convex sparse spectral clustering: single-view to multi-view,'' \emph{IEEE Trans. Image Process.}, vol.~25, no.~6, pp. 2833--2843, 2016.

\bibitem{Wang-Liu-Chen-Ma-Xue-Zhao22}
Z.~Wang, B.~Liu, S.~Chen, S.~Ma, L.~Xue, and H.~Zhao, ``A manifold proximal linear method for sparse spectral clustering with application to single-cell {RNA} sequencing data analysis,'' \emph{INFORMS J. on Optim.}, vol.~4, no.~2, pp. 200--214, 2022.

\bibitem{Kume-Yamada24A}
K.~Kume and I.~Yamada, ``A variable smoothing for weakly convex composite minimization with manifold constraint via parametrization,'' \emph{arXiv (2412.04225v4)}, 2024.

\bibitem{Duchi-Ruan18}
J.~C. Duchi and F.~Ruan, ``{Solving (most) of a set of quadratic equalities: composite optimization for robust phase retrieval},'' \emph{Information and Inference: A Journal of the IMA}, vol.~8, no.~3, pp. 471--529, 2018.

\bibitem{Zhen-Ma-Xue24}
Z.~Zheng, S.~Ma, and L.~Xue, ``A new inexact proximal linear algorithm with adaptive stopping criteria for robust phase retrieval,'' \emph{IEEE Trans. Signal Process.}, vol.~72, pp. 1081--1093, 2024.

\bibitem{Yazawa-Kume-Yamada26}
K.~Yazawa, K.~Kume, and I.~Yamada, ``A {DC} composite optimization via variable smoothing for robust phase retrieval with nonconvex loss functions,'' \emph{arXiv (2604.07686v3)}, 2026.

\bibitem{Charisopoulos-Chen-Davis-Diaz-Ding-Drusvyatskiy21}
V.~Charisopoulos, Y.~Chen, D.~Davis, M.~D{\'i}az, L.~Ding, and D.~Drusvyatskiy, ``Low-rank matrix recovery with composite optimization: Good conditioning and rapid convergence,'' \emph{Foundations of Computational Mathematics}, vol.~21, no.~6, pp. 1505--1593, 2021.

\bibitem{Wang-So-Zoubir23}
Z.-Y. Wang, H.~C. So, and A.~M. Zoubir, ``Robust low-rank matrix recovery via hybrid ordinary-{W}elsch function,'' \emph{IEEE Trans. Signal Process.}, vol.~71, pp. 2548--2563, 2023.

\bibitem{Yang-Chen-Ma-Chen-Gu-So19}
C.~Yang, X.~Shen, H.~Ma, B.~Chen, Y.~Gu, and H.~C. So, ``Weakly convex regularized robust sparse recovery methods with theoretical guarantees,'' \emph{IEEE Trans. Signal Process.}, vol.~67, no.~19, pp. 5046--5061, 2019.

\bibitem{Suzuki-Yukawa23}
K.~Suzuki and M.~Yukawa, ``Sparse stable outlier-robust signal recovery under {G}aussian noise,'' \emph{IEEE Trans. Signal Process.}, vol.~71, pp. 372--387, 2023.

\bibitem{Zhang10}
C.-H. Zhang, ``{Nearly unbiased variable selection under minimax concave penalty},'' \emph{Ann. Stat.}, vol.~38, no.~2, pp. 894 -- 942, 2010.

\bibitem{Fan-Li01}
J.~Fan and R.~Li, ``\BIBforeignlanguage{English (US)}{Variable selection via nonconcave penalized likelihood and its oracle properties},'' \emph{\BIBforeignlanguage{English (US)}{J. Am. Stat. Assoc.}}, vol.~96, no. 456, pp. 1348--1360, 2001.

\bibitem{Lewis-Wright16}
A.~S. Lewis and S.~J. Wright, ``A proximal method for composite minimization,'' \emph{Math. Program.}, vol. 158, no.~1, pp. 501--546, 2016.

\bibitem{Drusvyatskiy-Paquette19}
D.~Drusvyatskiy and C.~Paquette, ``{Efficiency of minimizing compositions of convex functions and smooth maps},'' \emph{Math. Program.}, vol. 178, no. 1-2, pp. 503--558, 2019.

\bibitem{DeMarchi-Jia-Kanzow-Mehlitz23}
A.~De~Marchi, X.~Jia, C.~Kanzow, and P.~Mehlitz, ``Constrained composite optimization and augmented {L}agrangian methods,'' \emph{Math. Program.}, vol. 201, no.~1, pp. 863--896, Sep 2023.

\bibitem{Bourkhissi-Necoara-Patrinos-TranDinh25}
L.~E. Bourkhissi, I.~Necoara, P.~Patrinos, and Q.~Tran-Dinh, ``Complexity of linearized perturbed augmented {L}agrangian methods for nonsmooth nonconvex optimization with nonlinear equality constraints,'' \emph{arXiv (2503.01056v1)}, 2025.

\bibitem{Davis-Drusvyatskiy-MacPhee18}
D.~Davis, D.~Drusvyatskiy, K.~J. MacPhee, and C.~Paquette, ``Subgradient methods for sharp weakly convex functions,'' \emph{J. Optim. Theory Appl.}, vol. 179, no.~3, pp. 962--982, 2018.

\bibitem{Li-Zhu-So-Man-Cho-Lee19}
X.~Li, Z.~Zhu, A.~M.-C. So, and J.~D. Lee, ``Incremental methods for weakly convex optimization,'' \emph{arXiv (1907.11687v2)}, pp. 1--25, 2019.

\bibitem{Zhu-Zhao-Zhang23}
D.-L. Zhu, L.~Zhao, and S.-Z. Zhang, ``A unified analysis on the subgradient upper bounds for the subgradient methods minimizing composite nonconvex, nonsmooth, and non-{L}ipschitz functions,'' \emph{J. Oper. Res. Soc. China.}, 2025.

\bibitem{Liu-Xia24}
Y.~Liu and F.~Xia, ``{Proximal variable smoothing method for three-composite nonconvex nonsmooth minimization with a linear operator},'' \emph{Numer. Algorithms}, vol.~96, no.~1, pp. 237--266, 2024.

\bibitem{Haines-Loeppky-Tseng-Wang13}
S.~Haines, J.~Loeppky, P.~Tseng, and X.~Wang, ``Convex relaxations of the weighted maxmin dispersion problem,'' \emph{SIAM J. Optim.}, vol.~23, no.~4, pp. 2264--2294, 2013.

\bibitem{Pan-Shen-Xu21}
W.~Pan, J.~Shen, and Z.~Xu, ``An efficient algorithm for nonconvex-linear minimax optimization problem and its application in solving weighted maximin dispersion problem,'' \emph{Comput. Optim. Appl.}, vol.~78, no.~1, pp. 287--306, 2021.

\bibitem{Lopes-Peres-Bilches25}
S.~L{\'o}pez-Rivera, P.~P{\'e}rez-Aros, and E.~Vilches, ``A projected variable smoothing for weakly convex optimization and supremum functions,'' \emph{J. Optim. Theory Appl.}, vol. 208, no.~3, p. 115, Feb 2026.

\bibitem{Domingos-Xavier24}
J.~Domingos and J.~Xavier, ``Robust target localization in {2D}: A value-at-risk approach,'' \emph{IEEE Trans. Signal Process.}, vol.~72, pp. 3028--3042, 2024.

\bibitem{Domingos-Xavier25}
------, ``Outlier-resilient model fitting via percentile losses: Methods for general and convex residuals,'' \emph{IEEE Signal Process. Lett.}, vol.~32, pp. 931--935, 2025.

\bibitem{Yang-Hanzo15}
S.~Yang and L.~Hanzo, ``Fifty years of {MIMO} detection: The road to large-scale {MIMO}s,'' \emph{IEEE Commun. Surveys Tuts.}, vol.~17, no.~4, pp. 1941--1988, 2015.

\bibitem{Chen19}
J.~C. Chen, ``Computationally efficient data detection algorithm for massive {MU-MIMO} systems using {PSK} modulations,'' \emph{IEEE Commun. Lett.}, vol.~23, no.~6, pp. 983--986, 2019.

\bibitem{Hayakawa-Hayashi20}
R.~Hayakawa and K.~Hayashi, ``Asymptotic performance of discrete-valued vector reconstruction via box-constrained optimization with sum of $\ell_{1}$ regularizers,'' \emph{IEEE Trans. Signal Process.}, vol.~68, pp. 4320--4335, 2020.

\bibitem{Shoji-Yata-Kume-Yamada25}
S.~Shoji, W.~Yata, K.~Kume, and I.~Yamada, ``An {LiGME} regularizer of designated isolated minimizers -- an application to discrete-valued signal estimation,'' \emph{IEICE Trans. Fundam.}, vol. E108.A, no.~12, pp. 1629--1642, 2025.

\bibitem{Bohm-Wright21}
A.~B{\"o}hm and S.~J. Wright, ``Variable smoothing for weakly convex composite functions,'' \emph{J. Optim. Theory Appl.}, vol. 188, no.~3, pp. 628--649, 2021.

\bibitem{Bot-Hendrich15}
R.~I. Bo{\c t} and C.~Hendrich, ``A variable smoothing algorithm for solving convex optimization problems,'' \emph{TOP}, vol.~23, no.~1, pp. 124--150, 2015.

\bibitem{Kume-Yamada24}
K.~Kume and I.~Yamada, ``A variable smoothing for nonconvexly constrained nonsmooth optimization with application to sparse spectral clustering,'' in \emph{IEEE ICASSP}, 2024, pp. 9296--9300.

\bibitem{Kume-Yamada21}
------, ``A global {C}ayley parametrization of {S}tiefel manifold for direct utilization of optimization mechanisms over vector spaces,'' in \emph{IEEE ICASSP}, 2021, pp. 5554--5558.

\bibitem{Kume-Yamada22}
------, ``Generalized left-localized {C}ayley parametrization for optimization with orthogonality constraints,'' \emph{Optimization}, vol.~73, no.~4, pp. 1113--1159, 2022.

\bibitem{Levin-Kileel-Boumal24}
E.~Levin, J.~Kileel, and N.~Boumal, ``The effect of smooth parametrizations on nonconvex optimization landscapes,'' \emph{Math. Program.}, vol. 209, no.~1, pp. 63--111, 2025.

\bibitem{Rockafellar-Wets98}
R.~Rockafellar and R.~J.-B. Wets, \emph{Variational Analysis}, 3rd~ed.\hskip 1em plus 0.5em minus 0.4em\relax Springer, 2010.

\bibitem{Li-So-Ma20}
J.~Li, A.~M.~C. So, and W.~K. Ma, ``Understanding notions of stationarity in nonsmooth optimization: A guided tour of various constructions of subdifferential for nonsmooth functions,'' \emph{IEEE Signal Process. Mag.}, vol.~37, no.~5, pp. 18--31, 2020.

\bibitem{Bauschke-Moursi-Wang21}
H.~H. Bauschke, W.~M. Moursi, and X.~Wang, ``Generalized monotone operators and their averaged resolvents,'' \emph{Math. Program.}, vol. 189, no.~1, pp. 55--74, 2021.

\bibitem{Yamada-Yukawa-Yamagishi11}
I.~Yamada, M.~Yukawa, and M.~Yamagishi, ``Minimizing the {M}oreau envelope of nonsmooth convex functions over the fixed point set of certain quasi-nonexpansive mappings,'' in \emph{Fixed-Point Algorithms for Inverse Problems in Science and Engineering}, H.~H. Bauschke, R.~S. Burachik, P.~L. Combettes, V.~Elser, D.~R. Luke, and H.~Wolkowicz, Eds.\hskip 1em plus 0.5em minus 0.4em\relax Springer New York, 2011, pp. 345--390.

\bibitem{Bauschke-Combettes17}
H.~H. Bauschke and P.~L. Combettes, \emph{Convex analysis and monotone operator theory in Hilbert spaces}, 2nd~ed.\hskip 1em plus 0.5em minus 0.4em\relax Springer, 2017.

\bibitem{Burke-Hoheisel17}
J.~V. Burke and T.~Hoheisel, ``Epi-convergence properties of smoothing by infimal convolution,'' \emph{Set-Valued Var. Anal.}, vol.~25, no.~1, pp. 1--23, 2017.

\bibitem{Hoheisel-Laborde-Oberman20}
T.~Hoheisel, L.~Maxime, and O.~Adam, ``A regularization interpretation of the proximal point method for weakly convex functions,'' \emph{J. Dyn. Games}, vol.~7, no.~1, pp. 79--96, 2020.

\bibitem{Beck17}
A.~Beck, \emph{First-Order Methods in Optimization}.\hskip 1em plus 0.5em minus 0.4em\relax SIAM, 2017.

\bibitem{Themelis-Stella-Patrinos18}
A.~Themelis, L.~Stella, and P.~Patrinos, ``Forward-backward envelope for the sum of two nonconvex functions: Further properties and nonmonotone linesearch algorithms,'' \emph{SIAM J. Optim.}, vol.~28, no.~3, pp. 2274--2303, 2018.

\bibitem{Pugh15}
C.~C. Pugh, \emph{Real mathematical analysis}, 2nd~ed.\hskip 1em plus 0.5em minus 0.4em\relax Springer, 2015.

\bibitem{Nagahara15}
M.~Nagahara, ``Discrete signal reconstruction by sum of absolute values,'' \emph{IEEE Signal Process. Lett.}, vol.~22, no.~10, pp. 1575--1579, 2015.

\bibitem{Ram^^c3^^adrez-Santamar^^c3^^ada06}
D.~Ram^^c3^^adrez and I.~Santamar^^c3^^ada, ``Regularised approach to detection of constant modulus signals in {MIMO} channels,'' \emph{Electron. Lett.}, vol.~42, pp. 184--186, 2006.

\bibitem{Condat13}
L.~Condat, ``A primal^^e2^^80^^93dual splitting method for convex optimization involving {L}ipschitzian, proximable and linear composite terms,'' \emph{J. Optim. Theory Appl.}, vol. 158, no.~2, pp. 460--479, 2013.

\end{thebibliography}


\begin{thebibliography}{1}
\providecommand{\url}[1]{#1}
\csname url@samestyle\endcsname
\providecommand{\newblock}{\relax}
\providecommand{\bibinfo}[2]{#2}
\providecommand{\BIBentrySTDinterwordspacing}{\spaceskip=0pt\relax}
\providecommand{\BIBentryALTinterwordstretchfactor}{4}
\providecommand{\BIBentryALTinterwordspacing}{\spaceskip=\fontdimen2\font plus
\BIBentryALTinterwordstretchfactor\fontdimen3\font minus \fontdimen4\font\relax}
\providecommand{\BIBforeignlanguage}[2]{{%
\expandafter\ifx\csname l@#1\endcsname\relax
\typeout{** WARNING: IEEEtran.bst: No hyphenation pattern has been}%
\typeout{** loaded for the language `#1'. Using the pattern for}%
\typeout{** the default language instead.}%
\else
\language=\csname l@#1\endcsname
\fi
#2}}
\providecommand{\BIBdecl}{\relax}
\BIBdecl

\bibitem{Rockafellar-Wets98}
R.~Rockafellar and R.~J.-B. Wets, \emph{Variational Analysis}, 3rd~ed.\hskip 1em plus 0.5em minus 0.4em\relax Springer, 2010.

\bibitem{Bauschke-Combettes17}
H.~H. Bauschke and P.~L. Combettes, \emph{Convex analysis and monotone operator theory in Hilbert spaces}, 2nd~ed.\hskip 1em plus 0.5em minus 0.4em\relax Springer, 2017.

\bibitem{Kume-Yamada24A}
K.~Kume and I.~Yamada, ``A variable smoothing for weakly convex composite minimization with manifold constraint via parametrization,'' \emph{arXiv (2412.04225v4)}, 2024.

\end{thebibliography}

\appendices
\def\thesection{Appendix \Alph{section}}

\newcounter{appnum}

\setcounter{theorem}{0}
\renewcommand{\thetheorem}{\Alph{appnum}.\arabic{theorem}}

\setcounter{equation}{0}
\renewcommand{\theequation}{\Alph{appnum}.\arabic{equation}}

\stepcounter{appnum}
\setcounter{equation}{0}
\section{Proof of Lemma~\ref{lemma:optimality}} \label{appendix:lemma:optimality}
\ref{enum:lemma:optimality:existence}
Since
$g$
is Lipschitz continuous, and
$h$
and
$\mathfrak{S}$
are locally Lipschitz continuous (or equivalently, strictly continuous) from continuous differentiabilities~\cite[Thm. 9.7]{Rockafellar-Wets98},
$\cost=h+g\circ\mathfrak{S}$
is also locally Lipschitz continuous over
$\mathcal{X}$.
Hence,
by~\cite[Thm. 9.13]{Rockafellar-Wets98},
$\Lsubdiff \cost(\widebar{\bm{x}}) \subset \mathcal{X}$
is nonempty and compact.
Then, there exists
a minimizer
$\widebar{\bm{v}} \in \Lsubdiff \cost(\widebar{\bm{x}})$
of the continuous function
$\mathcal{X}\ni\bm{v}\mapsto\stdnorm{\gamma^{-1}(\widebar{\bm{x}}-\prox{\gamma\phi}(\widebar{\bm{x}}-\gamma\bm{v}))}$
over the compact set
$\Lsubdiff \cost(\widebar{\bm{x}})$,
from which
$\mathcal{M}_{\gamma}^{\cost,\phi}(\widebar{\bm{x}}) = \stdnorm{\gamma^{-1}(\widebar{\bm{x}}-\prox{\gamma\phi}(\widebar{\bm{x}}-\gamma\widebar{\bm{v}}))}$
holds by~\eqref{eq:measure}.

\ref{enum:lemma:optimality:stationarity}
By Lemma~\ref{lemma:optimality}~\ref{enum:lemma:optimality:existence}, we have the relation:
$\mathcal{M}_{\gamma}^{\cost,\phi}(\widebar{\bm{x}}) = 0$
if and only if there exists
$\widebar{\bm{v}} \in \Lsubdiff \cost(\widebar{\bm{x}})$
such that
$\widebar{\bm{x}} = \prox{\gamma\phi}(\widebar{\bm{x}} - \gamma \widebar{\bm{v}})$.
Together with
$\prox{\gamma\phi} = (\Id+ \gamma \subdiff\phi )^{-1}$~\cite[Prop. 16.44]{Bauschke-Combettes17},
we get
$\mathcal{M}_{\gamma}^{\cost,\phi}(\widebar{\bm{x}}) = 0 \overset{\hphantom{(\clubsuit)}}{\Leftrightarrow} \exists \widebar{\bm{v}} \in \Lsubdiff  \cost(\widebar{\bm{x}})\ \mathrm{s.t.}\ \widebar{\bm{x}} = \prox{\gamma\phi}(\widebar{\bm{x}} - \gamma \widebar{\bm{v}})    \Leftrightarrow \exists \widebar{\bm{v}} \in \Lsubdiff  \cost(\widebar{\bm{x}})\ \mathrm{s.t.}\  (\Id + \gamma \subdiff\phi )(\widebar{\bm{x}}) \ni \widebar{\bm{x}} - \gamma \widebar{\bm{v}} \Leftrightarrow \Lsubdiff  \cost(\widebar{\bm{x}}) + \subdiff\phi (\widebar{\bm{x}}) \ni \bm{0} \overset{\hphantom{\;}\eqref{eq:sum_rule}\hphantom{\;}}{\Leftrightarrow} \Lsubdiff  (\cost+\phi)(\widebar{\bm{x}}) \ni \bm{0}$.

\stepcounter{appnum}
\setcounter{equation}{0}
\section{Proof of Proposition~\ref{proposition:stationary_relation}} \label{appendix:stationary_relation}
By~\cite[Prop. 16.44]{Bauschke-Combettes17}, we have
$\widebar{\bm{p}} = \prox{\gamma\phi}(\widebar{\bm{x}}-\gamma\nabla\cost^{\langle \mu \rangle}(\widebar{\bm{x}})) = (\Id + \gamma\subdiff\phi)^{-1}(\widebar{\bm{x}}-\gamma\nabla\cost^{\langle \mu \rangle}(\widebar{\bm{x}}))$,
from which
$\gamma^{-1}(\widebar{\bm{x}}-\widebar{\bm{p}}) - \nabla \cost^{\langle\mu\rangle}(\widebar{\bm{x}}) = \gamma^{-1}( \widebar{\bm{x}} - \gamma \nabla\cost^{\langle \mu \rangle} - \widebar{\bm{p}} ) \in \subdiff \phi(\widebar{\bm{p}})$
holds.
Then, we have
  {
    \thickmuskip=0.1\thickmuskip
    \medmuskip=0.1\medmuskip
    \thinmuskip=0.1\thinmuskip
    \arraycolsep=0.1\arraycolsep
    \begin{align}
       &
      \hspace{-1em}
       \dist(\bm{0},\nabla\cost^{\langle \mu \rangle}(\widebar{\bm{p}}) + \subdiff \phi (\widebar{\bm{p}})) \label{eq:surrogate_subdiff}                                                \\
       & \hspace{-1em}
       \leq \stdnorm{\nabla\cost^{\langle \mu \rangle}(\widebar{\bm{p}})+ \left(\gamma^{-1}(\widebar{\bm{x}}-\widebar{\bm{p}})- \nabla \cost^{\langle\mu\rangle}(\widebar{\bm{x}})\right)} \\
       & \hspace{-1em}\leq \stdnorm{\gamma^{-1}(\widebar{\bm{x}}-\widebar{\bm{p}})} + \stdnorm{\nabla\cost^{\langle \mu \rangle}(\widebar{\bm{p}}) - \nabla \cost^{\langle\mu\rangle}(\widebar{\bm{x}})}     \\
       & \hspace{-1em}\leq \stdnorm{\gamma^{-1}( \widebar{\bm{x}}-\widebar{\bm{p}} )} + L_{\nabla \cost^{\langle \mu\rangle}}\stdnorm{\widebar{\bm{p}} - \widebar{\bm{x}}}                                   \\
       & \hspace{-1em}= (1+\gamma L_{\nabla \cost^{\langle \mu \rangle}})\stdnorm{\gamma^{-1}(\widebar{\bm{x}} - \widebar{\bm{p}})}
       \overset{\eqref{eq:measure_smooth}}{=} (1+\gamma L_{\nabla \cost^{\langle \mu \rangle}})\mathcal{M}_{\gamma}^{\cost^{\langle \mu \rangle},\phi}(\widebar{\bm{x}}). \hspace{-3em}
    \end{align}}%
    From~\eqref{eq:chain_rule} and
$\nabla \moreau{g}{\mu}(\mathfrak{S}(\widebar{\bm{p}})) \in \Lsubdiff g(\prox{\mu g}(\mathfrak{S}(\widebar{\bm{p}})))$
by~\cite[Lemma 3.2]{Bohm-Wright21},
we get
$\nabla \cost^{\langle \mu \rangle}(\widebar{\bm{p}}) =
  \nabla h(\widebar{\bm{p}}) + (\mathrm{D}\mathfrak{S}(\widebar{\bm{p}}))^{*}[\nabla \moreau{g}{\mu}(\mathfrak{S}(\widebar{\bm{p}}))]
  \in \nabla h(\widebar{\bm{p}}) + (\mathrm{D}\mathfrak{S}(\widebar{\bm{p}}))^{*}[\Lsubdiff g(\prox{\mu g}(\mathfrak{S}(\widebar{\bm{p}})))]$.
  Hence,~\eqref{eq:surrogate_subdiff} is bounded below by $\dist(\bm{0}, \nabla h(\widebar{\bm{p}}) + (\mathrm{D}\mathfrak{S}(\widebar{\bm{p}}))^{*}[\Lsubdiff g(\prox{\mu g}(\mathfrak{S}(\widebar{\bm{p}})))] + \subdiff \phi(\widebar{\bm{p}}))$.
Together with~\eqref{eq:surrogate_subdiff}, we get the inequality~\eqref{eq:approximate_stationarity_inequality}.
The inequality~\eqref{eq:bound_prox} follows from
$\nabla \moreau{g}{\mu}(\bm{z}) = \mu^{-1}(\bm{z}-\prox{\mu g}(\bm{z}))$
in~\eqref{eq:Moreau_grad} and
$\sup_{\bm{z}\in\mathcal{Z}} \stdnorm{\nabla \moreau{g}{\mu}(\bm{z})} \leq L_{g}$~\cite[Lemma 3.3]{Bohm-Wright21}

\stepcounter{appnum}
\setcounter{equation}{0}
\section{Proof of Proposition~\ref{proposition:measure}} \label{appendix:lemma:set_convergence}
Our proof of Proposition~\ref{proposition:measure} exploits the {\em outer limit} defined originally for a set sequence.
Since we only use the outer limit for a vector sequence
$(\bm{v}_{n})_{n=1}^{\infty}\subset \mathcal{X}$,
we present its simple definition (see~\cite[Def. 4.1]{Rockafellar-Wets98} for a general definition):
\begin{equation}
  \Limsup_{n\to\infty}\bm{v}_{n}\coloneqq \{\bm{v} \in \mathcal{X} \mid \bm{v}\ \mathrm{is\ a\ cluster\ point\ of}\ (\bm{v}_{n})_{n=1}^{\infty}\},
\end{equation}
where
$\bm{v} \in \mathcal{X}$
is said to be a {\em cluster point} of
$(\bm{v}_{n})_{n=1}^{\infty} \subset \mathcal{X}$
if there exists a subsequence of
$(\bm{v}_{n})_{n=1}^{\infty}$
converging to
$\bm{v}$.
Here,
``$\Limsup$''
is used for the outer limit in order to avoid confusion with the limit superior
``$\limsup$''
for a real sequence.
We also exploit the following useful properties of outer limit.

\begin{fact}
  \label{fact:outer}
  For a sequence
  $(\bm{v}_{n})_{n=1}^{\infty}\subset \mathcal{X}$,
  the following hold:
  \begin{enumerate}[label=(\alph*)]
    \item \label{enum:outer_distance}
      {
      \thickmuskip=0.2\thickmuskip
      \medmuskip=0.2\medmuskip
      \thinmuskip=0.2\thinmuskip
      \arraycolsep=0.2\arraycolsep
    $\liminf\limits_{n\to\infty} \stdnorm{\bm{v}_{n}} = \dist(\bm{0},\Limsup_{n\to\infty}\bm{v}_{n})$}
          holds~\cite[Exe. 4.8]{Rockafellar-Wets98}.
    \item \label{enum:composite_outer_convergence}
          If
          $\mathcal{F}:\mathcal{X} \to \mathcal{X}$
          is continuous and
          $(\bm{v}_{n})_{n=1}^{\infty}$
          is bounded,
          then
          $\Limsup_{n\to\infty} \mathcal{F}(\bm{v}_{n}) = \mathcal{F}(\Limsup_{n\to \infty}\bm{v}_{n})$~\cite[Thm. 4.26]{Rockafellar-Wets98}.
  \end{enumerate}
\end{fact}

\begin{lemma} \label{lemma:smoothed_gradient}
  Consider Problem~\ref{problem:origin}.
  Let
  $\widebar{\bm{x}} \in\mathcal{X}$
  and
  $(\bm{x}_{n})_{n=1}^{\infty} \subset \mathcal{X}$
  satisfy
  $\lim_{n\to\infty}\bm{x}_{n}=\widebar{\bm{x}}$,
  and
  $(\mu_{n})_{n=1}^{\infty} (\subset (0,\eta^{-1}))\searrow 0$,
  Then, for
  $\cost = h +g\circ\mathfrak{S}$
  and
  $\cost^{\langle \mu \rangle}$
  in~\eqref{eq:surrogate_function} with the expression of
  $\nabla \cost^{\langle \mu \rangle}$
  in~\eqref{eq:gradient_expression},
  the following hold:
  \begin{enumerate}[label=(\alph*)]
    \item \label{enum:gradient_bounded}
          $(\nabla \cost^{\langle \mu_{n}\rangle}(\bm{x}_{n}))_{n=1}^{\infty}$
          is a bounded sequence.
    \item \label{enum:gradient_consistency}
          $\Lsubdiff \cost(\widebar{\bm{x}}) \supset \Limsup_{n\to\infty} \nabla \cost^{\langle \mu_{n}\rangle}(\bm{x}_{n})$
          holds.
  \end{enumerate}
\end{lemma}
\begin{proof}
  (a)
  For
  $n\in\mathbb{N}$,
  by the triangle inequality, we have
    { %
      \thickmuskip=0.2\thickmuskip
      \medmuskip=0.2\medmuskip
      \thinmuskip=0.2\thinmuskip
      \arraycolsep=0.2\arraycolsep
      \begin{equation}
        \begin{array}{l}
          \stdnorm{\nabla \cost^{\langle \mu_{n}\rangle}(\bm{x}_{n})}
          \leq \stdnorm{\nabla h(\bm{x}_{n})} + \stdnorm{(\mathrm{D}\mathfrak{S}(\bm{x}_{n}))^{*}[\nabla \moreau{g}{\mu_{n}}(\mathfrak{S}(\bm{x}_{n}))]}                \\
          \leq \stdnorm{\nabla h(\bm{x}_{n})} + \opnorm{(\mathrm{D}\mathfrak{S}(\bm{x}_{n}))^{*}}\stdnorm{\nabla \moreau{g}{\mu_{n}}(\mathfrak{S}(\bm{x}_{n}))]} \\
          \leq \stdnorm{\nabla h(\bm{x}_{n})} + \opnorm{(\mathrm{D}\mathfrak{S}(\bm{x}_{n}))^{*}}L_{g},
        \end{array}
      \end{equation}}%
  where the last inequality holds by $\sup_{\bm{z}\in\mathcal{Z}} \stdnorm{\nabla \moreau{g}{\mu}(\bm{z})} \leq L_{g}$~\cite[Lemma 3.3]{Bohm-Wright21}
  with the Lipschitz constant
  $L_{g}$
  of
  $g$.
  Since
  $\nabla h$
  and
  $\mathrm{D}\mathfrak{S}$
  are continuous and
  $\lim_{n\to\infty}\bm{x}_{n} =\widebar{\bm{x}}$,
  we get
  $\sup_{n\in\mathbb{N}}\stdnorm{\nabla h(\bm{x}_{n})} < +\infty$
  and
  $\sup_{n\in \mathbb{N}}\opnorm{(\mathrm{D}\mathfrak{S}(\bm{x}_{n}))^{*}} < +\infty$,
  implying thus
  $(\nabla \cost^{\langle \mu_{n}\rangle}(\bm{x}_{n}))_{n=1}^{\infty}$
  is bounded.

  (b)
  See~\cite[Thm. 4.4 (a)]{Kume-Yamada24A};
  a similar result under the convexity of
  $g$
  is found in~\cite[Thm. 5.3 d) with
    $\omega\coloneqq \frac{\stdnorm{\cdot}^{2}}{2}$]{Burke-Hoheisel17}.
\end{proof}

\begin{proof}[Proof of Proposition~\ref{proposition:measure}]
  Let
  $E\coloneqq \widebar{\bm{x}} - \gamma \Lsubdiff\cost(\widebar{\bm{x}}) \subset \mathcal{X}$
  and
  $\bm{v}_{n}\coloneqq \bm{x}_{n} - \gamma \nabla \cost_{n}(\bm{x}_{n})\in \mathcal{X}\ (n\in \mathbb{N})$.
  By
  $\lim_{n\to\infty}\bm{x}_{n}=\widebar{\bm{x}}$
  and
  Lemma~\ref{lemma:smoothed_gradient}~\ref{enum:gradient_consistency},
  we have
  \begin{equation}
    \Limsup_{n\to\infty} \bm{v}_{n}
    = \widebar{\bm{x}} - \gamma \Limsup_{n\to\infty}\nabla \cost_{n}(\bm{x}_{n})
    \subset \widebar{\bm{x}} - \gamma \Lsubdiff\cost(\widebar{\bm{x}}) = E. \label{eq:E_n_E}
  \end{equation}
  By Lemma~\ref{lemma:smoothed_gradient}~\ref{enum:gradient_bounded} with
  $\lim_{n\to\infty}\bm{x}_{n}=\widebar{\bm{x}}$,
  $(\nabla \cost_{n}(\bm{x}_{n}))_{n=1}^{\infty}$
  is a bounded sequence, and thus
  $(\bm{v}_{n})_{n=1}^{\infty}=(\bm{x}_{n} - \gamma \nabla \cost_{n}(\bm{x}_{n}))_{n=1}^{\infty}$
  is also a bounded sequence.
  Then, Fact~\ref{fact:outer}~\ref{enum:composite_outer_convergence} with the continuity of
  $\prox{\gamma\phi}$
  (see, e.g.,~\cite[Prop. 12.28]{Bauschke-Combettes17})
  yields
  \begin{equation}
    \begin{array}{ll}
      \Limsup_{n\to \infty} \left(\prox{\gamma\phi}(\bm{v}_{n})\right)
       & \overset{\rm Fact~\ref{fact:outer}~\ref{enum:composite_outer_convergence}}{=}\prox{\gamma\phi}(\Limsup_{n\to \infty} \bm{v}_{n}) \\
       & \overset{\hphantom{aaaa}\eqref{eq:E_n_E}\hphantom{aaa}}{\subset} \prox{\gamma\phi}(E).\vspace{-1em}
      \label{eq:limsup_prox}
    \end{array}
    \vspace{1em}
  \end{equation}
  Hence, the desired inequality is obtained by
  \begin{align}
     &
    \liminf_{n\to\infty}\mathcal{M}_{\gamma}^{\cost_{n},\phi}(\bm{x}_{n})
    \overset{\eqref{eq:measure_smooth}}{=} \gamma^{-1}\liminf_{n\to\infty}\stdnorm{\bm{x}_{n}-\prox{\gamma\phi}(\bm{v}_{n})}                                                                                               \\
     & \overset{\rm Fact~\ref{fact:outer}~\ref{enum:outer_distance}}{=} \gamma^{-1} \dist(\bm{0}, \Limsup_{n\to\infty}(\bm{x}_{n}- \prox{\gamma\phi}(\bm{v}_{n})))                 \\
     & \overset{\hphantom{\rm Fact~\ref{fact:outer}~\ref{enum:outer_distance}}}{=} \gamma^{-1} \dist(\bm{0}, \widebar{\bm{x}}-\Limsup_{n\to\infty}(\prox{\gamma\phi}(\bm{v}_{n}))) \\
     & \overset{\hphantom{aaaa}\rm\eqref{eq:limsup_prox}\hphantom{aaa}}{\geq} \gamma^{-1} \dist(\bm{0}, \widebar{\bm{x}}-\prox{\gamma\phi}(E))
    \overset{\eqref{eq:measure}}{=} \mathcal{M}_{\gamma}^{\cost,\phi}(\widebar{\bm{x}}),
  \end{align}
  where we used
  $E=\widebar{\bm{x}} - \gamma \Lsubdiff\cost(\widebar{\bm{x}})$
  in the last equality.
\end{proof}

\stepcounter{appnum}
\setcounter{equation}{0}
\section{Proof of Theorem~\ref{theorem:convergence_extension}}\label{appendix:convergence}
The proof is partially inspired by that for~\cite[Thm. 4.8]{Kume-Yamada24A}.

(a)
For
$n\in \mathbb{N}$,
recall
$\bm{x}_{n+1} = \prox{\gamma_{n}\phi}(\bm{x}_{n}-\gamma_{n}\nabla \cost_{n}(\bm{x}_{n}))$.
Since the Armijo condition~\eqref{eq:Armijo} is satisfied with a stepsize
$\gamma\coloneqq\gamma_{n}$
in both~\ref{enum:case:diminishing} and~\ref{enum:case:backtracking} (see Examples~\ref{example:stepsize:diminishing} and~\ref{example:stepsize:backtracking}),
we have the following inequality (see also~\eqref{eq:residual}):
\begin{equation}
  (\cost_{n}+\phi)(\bm{x}_{n+1})
  \leq (\cost_{n}+\phi)(\bm{x}_{n}) - c\gamma_{n} \left(\mathcal{M}_{\gamma_{n}}^{\cost_{n},\phi}(\bm{x}_{n})\right)^{2}.\hspace{-0.5em}
  \label{eq:extension:Armijo}
\end{equation}
The inequality~\eqref{eq:smoothed_uniformly_different_parameter} with
$\mu_{n+1}\leq\mu_{n}\ (n\in \mathbb{N})$
from~\eqref{eq:nonsummable} yields the following inequality, for every
$\bm{x} \in \mathcal{X}$,
  {
    \thickmuskip=0.3\thickmuskip
    \medmuskip=0.3\medmuskip
    \thinmuskip=0.3\thinmuskip
    \arraycolsep=0.3\arraycolsep
    \begin{equation}
      (\cost_{n+1}+\phi)(\bm{x})
      \leq  (\cost_{n} + \phi)(\bm{x}) + L_{g}^{2}(\mu_{n} - \mu_{n+1}). \label{eq:extension:f_k_moreau}
    \end{equation}}%
By combining~\eqref{eq:extension:Armijo} with~\eqref{eq:extension:f_k_moreau}, for
$n\in \mathbb{N}$,
we obtain
  {
    \thickmuskip=0.3\thickmuskip
    \medmuskip=0.3\medmuskip
    \thinmuskip=0.3\thinmuskip
    \arraycolsep=0.3\arraycolsep
    \begin{align}
       & \hspace{-2em}(\cost_{n+1}+\phi)(\bm{x}_{n+1})
      \overset{\eqref{eq:extension:f_k_moreau}}{\leq} (\cost_{n}+\phi)(\bm{x}_{n+1}) + L_{g}^{2}(\mu_{n} - \mu_{n+1})                                                                                                                            \\
       & \hspace{-2em}\overset{\eqref{eq:extension:Armijo}}{\leq} (\cost_{n}+\phi)(\bm{x}_{n}) - c\gamma_{n}\left(\mathcal{M}_{\gamma_{n}}^{\cost_{n},\phi}(\bm{x}_{n})\right)^{2} + L_{g}^{2}(\mu_{n}-\mu_{n+1}), \label{eq:extension:gradient_inequality_1}
    \end{align}}%
and thus
\begin{equation}
  (\cost_{n+1}+\phi)(\bm{x}_{n+1})
  \leq (\cost_{n}+\phi)(\bm{x}_{n})  + L_{g}^{2}(\mu_{n}-\mu_{n+1}). \label{eq:extension:gradient_inequality}
\end{equation}
By summing~\eqref{eq:extension:gradient_inequality_1} up,
we have
  {
    \thickmuskip=0.3\thickmuskip
    \medmuskip=0.3\medmuskip
    \thinmuskip=0.3\thinmuskip
    \arraycolsep=0.3\arraycolsep
    \begin{align}
       & c\sum_{n=\underline{k}}^{\overline{k}}  \gamma_{n}\left(\mathcal{M}_{\gamma_{n}}^{\cost_{n},\phi}(\bm{x}_{n})\right)^{2}                                                 \\
       & \leq (\cost_{\underline{k}}+\phi)(\bm{x}_{\underline{k}}) - (\cost_{\overline{k}+1}+\phi)(\bm{x}_{\overline{k}+1}) + L_{g}^{2}(\mu_{\underline{k}}-\mu_{\overline{k}+1}) \\
       & \leq (\cost_{\underline{k}}+\phi)(\bm{x}_{\underline{k}}) - \inf_{\bm{x}\in\mathcal{X}} (\cost_{1}+\phi)(\bm{x}) + L_{g}^{2}(\mu_{\underline{k}}-\mu_{\overline{k}+1})   \\
       & \leq (\cost_{1}+\phi)(\bm{x}_{1}) - \inf_{\bm{x}\in\mathcal{X}} (\cost_{1}+\phi)(\bm{x}) + L_{g}^{2}(\mu_{1}-\mu_{\overline{k}+1})                                       \\
       & \leq (\cost_{1}+\phi)(\bm{x}_{1}) - \inf_{\bm{x}\in\mathcal{X}} (\cost_{1}+\phi)(\bm{x}) + L_{g}^{2}\mu_{1} < +\infty, \label{eq:extension:sum_gradient_bound}
    \end{align}}%
where the second inequality follows by
$\inf_{\bm{x}\in \mathcal{X}}(\cost_{1}+\phi)(\bm{x}) \overset{\eqref{eq:smoothed_uniformly_different_parameter}}{\leq} \inf_{\bm{x}\in \mathcal{X}}(\cost_{\overline{k}+1}+\phi)(\bm{x}) \leq (\cost_{\overline{k}+1}+\phi)(\bm{x}_{\overline{k}+1})$,
and the third inequality follows by using~\eqref{eq:extension:gradient_inequality} recursively.
The last inequality in~\eqref{eq:extension:sum_gradient_bound} follows from
$\bm{x}_{1} \in \dom{\phi}$
and
$-\infty < \inf_{\bm{x}\in \mathcal{X}}(\cost_{1}+\phi)(\bm{x})$,
where the latter can be verified by
(i)
$(\cost+\phi)(\widebar{\bm{x}}) - \mu_{1}L_{g} \overset{\eqref{eq:smoothed_uniformly}}{\leq} (\cost_{1} + \phi)(\widebar{\bm{x}}) \ (\widebar{\bm{x}} \in \mathcal{X})$
holds with
$\cost=h+g\circ\mathfrak{S}$;
and
(ii)
$-\infty < \inf_{\bm{x}\in\mathcal{X}}(\cost+\phi)(\bm{x})$
from the condition~\ref{enum:problem:origin:minimizer} in Problem~\ref{problem:origin}.

Recall
$L_{\nabla \cost_{n}} = \varpi_{1} + \varpi_{2}\mu_{n}^{-1}$
is the Lipschitz constant of
$\nabla \cost_{n}$
in Assumption~\ref{assumption:Lipschitz_smoothed}.
By letting
$\beta \coloneqq\min\left\{\gamma_{\rm init}L_{\nabla \cost_{1}}, 2\rho(1-c)\right\} \in \mathbb{R}_{++}$,
we get
$\gamma_{n} \geq \beta L_{\nabla \cost_{n}}^{-1}\ (n\in\mathbb{N})$
for~\ref{enum:case:diminishing} and~\ref{enum:case:backtracking} respectively.
Indeed, for~\ref{enum:case:diminishing}, we have
$\gamma_{n} = 2(1-c) L_{\nabla \cost_{n}}^{-1} \geq 2\rho(1-c) L_{\nabla \cost_{n}}^{-1} \geq \beta L_{\nabla \cost_{n}}^{-1}$.
For~\ref{enum:case:backtracking}, the inequality~\eqref{eq:stepsize_lower_bound} implies
$\gamma_{n} \geq \min\{\gamma_{\rm init}, 2\rho(1-c)L_{\nabla \cost_{n}}^{-1}\} = \min\{\gamma_{\rm init}L_{\nabla \cost_{n}}, 2\rho(1-c)\}L_{\nabla \cost_{n}}^{-1} \geq \beta L_{\nabla \cost_{n}}^{-1}$,
where the last inequality follows from
$L_{\nabla \cost_{n}} = \varpi_{1}+\varpi_{2}\mu_{n}^{-1} \geq L_{\nabla \cost_{1}}$
by
$\mu_{n} \leq \mu_{1}$
(see~\eqref{eq:nonsummable}).
Hence, for both~\ref{enum:case:diminishing} and~\ref{enum:case:backtracking}, we have
\begin{equation}
  \thickmuskip=0.3\thickmuskip
  \medmuskip=0.3\medmuskip
  \thinmuskip=0.3\thinmuskip
  \arraycolsep=0.3\arraycolsep
  \gamma_{n}
  \geq \beta L_{\nabla \cost_{n}}^{-1}
  = \frac{\mu_{n}\beta}{\varpi_{1}\mu_{n} + \varpi_{2}}
  =\frac{\eta\mu_{n}\beta}{\varpi_{1}\eta\mu_{n} + \eta\varpi_{2}}
  \geq \frac{\eta\beta}{\varpi_{1}+\eta\varpi_{2}}\mu_{n},
\end{equation}
where the last inequality follows by
$\eta\mu_{n} \in (0,2^{-1}]$.

Then, the LHS in~\eqref{eq:extension:sum_gradient_bound} is bounded below as
\begin{align}
   & \hspace{-1em} c\sum_{n=\underline{k}}^{\overline{k}} \gamma_{n}
  \left(\mathcal{M}_{\gamma_{n}}^{\cost_{n},\phi}(\bm{x}_{n})\right)^{2}
  \geq \frac{c\eta\beta}{\varpi_{1} + \eta\varpi_{2}}\sum_{n=\underline{k}}^{\overline{k}} \mu_{n}
  \left(\mathcal{M}_{\gamma_{n}}^{\cost_{n},\phi}(\bm{x}_{n})\right)^{2}\hspace{-0.5em}                            \\
   & \hspace{-1em} \geq \frac{c\eta\beta}{\varpi_{1} + \eta\varpi_{2}}\min_{\underline{k}\leq n \leq \overline{k}}
  \left(\mathcal{M}_{\gamma_{n}}^{\cost_{n},\phi}(\bm{x}_{n})\right)^{2}
  \sum_{n=\underline{k}}^{\overline{k}} \mu_{n}                                                                    \\
   & \hspace{-1em} \geq \frac{c\eta\beta}{\varpi_{1} + \eta\varpi_{2}}\min_{\underline{k}\leq n \leq \overline{k}}
  \left(\mathcal{M}_{\widebar{\gamma}}^{\cost_{n},\phi}(\bm{x}_{n})\right)^{2}
  \sum_{n=\underline{k}}^{\overline{k}} \mu_{n}, \label{eq:tmp_convergence}
\end{align}
where the last inequality holds from~\eqref{eq:measure_smooth} with
{ \small
    \thickmuskip=0.0\thickmuskip
    \medmuskip=0.0\medmuskip
    \thinmuskip=0.0\thinmuskip
    \arraycolsep=0.0\arraycolsep
$\stdnorm{\frac{\bm{x}_{n} - \prox{\gamma_{n}\phi}(\bm{x}_{n}-\gamma_{n}\nabla \cost_{n}(\bm{x}_{n}))}{\gamma_{n}}} \geq \stdnorm{\frac{\bm{x}_{n} - \prox{\widebar{\gamma}\phi}(\bm{x}_{n}-\widebar{\gamma}\nabla \cost_{n}(\bm{x}_{n}))}{\widebar{\gamma}}}$}
with
$\widebar{\gamma} \geq \gamma_{n}\ (n\in \mathbb{N})$
by~\cite[Thm. 10.9]{Beck17} (or~\cite[Lemma 6]{Liu-Xia24}).
By~\eqref{eq:extension:sum_gradient_bound} with~\eqref{eq:tmp_convergence} and
$\beta =\min\{\gamma_{\rm init}L_{\nabla \cost_{1}}, 2\rho(1-c)\}$,
we get the desired inequality~\eqref{eq:convergence_rate}.

(b)
Let
$a_{n}\coloneqq \mathcal{M}_{\widebar{\gamma}}^{\cost_{n},\phi}(\bm{x}_{n})$,
and
$b_{n} \coloneqq  \inf\{a_{m} \mid m \geq n\}$.
For every
$\underline{k} \in \mathbb{N}$,
the inequality
\eqref{eq:convergence_rate} with
$\overline{k}\to\infty$
implies
$b_{\underline{k}} = 0$
by
$\sum_{n=1}^{\infty} \mu_{n}=+\infty$
(see the condition (ii) in~\eqref{eq:nonsummable}).
This yields
$\liminf_{n\to\infty}\mathcal{M}_{\widebar{\gamma}}^{\cost_{n},\phi}(\bm{x}_{n}) = \liminf_{n\to\infty} a_{n} = \sup_{n\in\mathbb{N}} b_{n} = 0$.

(c)
From the definition of
$m$,
$\lim_{l\to\infty} \mathcal{M}_{\widebar{\gamma}}^{\cost_{m(l)},\phi}(\bm{x}_{m(l)}) = 0$
holds.
Together with
Proposition~\ref{proposition:measure},
every cluster point
$\widebar{\bm{x}}$
of
$(\bm{x}_{m(l)})_{l=1}^{\infty}$
is a stationary point of
$\cost+\phi$.

(d)
By letting
$\gamma\coloneqq \gamma_{n}=2(1-c)L_{\nabla \cost_{n}}^{-1}$,
$\widebar{\bm{x}} \coloneqq \bm{x}_{n}$,
and
$\widebar{\bm{p}} \coloneqq \bm{x}_{n+1} = \prox{\gamma_{n} \phi}(\bm{x}_{n} - \gamma_{n} \nabla \cost_{n}(\bm{x}_{n}))$,
the inequalities~\eqref{eq:approximate_stationarity_inequality} and~\eqref{eq:bound_prox} in Proposition~\ref{proposition:stationary_relation} yield, for each
$n\in\mathbb{N}$,
{
    \thickmuskip=0.0\thickmuskip
    \medmuskip=0.0\medmuskip
    \thinmuskip=0.0\thinmuskip
    \arraycolsep=0.0\arraycolsep
    \begin{align}
       & \dist(\bm{0}, \nabla h(\bm{x}_{n+1}) + (\mathrm{D}\mathfrak{S}(\bm{x}_{n+1}))^{*}[\Lsubdiff g(\prox{\mu_{n} g}(\mathfrak{S}(\bm{x}_{n+1})))]\hspace{-1em} \\
       & \hspace{9em}+ \subdiff \phi(\bm{x}_{n+1}))
      \leq (3-2c) \mathcal{M}_{\gamma_{n}}^{\cost_{n}, \phi}(\bm{x}_{n}), \label{eq:approximate_stationarity_inequality_analysis}                     \\
       & \stdnorm{\prox{\mu_{n} g}(\mathfrak{S}(\bm{x}_{n+1})) - \mathfrak{S}(\bm{x}_{n+1})} \leq \mu_{n} L_{g}. \label{eq:bound_prox_analysis}
    \end{align}}%
Hence,
$\bm{x}_{\widehat{n}+1}$
is guaranteed to be an $\epsilon$-stationary point if the following inequalities are jointly satisfied with some
$\widehat{n}\in \mathbb{N}$
\begin{equation}
  (3-2c)\mathcal{M}_{\gamma_{\widehat{n}}}^{\cost_{\widehat{n}},\phi}(\bm{x}_{\widehat{n}}) \leq \epsilon,
  \ \mathrm{and}\
  \mu_{\widehat{n}}L_{g} \leq \epsilon. \label{eq:goal_rate}
\end{equation}

By
$\mu_{n} = \tau n^{-1/\alpha}$
with
$\alpha > 1$,
we have
$\mu_{n}L_{g} \leq \epsilon$
for all
$n \geq \underline{k} \coloneqq \lceil (\tau L_{g}\epsilon^{-1})^{\alpha} \rceil$.
We derive
$\overline{k} \in \mathbb{N}$
below such that
$\overline{k} \geq \underline{k}$
and
$(3-2c)\min_{\underline{k}\leq n \leq \overline{k}} \mathcal{M}_{\gamma_{n}}^{\cost_{n},\phi}(\bm{x}_{n}) \leq \epsilon$.
By~\eqref{eq:convergence_rate}
and
$\Delta =  \chi\tau^{-1}(1-\alpha^{-1})(3-2c)^{2} > 0$,
we have
\begin{align}
   (3-2c)\min_{\underline{k}\leq n \leq \overline{k}} \mathcal{M}_{\gamma_{n}}^{\cost_{n},\phi}(\bm{x}_{n})
  &  \leq \sqrt{\frac{(3-2c)^{2}\chi}{\tau\sum_{n=\underline{k}}^{\overline{k}}n^{-1/\alpha}}} \\
  & 
   \leq
  \sqrt{\frac{\Delta}{(\overline{k}+1)^{1-\alpha^{-1}} - \underline{k}^{1-\alpha^{-1}}}},
\end{align}
where we used
$\sum_{n=\underline{k}}^{\overline{k}} n^{-1/\alpha} = \sum_{n=\underline{k}}^{\overline{k}} \int_{n}^{n+1}n^{-1/\alpha}dx \geq \int_{\underline{k}}^{\overline{k}+1} x^{-1/\alpha} dx =  \frac{(\overline{k}+1)^{1-\alpha^{-1}}-\underline{k}^{1-\alpha^{-1}}}{1-\alpha^{-1}}$.
Hence,
$(3-2c)\min_{\underline{k}\leq n \leq \overline{k}} \mathcal{M}_{\gamma_{n}}^{\cost_{n},\phi}(\bm{x}_{n}) \leq \epsilon$
holds if
$\overline{k}$
satisfies
\begin{equation}
  \overline{k}
  \geq \left(\Delta\epsilon^{-2}  + \underline{k}^{1-\alpha^{-1}} \right)^{\frac{\alpha}{\alpha-1}} - 1.
  \label{eq:lower_iteration}
\end{equation}
In the following, we derive an upper bound, with
$\epsilon$,
of the RHS in~\eqref{eq:lower_iteration}.
By the subadditivity
$(a+b)^{p} \leq a^{p} + b^{p}\ (p \in (0,1), a,b\in \mathbb{R}_{++})$, e.g.~\cite[Footnote~12]{Kume-Yamada24A},
we get
$\underline{k}^{1-\alpha^{-1}}
  = \lceil (\tau L_{g} \epsilon^{-1})^{\alpha} \rceil^{1-\alpha^{-1}}
  \leq ((\tau L_{g} \epsilon^{-1})^{\alpha}+1)^{1-\alpha^{-1}}
  \leq (\tau L_{g} \epsilon^{-1})^{\alpha-1} + 1
  = (\tau L_{g})^{\alpha-1}\epsilon^{1-\alpha} + 1$.
From
$\max\{\epsilon^{-2}, \epsilon^{1-\alpha}, 1\} \leq \epsilon^{\min\{-2,1-\alpha\}}$
by
$\epsilon \in (0,1)$
and
$\alpha > 1$,
we get
$\underline{k}^{1-\alpha^{-1}} \leq ((\tau L_{g})^{\alpha-1} + 1)\epsilon^{1-\alpha}$.
Then, an upper bound of the RHS of~\eqref{eq:lower_iteration} is derived as
       $\left(\Delta\epsilon^{-2}  + \underline{k}^{1-\alpha^{-1}}\right)^{\frac{\alpha}{\alpha-1}} \leq \left(\left(\Delta + (\tau L_{g})^{\alpha-1}+1\right)\epsilon^{\min\{-2,1-\alpha\}}\right)^{\frac{\alpha}{\alpha-1}} \leq \left\lceil\left(\Delta + (\tau L_{g})^{\alpha-1} + 1\right)^{\frac{\alpha}{\alpha-1}} \epsilon^{\min\left\{\frac{-2\alpha}{\alpha-1},-\alpha\right\}}\right\rceil\eqqcolon n_{0}.$ 
Hence,~\eqref{eq:lower_iteration} is achieved by letting
$\overline{k}\coloneqq n_{0}$.
This implies
$(3-2c)\min_{\underline{k}\leq n \leq n_{0}} \mathcal{M}_{\gamma_{n}}^{\cost_{n},\phi}(\bm{x}_{n}) \leq \epsilon$,
i.e.,
there exists
$\widehat{n} \in [\underline{k},n_{0}]\cap \mathbb{N}$
such that
$(3-2c)\mathcal{M}_{\gamma_{\widehat{n}}}^{\cost_{\widehat{n}},\phi}(\bm{x}_{\widehat{n}}) \leq \epsilon$.
Since
$\mu_{n}L_{g} \leq \epsilon$
holds for all
$n \geq \underline{k}$,
we get
$\mu_{\widehat{n}}L_{g}\leq \epsilon$.
Therefore, the condition~\eqref{eq:goal_rate}
is achieved within at most
$n_{0}$
iterations.

\end{document}


\title{\Large Supplementary Material for \\ ``A Proximal Variable Smoothing for Minimization of Nonlinearly Composite Nonsmooth Function   -- Finite-Max Minimization and MIMO Applications''}
\author{Keita Kume \IEEEmembership{Member, IEEE}, Isao Yamada \IEEEmembership{Fellow, IEEE}}

\date{}

\maketitle
\setcounter{section}{0}
\setcounter{equation}{0}
\setcounter{figure}{0}
\setcounter{table}{0}

\renewcommand{\thesection}{S\arabic{section}}
\renewcommand{\thesubsection}{\thesection.\arabic{subsection}}

\renewcommand{\theequation}{S.\arabic{equation}}
\renewcommand{\thefigure}{S.\arabic{figure}}
\renewcommand{\thetable}{S.\arabic{table}}
In this supplementary material, we present complete proofs of results in application parts in Section~\ref{sec:application} of the main paper.
\section{Proof of Lemma~\ref{lemma:application_constant}}
For each
  $j \in \{1,2,\ldots, m\}$,
  we get expressions of the gradient and the Hessian matrix of
  $\mathfrak{S}_{j}(\bm{x}) \coloneqq (y_{j}^{2}-\stdnorm{\bm{x}-\bm{u}_{j}}^{2})^{2}$
  with
  $\bm{q}_{j}(\bm{x}) \coloneqq \bm{x}-\bm{u}_{j}$
  as
  \begin{align}
    (\bm{x} \in \mathcal{X}) & \quad \nabla \mathfrak{S}_{j}(\bm{x}) = -4(y_{j}^{2} - \stdnorm{\bm{q}_{j}(\bm{x})}^{2})\bm{q}_{j}(\bm{x}), \label{eq:S_j_grad} \\
    (\bm{x} \in \mathcal{X}) & \quad \nabla^{2} \mathfrak{S}_{j}(\bm{x}) = 8\bm{q}_{j}(\bm{x})\bm{q}_{j}(\bm{x})^{\TT} + 4(\stdnorm{\bm{q}_{j}(\bm{x})}^{2}-y_{j}^{2})\Id. \hspace{5em} \label{eq:S_j_Hesse}
  \end{align}
  For
  $\bm{x} \in C \coloneqq B(\bm{0}, r)\coloneqq \{\bm{x}\in\mathcal{X}\mid \stdnorm{\bm{x}}\leq r\}$,
  we have a bound
  \begin{equation}
    q_{j,{\rm low}}\coloneqq \max\{0, \stdnorm{\bm{u}_{j}} - r\} \leq \stdnorm{\bm{q}_{j}(\bm{x})} \leq \stdnorm{\bm{u}_{j}} + r\eqqcolon q_{j,{\rm up}}. \label{eq:q_bound}
  \end{equation}
  By using the above expressions, we prove Lemma~\ref{lemma:application_constant}.

  (For
  $L_{\nabla\mathfrak{S}_{j}}$)
  From the twice continuous differentiability of
  $\mathfrak{S}_{j}$,
  we can take an upper bound of
  $\max_{\bm{x}\in C}\opnorm{\nabla^{2}\mathfrak{S}_{j}(\bm{x})}$
  as a Lipschitz constant
  $L_{\nabla\mathfrak{S}_{j}}$
  of
  $\nabla \mathfrak{S}_{j}$
  over
  $C$
  by~\cite[Thm. 9.7]{Rockafellar-Wets98}.
  Consider the case
  $\bm{q}_{j}(\bm{x})\neq \bm{0}$.
  Since we have
  {
      \thickmuskip=0.1\thickmuskip
      \medmuskip=0.1\medmuskip
      \thinmuskip=0.1\thinmuskip
      \arraycolsep=0.1\arraycolsep
  \begin{align}
    \nabla^{2} \mathfrak{S}_{j}(\bm{x})[\bm{q}_{j}(\bm{x})] 
    & = 8\stdnorm{\bm{q}_{j}(\bm{x})}^{2}\bm{q}_{j}(\bm{x}) + 4(\stdnorm{\bm{q}_{j}(\bm{x})}^{2}-y_{j}^{2})\bm{q}_{j}(\bm{x}) \hspace{-1em}\\
    & =4(3\stdnorm{\bm{q}_{j}(\bm{x})}^{2} - y_{j}^{2})\bm{q}_{j}(\bm{x}), \label{eq:supp:eigen_1}
\end{align}}%
  $\bm{q}_{j}(\bm{x})$
  is an eigenvector of
  $\nabla^{2} \mathfrak{S}_{j}(\bm{x})$
  with eigenvalue
  $4(3\stdnorm{\bm{q}_{j}(\bm{x})}^{2} - y_{j}^{2})$.
  On the other hand, since we have
  \begin{equation}
      \thickmuskip=0.2\thickmuskip
      \medmuskip=0.2\medmuskip
      \thinmuskip=0.2\thinmuskip
      \arraycolsep=0.2\arraycolsep
    (\bm{d}\in \mathcal{X}\ \mathrm{s.t.}\ \inprod{\bm{d}}{\bm{q}_{j}(\bm{x})}=0) \ 
    \nabla^{2} \mathfrak{S}_{j}(\bm{x})[\bm{d}] 
    = 4(\stdnorm{\bm{q}_{j}(\bm{x})}^{2}-y_{j}^{2})\bm{d}, \label{eq:supp:eigen_2}
\end{equation}
  every
  $\bm{d} \in \mathcal{X}\setminus\{\bm{0}\}$
  such that
  $\inprod{\bm{d}}{\bm{q}_{j}(\bm{x})} = 0$
  is also an eigenvector of
  $\nabla^{2} \mathfrak{S}_{j}(\bm{x})$
  with eigenvalue
  $4(\stdnorm{\bm{q}_{j}(\bm{x})}^{2}-y_{j}^{2})$.
  By combining them,
  with
  $\Upsilon_{j}^{(1)}(q)\coloneqq \max\{\abs{3q^{2} - y_{j}^{2}}, \abs{q^{2}-y_{j}^{2}}\}$,
  we have
  $\opnorm{\nabla^{2}\mathfrak{S}_{j}(\bm{x})} \leq 4\Upsilon_{j}^{(1)}(\stdnorm{\bm{q}_{j}(\bm{x})}) \leq 
  4\max_{q_{j,{\rm low}} \leq q \leq q_{j,{\rm up}}}\Upsilon_{j}^{(1)}(q) =
  4\max\{\Upsilon_{j}^{(1)}(q_{j,{\rm low}}), \Upsilon_{j}^{(1)}(q_{j,{\rm up}})\}$
  by the bound~\eqref{eq:q_bound},
  where the last equality follows from the fact that
  (i)
  $\abs{3q^{2} - y_{j}^{2}}$
  is monotonically decreasing on
  $[0,\abs{y_{j}}/\sqrt{3}]$
  and is monotonically increasing on
  $[\abs{y_{j}}/\sqrt{3},+\infty)$;
  (ii)
  $\abs{q^{2}-y_{j}^{2}}$
  is monotonically decreasing on
  $[0,\abs{y_{j}}]$
  and is monotonically increasing on
  $[\abs{y_{j}},+\infty)$.
  Consider the other case
  $\bm{q}_{j}(\bm{x}) = \bm{0}$,
  hence
  $q_{j,{\rm low}} = 0$.
  Then, we have
  $\nabla^{2}\mathfrak{S}_{j}(\bm{x}) = -4y_{j}^{2}\Id$,
  i.e.,
  $\opnorm{\nabla^{2}\mathfrak{S}_{j}(\bm{x})} =4y_{j}^{2} = 4\Upsilon_{j}^{(1)}(0)=4\Upsilon_{j}^{(1)}(q_{j,{\rm low}})$.
  Consequently,
  we get
  $\max_{\bm{x}\in C}\opnorm{\nabla^{2}\mathfrak{S}_{j}(\bm{x})} \leq 4\max\{\Upsilon_{j}^{(1)}(q_{j,{\rm low}}), \Upsilon_{j}^{(1)}(q_{j,{\rm up}})\}$,
  i.e.,
$4\max\{\Upsilon_{j}^{(1)}(q_{j,{\rm low}}), \Upsilon_{j}^{(1)}(q_{j,{\rm up}})\}$
is a Lipschitz constant
$L_{\nabla\mathfrak{S}_{j}}$ 
of
$\nabla \mathfrak{S}_{j}$
over
$C$.

  (For $\widetilde{\eta}$)
  By~\eqref{eq:supp:eigen_1} and~\eqref{eq:supp:eigen_2} and by~\eqref{eq:S_j_Hesse} with
  the case
  $\bm{q}_{j}(\bm{x}) = \bm{0}$,
  the smallest eigenvalue of
  $\nabla^{2}\mathfrak{S}_{j}(\bm{x})$
  for
  $\bm{x} \in C$
  is bounded below by
  $\min\{4(3\stdnorm{\bm{q}_{j}(\bm{x})}^{2}- y_{j}^{2}), 4(\stdnorm{\bm{q}_{j}(\bm{x})}^{2}-y_{j}^{2}), -4y_{j}^{2}\} = -4y_{j}^{2}$.
  Then,
  $\mathfrak{S}_{j} + \frac{4y_{j}^{2}}{2}\stdnorm{\cdot}^{2}$
  is convex~\cite[Thm. 2.14]{Rockafellar-Wets98} because its Hessian matrix
  $\nabla^{2}(\mathfrak{S}_{j}+\frac{4y_{j}^{2}}{2}\stdnorm{\cdot}^{2})(\bm{x}) = \nabla^{2}\mathfrak{S}_{j}(\bm{x})+4y_{j}^{2}\Id$
  is positive semidefinite.
  Hence, for every
  $j = 1,2,\ldots,m$,
  $\mathfrak{S}_{j}+\frac{4\max_{k=1,2,\ldots,m}y_{k}^{2}}{2}\stdnorm{\cdot}^{2}$
  is also convex.
  Since the maximum of convex functions is convex~\cite[Prop. 8.16]{Bauschke-Combettes17},
  the function
  $\max_{j=1,2,\ldots,m}(\mathfrak{S}_{j}+\frac{4\max_{k=1,2,\ldots,m}y_{k}^{2}}{2}\stdnorm{\cdot}^{2})(\bm{x})=\cost(\bm{x}) + \frac{4\max_{k=1,2,\ldots,m}y_{k}^{2}}{2}\stdnorm{\bm{x}}^{2}$
  is convex with
  $\cost(\bm{x}) = \max_{j=1,2,\ldots,m}\mathfrak{S}_{j}(\bm{x})$
  from which
  $\cost$
  is
  $\widetilde{\eta}$-weakly convex with
  $\widetilde{\eta}\coloneqq (4\max_{k=1,2,\ldots,m}y_{k}^{2})$.

  (For
  $\kappa_{\mathfrak{S}_{j}}$)
  By~\eqref{eq:S_j_grad}, we get
  $\stdnorm{\nabla \mathfrak{S}_{j}(\bm{x})} \leq 4\stdnorm{\bm{q}_{j}(\bm{x})}\times\abs{y_{j}^{2} - \stdnorm{\bm{q}_{j}(\bm{x})}^{2}} = 4\Upsilon_{j}^{(2)}(\stdnorm{\bm{q}_{j}(\bm{x})})$
  with
  $\Upsilon_{j}^{(2)}(q)\coloneqq q\abs{y_{j}^{2} - q^{2}}$.
  By using the bound~\eqref{eq:q_bound}, we obtain
  $\kappa_{\mathfrak{S}_{j}} \leq 4\max_{\bm{x}\in C}\Upsilon_{j}^{(2)}(\stdnorm{\bm{q}_{j}(\bm{x})}) \leq 4\max_{q_{j,{\rm low}} \leq q \leq q_{j,{\rm up}}} \Upsilon_{j}^{(2)}(q)$.

The function
$\Upsilon_{j}^{(2)}$
is maximized over
$[q_{j,{\rm low}},q_{j,{\rm up}}]\subset[0,+\infty)$
at one of
$q\in Q_{j}$,
where
  $Q_{j} \coloneqq \{ \abs{y_{j}}/\sqrt{3},\ q_{j,{\rm up}}\}$
  if
  $\abs{y_{j}}/\sqrt{3} \in [q_{j,{\rm low}}, q_{j,{\rm up}}]$;
  $Q_{j} \coloneqq \{q_{j,{\rm low}},\ q_{j,{\rm up}}\}$
  otherwise.
Indeed,
as illustrated in Fig.~S.1,
$\Upsilon_{j}^{(2)}$
is 
(i)
monotonically increasing on
$[0,\abs{y_{j}}/\sqrt{3}]$;
(ii)
monotonically decreasing on
$(\abs{y_{j}}/\sqrt{3},\abs{y_{j}}]$;
(iii)
monotonically increasing on
$(\abs{y_{j}},+\infty)$.
Hence,
we get
  $\kappa_{\mathfrak{S}_{j}} \leq 4\max_{q \in Q_{j}} \Upsilon_{j}^{(2)}(q)$.

  \begin{figure}[t]
    \centering
\begin{tikzpicture}
\begin{axis}[
    width=7cm,
    height=4cm,
    axis lines=middle,
    xlabel={$q$},
    ylabel={$q|y_j^2-q^2|$},
    xmin=-0.2,
    xmax=1.35,
    ymin=-0.3,
    ymax=1,
    samples=300,
    domain=-0.2:1.35,
    smooth,
    thick,
    xtick={0,0.5773502692,1},
    xticklabels={
        {$0$},
        {$\abs{y_{j}}/\sqrt{3}$},
        {$\abs{y_j}$}
    },
    ytick=\empty,
    tick label style={font=\small},
    label style={font=\small},
    clip=false,
]

\addplot[blue, very thick]
    {x*abs(1 - x^2)};

\addplot[
    only marks,
    mark=*,
    mark size=2pt,
    blue
]
coordinates {
    (0,0)
    (0.5773502692,0.3849001795)
    (1,0)
};

\addplot[dashed, gray] coordinates {(0.5773502692,0) (0.5773502692,0.3849001795)};
\addplot[dashed, gray] coordinates {(1,0) (1,0.05)};

\node[above right, font=\small] at (axis cs:0.5773502692,0.3849001795)
{$\dfrac{2}{3\sqrt{3}}\abs{y_j}^3$};

\end{axis}
\end{tikzpicture}
\label{fig:func_plot1}
\caption{Illustration of the graph of $\Upsilon_{j}^{(2)}(q)\coloneq q\abs{y_{j}^{2}-q^{2}}$.}
\end{figure}

\section{Proof of Lemma~\ref{lemma:application_mimo_constant}}
(For $\phi$)
Since
$C=[\underline{r},1]^{U}\times \mathbb{R}^{U}$
is closed convex,
$\phi=\iota_{C}$
is proper lower semicontinuous convex.
Clearly,
$\phi$
is prox-friendly since
$\prox{\gamma \phi}(\bm{r},\bm{\theta})=P_{C}(\bm{r},\bm{\theta})$
can be expressed as
$(P_{[\underline{r},1]^{U}}(\bm{r}),\bm{\theta})$
with the projection
$P_{[\underline{r},1]^{U}}$
onto the box constraint
$[\underline{r},1]^{U}$.

(For
$h$)
Recall that
$h(\bm{r},\bm{\theta}) = q_{\bm{y},\bm{H}}\circ \param(\bm{r},\bm{\theta}) + \lambda_{r}\sum_{u=1}^{U} \mathfrak{d}([\bm{r}]_{u})$
is given by the continuously differentiable mappings
$q_{\bm{y},\bm{H}}:\mathbb{R}^{2U}\to\mathbb{R}:\bm{s}\mapsto \frac{1}{2}\stdnorm{\bm{y}-\bm{H}\bm{s}}^{2}$
and
$\param:\mathbb{R}^{2U}\to\mathbb{R}^{2U}:(\bm{r},\bm{\theta}) \mapsto \begin{bmatrix}
    \bm{r}\odot \mathrm{\bf cos}(\bm{\theta}) \\ \bm{r} \odot \mathrm{\bf sin}(\bm{\theta})
  \end{bmatrix}$,
  where
$\odot$
is the entry-wise (Hadamard) product,
and
$\mathfrak{d}:\mathbb{R}\to\mathbb{R}$
is a continuously differentiable mapping satisfying
$\mathfrak{d}(t) = t^{-1}$
for
$t\in [\underline{r},1]$.
Hence,
$h$
is also continuously differentiable.
For the Lipschitz continuity of
$\nabla h$
over
$C$,
it is enough to show the Lipschitz continuity of
$\nabla (q_{\bm{y},\bm{H}}\circ\param)$
over
$C$
because
$\nabla \mathfrak{d}(t)=-\frac{1}{t^{2}}\ (t\in [\underline{r},1])$
is clearly Lipschitz continuous over
$[\underline{r},1] \subset (0,1]$.
By, e.g.,~\cite[Lemma C.1]{Kume-Yamada24A},
$\nabla (q_{\bm{y},\bm{H}}\circ\param)$
is Lipschitz continuous over
$C$
if the following four conditions are jointly satisfied:
\begin{enumerate}[label=(\alph*)]
  \item
$\sup_{(\bm{r},\bm{\theta})\in C}\opnorm{\mathrm{D}\param(\bm{r},\bm{\theta})} < +\infty$;
\item
$\mathrm{D}\param$
is Lipschitz continuous over
$C$;
\item
$\sup_{\bm{s} \in \param(C)}\stdnorm{\nabla q_{\bm{y},\bm{H}}(\bm{s})} < + \infty$;
\item
$\nabla q_{\bm{y},\bm{H}}$
is Lipschitz continuous over
$\param(C)$,
\end{enumerate}
where
the G\^{a}teaux derivative of
$\param$
is defined by
\begin{equation}
  \mathrm{D}\param(\bm{r},\bm{\theta})[\bm{d}_{r},\bm{d}_{\theta}] \coloneqq \lim_{\mathbb{R}\setminus\{0\} \ni t \to 0} \frac{\param(\bm{r}+t\bm{d}_{r},\bm{\theta}+t\bm{d}_{\theta}) - \param(\bm{r},\bm{\theta})}{t},
\end{equation}
and we have the expressions:
\begin{align}
  & ((\bm{r},\bm{\theta}),(\bm{d}_{r},\bm{d}_{\theta}) \in \mathcal{X})  \\
  & \quad
  \mathrm{D}\param(\bm{r},\bm{\theta})[\bm{d}_{r},\bm{d}_{\theta}] = \begin{bmatrix} \bm{d}_{r}\odot \mathrm{\bf cos}(\bm{\theta}) - \bm{r}\odot \mathrm{\bf sin}(\bm{\theta})\odot \bm{d}_{\theta} \\ \bm{d}_{r} \odot \mathrm{\bf sin}(\bm{\theta}) + \bm{r}\odot \mathrm{\bf cos}(\bm{\theta}) \odot \bm{d}_{\theta} \end{bmatrix}, \\
  & (\bm{s}\in \mathbb{R}^{2U})  \quad \nabla q_{\bm{y},\bm{H}}(\bm{s})  = \bm{H}^{*}(\bm{H}\bm{s}-\bm{y}).
\end{align}
In what follows, we verify (a)-(d).

(a)
From the definition of the operator norm and the triangle inequality,
we have
for any
$(\bm{r},\bm{\theta}) \in C=[\underline{r},1]^{U}\times \mathbb{R}^{U}$,
{
      \thickmuskip=0.0\thickmuskip
      \medmuskip=0.0\medmuskip
      \thinmuskip=0.0\thinmuskip
      \arraycolsep=0.0\arraycolsep
\begin{align}
& \opnorm{\mathrm{D}\param(\bm{r},\bm{\theta})}
= \sup_{\stdnorm{(\bm{d}_{r},\bm{d}_{\theta})}
\leq 1}\stdnorm{\mathrm{D}\param(\bm{r},\bm{\theta})[\bm{d}_{r},\bm{d}_{\theta}]} \\
& \leq 
\sup_{\stdnorm{(\bm{d}_{r},\bm{d}_{\theta})}
\leq 1} \sum_{j=1}^{U} (\abs{[\bm{d}_{r}]_{j}\cos([\bm{\theta}]_{j}) - [\bm{r}]_{j}\sin([\bm{\theta}]_{j})[\bm{d}_{\theta}]_{j}}\\
& \quad\quad + \abs{[\bm{d}_{r}]_{j}\sin([\bm{\theta}]_{j}) + [\bm{r}]_{j}\cos([\bm{\theta}]_{j})[\bm{d}_{\theta}]_{j}}) \\
& \leq 
\sup_{\stdnorm{(\bm{d}_{r},\bm{d}_{\theta})}
\leq 1} \sum_{j=1}^{U} (\abs{[\bm{d}_{r}]_{j}\cos([\bm{\theta}]_{j})} + \abs{[\bm{r}]_{j}\sin([\bm{\theta}]_{j})[\bm{d}_{\theta}]_{j}}\\
& \quad\quad + \abs{[\bm{d}_{r}]_{j}\sin([\bm{\theta}]_{j})} + \abs{[\bm{r}]_{j}\cos([\bm{\theta}]_{j})[\bm{d}_{\theta}]_{j}}) \\
& \leq\hspace{-1em}
\sup_{(\bm{d}_{r},\bm{d}_{\theta})\in [-1,1]^{U}\times [-1,1]^{U}} \sum_{j=1}^{U} (\abs{[\bm{d}_{r}]_{j}\cos([\bm{\theta}]_{j})} + \abs{[\bm{r}]_{j}\sin([\bm{\theta}]_{j})[\bm{d}_{\theta}]_{j}} \\
& \quad\quad + \abs{[\bm{d}_{r}]_{j}\sin([\bm{\theta}]_{j})} + \abs{[\bm{r}]_{j}\cos([\bm{\theta}]_{j})[\bm{d}_{\theta}]_{j}}) \\
& \leq
\sum_{j=1}^{U} \left(\abs{\cos([\bm{\theta}]_{j})} + \abs{[\bm{r}]_{j}\sin([\bm{\theta}]_{j})} + \abs{\sin([\bm{\theta}]_{j})} + \abs{[\bm{r}]_{j}\cos([\bm{\theta}]_{j})}\right). \label{eq:op_norm_max}
\end{align}}%
From
$[\bm{r}]_{j} \in [\underline{r},1]$,
$\abs{\cos([\bm{\theta}]_{j})} \leq 1$
and
$\abs{\sin([\bm{\theta}]_{j})} \leq 1$,
we get
$\sup_{(\bm{r},\bm{\theta}) \in C} \opnorm{\mathrm{D}\param(\bm{r},\bm{\theta})} \leq \sum_{j=1}^{U} (1+1+1+1) = 4U$.

(b)
Let
$(\bm{r}_{1},\bm{\theta}_{1}), (\bm{r}_{2},\bm{\theta}_{2}) \in C=[\underline{r},1]^{U}\times \mathbb{R}^{U}$.
By the discussion similar to~\eqref{eq:op_norm_max}, we get
{
      \thickmuskip=0.1\thickmuskip
      \medmuskip=0.1\medmuskip
      \thinmuskip=0.1\thinmuskip
      \arraycolsep=0.1\arraycolsep
\begin{align}
& \hspace{-2em}\opnorm{\mathrm{D}\param(\bm{r}_{1},\bm{\theta}_{1})- \mathrm{D}\param(\bm{r}_{2},\bm{\theta}_{2})}\\
& \hspace{-2em}\leq \sum_{j=1}^{U} (\abs{\cos([\bm{\theta}_{1}]_{j}) - \cos([\bm{\theta}_{2}]_j)} + \abs{[\bm{r}_{1}]_{j}\sin([\bm{\theta}_{1}]_{j})- [\bm{r}_{2}]_{j}\sin([\bm{\theta}_{2}]_{j})} \\
&\hspace{-2em} \quad + \abs{\sin([\bm{\theta}_{1}]_{j}) - \sin([\bm{\theta}_{2}]_{j})} + \abs{[\bm{r}_{1}]_{j}\cos([\bm{\theta}_{1}]_{j})- [\bm{r}_{2}]_{j}\cos([\bm{\theta}_{2}]_{j})}). \label{eq:op_norm_max_Lip}
\end{align}}%
Since
$\cos$
and
$\sin$
are $1$-Lipschitz continuous, we get
$\abs{\cos([\bm{\theta}_{1}]_{j}) - \cos([\bm{\theta}_{2}]_j)} \leq \abs{[\bm{\theta}_{1}]_{j}-[\bm{\theta}_{2}]_{j}} \leq \stdnorm{(\bm{r}_{1},\bm{\theta}_{1}) - (\bm{r}_{2},\bm{\theta}_{2})}$
and
$\abs{\sin([\bm{\theta}_{1}]_{j}) - \sin([\bm{\theta}_{2}]_j)} \leq \stdnorm{(\bm{r}_{1},\bm{\theta}_{1}) - (\bm{r}_{2},\bm{\theta}_{2})}$.
For showing the Lipschitz continuity of
$\mathrm{D}\param$
over
$C$,
it is enough to show the Lipschitz continuities of
$\abs{[\bm{r}_{1}]_{j}\sin([\bm{\theta}_{1}]_{j})- [\bm{r}_{2}]_{j}\sin([\bm{\theta}_{2}]_{j})}$
and
$\abs{[\bm{r}_{1}]_{j}\cos([\bm{\theta}_{1}]_{j})- [\bm{r}_{2}]_{j}\cos([\bm{\theta}_{2}]_{j})}$
over
$C$.

We verify only the former
(i.e., the Lipschitz continuities of
$\abs{[\bm{r}_{1}]_{j}\sin([\bm{\theta}_{1}]_{j})- [\bm{r}_{2}]_{j}\sin([\bm{\theta}_{2}]_{j})}$
over
$C$)
because the same discussion can apply to the latter.
The triangle inequality yields
\begin{align}
& \abs{[\bm{r}_{1}]_{j}\sin([\bm{\theta}_{1}]_{j})- [\bm{r}_{2}]_{j}\sin([\bm{\theta}_{2}]_{j})}\\
&\leq \abs{[\bm{r}_{1}]_{j}(\sin([\bm{\theta}_{1}]_{j}) - \sin([\bm{\theta}_{2}]_{j}))} + \abs{([\bm{r}_{1}]_{j}-[\bm{r}_{2}]_{j})\sin([\bm{\theta}_{2}]_{j})}. \label{eq:supp:Lip}
\end{align}
By
$[\bm{r}_{1}]_{j} \in [\underline{r},1]$,
the first term of RHS in~\eqref{eq:supp:Lip} is bounded above as
\begin{align}
& \abs{[\bm{r}_{1}]_{j}(\sin([\bm{\theta}_{1}]_{j}) - \sin([\bm{\theta}_{2}]_{j}))}
 \leq \abs{\sin([\bm{\theta}_{1}]_{j}) - \sin([\bm{\theta}_{2}]_{j})} \\
& \leq \abs{[\bm{\theta}_{1}]_{j}-[\bm{\theta}_{2}]_{j}} 
 \leq \stdnorm{(\bm{r}_{1},\bm{\theta}_{1}) - (\bm{r}_{2},\bm{\theta}_{2})}.
\end{align}
On the other hand,
by
$\abs{\sin([\bm{\theta}_{2}]_{j})} \leq 1$,
the second term of RHS in~\eqref{eq:supp:Lip} is bounded above as
\begin{align}
\abs{([\bm{r}_{1}]_{j}-[\bm{r}_{2}]_{j})\sin([\bm{\theta}_{2}]_{j})} 
& \leq \abs{[\bm{r}_{1}]_{j}-[\bm{r}_{2}]_{j}} \\
& \leq \stdnorm{(\bm{r}_{1},\bm{\theta}_{1}) - (\bm{r}_{2},\bm{\theta}_{2})}.
\end{align}

Therefore,
$\mathrm{D}\param$
is Lipschitz continuous over
$C$
because all terms of RHS in~\eqref{eq:op_norm_max_Lip} are Lipschitz continuous over
$C$.

(c)
Since we have
$\bm{r} \in [\underline{r},1]^{U}$,
$\mathrm{\bf cos}(\bm{\theta}) \in [-1,1]^{U}$
and
$\mathrm{\bf sin}(\bm{\theta}) \in [-1,1]^{U}$
for
$(\bm{r},\bm{\theta}) \in C$,
we have
$\param(C) \subset [-1, 1]^{2U}$.
Since
$\nabla q_{\bm{y},\bm{H}}$
is continuous over the compact set
$[-1, 1]^{2U}$,
we get
$\sup_{\bm{s}\in \param(C)}\stdnorm{\nabla q_{\bm{y},\bm{H}}(\bm{s})} < + \infty$.

(d)
Since
$\nabla q_{\bm{y},\bm{H}}(\bm{s}) = \bm{H}^{*}(\bm{H}\bm{s}-\bm{y})$
is
$\opnorm{\bm{H}^{*}\bm{H}}$-Lipschitz continuous over
$\mathbb{R}^{2U}$,
we get the desired Lipschitz continuity of
$\nabla q_{\bm{y},\bm{H}}$
over
$\param(C)$.

(For $\mathfrak{S}$)
For
$\mathfrak{S}(\bm{r},\bm{\theta}) = \mathrm{\bf sin}\left(\frac{M\bm{\theta}}{2}\right)$,
we have
$\mathrm{D}\mathfrak{S}(\bm{r},\bm{\theta})[\bm{d}_{r},\bm{d}_{\theta}]=\frac{M}{2}\mathrm{\bf cos}(\frac{M\bm{\theta}}{2})\odot \bm{d}_{\theta}$.
Hence, a similar discussion to~\eqref{eq:op_norm_max} or~\eqref{eq:op_norm_max_Lip}
yields, for
$(\bm{r}_{1},\bm{\theta}_{1}), (\bm{r}_{2},\bm{\theta}_{2}) \in C$,
$\opnorm{\mathrm{D}\mathfrak{S}(\bm{r}_{1},\bm{\theta}_{1})-\mathrm{D}\mathfrak{S}(\bm{r}_{2},\bm{\theta}_{2})} \leq \frac{M}{2}\sum_{j=1}^{U}\abs{\cos(\frac{M[\bm{\theta}_{1}]_{j}}{2})-\cos(\frac{M[\bm{\theta}_{2}]_{j}}{2})} \leq \frac{M}{2}\sum_{j=1}^{U}\abs{\frac{M}{2}[\bm{\theta}_{1}]_{j}-\frac{M}{2}[\bm{\theta}_{2}]_{j}} \leq \frac{M^{2}}{4}\sum_{j=1}^{U}\stdnorm{(\bm{r}_{1},\bm{\theta}_{1})-(\bm{r}_{2},\bm{\theta}_{2})}$,
  i.e.,
$\mathrm{D}\mathfrak{S}(\bm{r},\bm{\theta})$
is Lipschitz continuous over
$C$.
Moreover,
$\sup_{(\bm{r},\bm{\theta})\in C} \opnorm{\mathrm{D}\mathfrak{S}(\bm{r},\bm{\theta})} \leq \sum_{j=1}^{U} \abs{\frac{M}{2}\cos(\frac{M[\bm{\theta}]_{j}}{2})} \leq \frac{M}{2}U < +\infty$.

(For $g$)
$g=\lambda_{\theta}\norm{\cdot}_{1}$
is clearly convex, i.e.,
$\eta$-weakly convex for any
$\eta > 0$,
and Lipschitz continuous over
$\mathcal{Z}$.
The proximity operator of
$g$
is given by the so-called soft-thresholding~\cite[Exm. 24.22]{Bauschke-Combettes17}.

(For $\cost+\phi$)
Since each term in
$\cost$
and
$\phi$
takes nonnegative values on
$\dom{\phi}=C$,
we have
$\inf_{\bm{x}\in \mathcal{X}}(\cost+\phi)(\bm{x}) \geq 0$.

\bibliographystyle{IEEEtran}
\bibliography{main}